\documentclass{amsart}
\usepackage{accents}
\usepackage{graphicx}
\usepackage{amsfonts}
\usepackage{amscd}
\usepackage{amssymb}
\usepackage{stmaryrd}
\usepackage{xypic}
\xyoption{all}
\usepackage{hyperref}
\usepackage{url}
\usepackage{tikz}
\usetikzlibrary{shapes}
\usetikzlibrary{decorations.markings}
\usepackage{multirow}
\usepackage[savepos]{zref}
\newtheorem{thm}{Theorem}[section]
\newtheorem{corollary}[thm]{Corollary}
\newtheorem{lemma}[thm]{Lemma}
\newtheorem{proposition}[thm]{Proposition}
\newtheorem{construction}[thm]{Construction}

\newtheorem{conj}[thm]{Conjecture}

\theoremstyle{definition}
\newtheorem{definition}[thm]{Definition}
\theoremstyle{remark}
\newtheorem{remark}[thm]{Remark}
\numberwithin{equation}{section}
\DeclareMathOperator{\checkdot}{\check{\cdot}}
\DeclareMathOperator{\checkodot}{\check{\odot}}
\DeclareMathOperator{\Fr}{Fr}
\DeclareMathOperator{\Irr}{Irr}
\DeclareMathOperator{\Par}{Par}

\DeclareMathOperator{\Gal}{Gal}

\DeclareMathOperator{\Ker}{Ker}

\DeclareMathOperator{\Hilb}{Hilb}
\DeclareMathOperator{\Leg}{Leg}
\DeclareMathOperator{\Char}{char}
\DeclareMathOperator{\rec}{rec}

\DeclareMathOperator{\Spec}{Spec}

\DeclareMathOperator{\Id}{Id}
\DeclareMathOperator{\Res}{Res}
\DeclareMathOperator{\Inc}{Inc}
\DeclareMathOperator{\ev}{ev}
\newdir{ >}{{}*!/-5pt/@{>}}

\newcommand{\from}{\colon}
\newcommand{\LVee}{ {}^\vee}

\renewcommand{\Im}{\mathrm{Im}}

\newcommand{\dual}{\diamond}

\newcommand{\Weyl}{\mathrm{W}}
\newcommand{\Braid}{\mathrm{B}}
\newcommand{\Weil}{\mathcal{W}}
\newcommand{\weil}{\mathbf{w}}

\newcommand{\Inertia}{\mathcal{I}}

\newcommand{\liminv}{\varprojlim}
\newcommand{\limdir}{\varinjlim}

\newcommand{\nt}[1]{{\mathring{#1}}}
\newcommand{\rt}[1]{\vec{#1}}
\newcommand{\lt}[1]{\reflectbox{\ensuremath{\vec{\reflectbox{\ensuremath{#1}}}}}}

\newcommand{\C}[1]{{\mathsf{\widehat{#1}}}}

\newcommand{\One}{\mathbf{1}}

\newcommand{\Cat}[1]{ {\mathsf{#1}} }

\newcommand{\Lie}[1]{ {\mathfrak{#1}} }

\newcommand{\alg}[1]{\mathbf{#1}}

\newcommand{\Cont}{{\mathcal{C}}}
\newcommand{\Hom}{ {\mbox{Hom}} }

\newcommand{\ZZ}{\mathbb Z}
\newcommand{\QQ}{\mathbb Q}
\newcommand{\RR}{\mathbb R}
\newcommand{\NN}{\mathbb N}
\newcommand{\CC}{\mathbb C}
\newcommand{\PP}{\mathcal P}
\newcommand{\HH}{\mathcal H}
\newcommand{\BB}{\mathcal B}

\newcommand{\UU}{\mathcal U}

\newcommand{\ident}{\equiv}

\newcommand{\FF}{\mathbb F}
\newcommand{\OO}{\mathcal{O}}

\newcommand{\isom}{\cong}

\begin{document}

\title{Split Metaplectic groups and their L-groups}%
\author{Martin H. Weissman}%
\date{\today}

\address{Dept. of Mathematics, University of California, Santa Cruz, CA 95064}
\email{weissman@ucsc.edu}%


\begin{abstract}
We adapt the conjectural local Langlands parameterization to split metaplectic groups over local fields.  When $\tilde G$ is a central extension of a split connected reductive group over a local field (arising from the framework of Brylinski and Deligne), we construct a dual group $\alg{\tilde G}\LVee$ and an L-group ${}^L \alg{\tilde G}\LVee$ as group schemes over $\ZZ$.  Such a construction leads to a definition of Weil-Deligne parameters (Langlands parameters) with values in this L-group, and to a conjectural parameterization of the irreducible genuine representations of $\tilde G$.  This conjectural parameterization is compatible with what is known about metaplectic tori, Iwahori-Hecke algebra isomorphisms between metaplectic and linear groups, and classical theta correspondences between $Mp_{2n}$ and special orthogonal groups.  
\end{abstract}

\maketitle

\tableofcontents

\section*{Introduction}

Let $\alg{G}$ be a connected reductive group over a local field $F$.  We write $G_F = \alg{G}(F)$.  Langlands, with later refinements by Deligne, Arthur, and Vogan, conjectured a finite-to-one parameterization of the irreducible admissible representations of $G_F$ (up to isomorphism) by certain homomorphisms from the Weil-Deligne group $\Weil_F'$ to the L-group ${}^L G^\vee$.  We refer to Vogan's article \cite{Vog} for precise versions of the ``local Langlands'' conjectures.

We adapt these conjectures to a large class of metaplectic groups, certain groups which are central extensions of $G_F$ by finite cyclic groups, when $\alg{G}$ is {\em split} over a local field $F$.  Namely, when $\tilde G_F$ is such a metaplectic group, one hopes for a L-group ${}^L \tilde G^\vee$ such that the smooth irreducible {\em genuine} complex representations of $\tilde G_F$ are naturally parameterized by equivalence classes of homomomorphisms from $\Weil_F'$ to ${}^L \tilde G^\vee$.

We describe the construction of such a group ${}^L \tilde G^\vee$, outlining the article below.  First, when $\Psi = (Y,X)$ is a root datum ($Y$ the cocharacter lattice, $X$ the character lattice), and $Q\colon Y \rightarrow \ZZ$ is a Weyl-invariant quadratic form, and $n$ is a positive integer, we explain a ``metaplectic modification'' of the root datum:  $\tilde \Psi = (\tilde Y, \tilde X)$.  This can also be found in the recent work of McNamara \cite{McN}, and is closely related to an earlier construction of Lusztig in \cite[\S 2.2.5]{Lus1}.  This $\tilde \Psi$ is again a root datum, and we write $\tilde \Psi^\vee$ for the dual root datum.  Using Lusztig's construction \cite{Lus} (or the existence theorem of \cite[Expos\'e XXV]{SGA3}) there exists a split reductive group scheme $\alg{\tilde G}\LVee$ over $\ZZ$ with root datum $\tilde \Psi^\vee$.  This construction of $\alg{\tilde G}\LVee$ -- the ``metaplectic dual group'' -- is the content of the first section of the article.

In the second section, we describe the ``incarnation'' of metaplectic groups by cocycles, where the cocycles come from linear algebraic gadgets we call ``bisectors''.  Here, a bisector is a $\ZZ$-bilinear form $C\colon Y \times Y \rightarrow \ZZ / 2 \ZZ$, for which $C(y,y) \ident Q(y)$ mod $2$ (with the same quadratic form $Q$ as before).  Given a bisector $C$ (satisfying another technical condition which we call fairness), and a surjective map from $\alg{K}_2(F)$ to $\mu_n$, one may construct a ``metaplectic group'' $1 \rightarrow \mu_n \rightarrow \tilde G_F \rightarrow G_F \rightarrow 1$.  Moreover, the bisectors are the objects of a category we call $\Cat{Bis}_Q^{f,sc}$, and the metaplectic groups are the objects of a category (due to Brylinski and Deligne \cite{B-D}), and this construction of metaplectic groups is functorial.  Section 2 covers this functorial ``incarnation'' of split metaplectic groups.

The third section is the technical heart of this paper, and it contains the construction of an L-group associated to a split metaplectic group over a local field $F \not \isom \CC$.  More precisely, we construct a group scheme ${}^L \alg{\tilde G}\LVee$ over $\ZZ$, together with morphisms of group schemes over $\ZZ$:
$$\alg{\tilde G}\LVee \rightarrow {}^L \alg{\tilde G}\LVee \rightarrow \alg{\Gamma}_F$$
in which $\alg{\tilde G}\LVee$ is the group scheme constructed in the first section and $\alg{\Gamma}_F$ is the group scheme over $\ZZ$ ``underlying'' the absolute Galois group of $F$.  The L-group ${}^L \alg{\tilde G}\LVee$ is constructed using the data of a bisector $C$, and it depends functorially on this data.  Thus the L-group ${}^L \alg{\tilde G}\LVee$ and the metaplectic group $\tilde G_F$ are constructed via functors from a common source.  The construction of the group scheme ${}^L \alg{\tilde G}\LVee$ heavily uses Lusztig's explicit description of the Hopf algebra (over $\ZZ$) of regular functions $\OO(\alg{\tilde G}\LVee)$.  It also exploits a construction, known in the literature on Hopf algebras and quantum groups as the ``cocycle bicrossedproduct'' \cite{Maj}.  We call this the ``double-twist'', as we twist not only the multiplication, but also the comultiplication, in a Hopf algebra.  A Hopf algebra over $\ZZ$, denoted $\OO({}^L \alg{\tilde G}\LVee)$ is constructed as a double-twist of the tensor product $\OO(\Gamma_F) \otimes \OO(\alg{\tilde G}\LVee)$ via a compatible pair of cocycles $(\tau, \chi)$.  The L-group is then defined as $\Spec(\OO({}^L \alg{\tilde G}\LVee))$.

When the covering degree $n$ of the metaplectic group is odd, this L-group is simply the direct product ${}^L \alg{\tilde G}\LVee = \alg{\Gamma}_F \times \alg{\tilde G}\LVee$ of group schemes over $\ZZ$.  When $n = 2$, this L-group is {\em isomorphic} to the direct product, but {\em not uniquely} and not without {\em extension of scalars} to $\ZZ[i]$ ($i$ a primitive fourth root of unity in $\CC$).  When $n=2$, we find that to give such an isomorphism ${}^L \alg{\tilde G}\LVee \isom \alg{\Gamma}_F \times \alg{\tilde G}\LVee$  (of group schemes over $\ZZ[i]$), it suffices to make the following choices:
\begin{itemize}
\item
A choice of a nontrivial continuous additive character $\psi\colon F \rightarrow \CC^\times$.
\item
A choice of a ``tetractor'' $\kappa\colon \bar Y \rightarrow \ZZ / 4 \ZZ$ of $C$.  
\end{itemize}
Here $\bar Y = \tilde Y / \left( (2Y \cap \tilde Y) + \ZZ[\tilde \Delta^\vee] \right)$, where $\tilde Y$ is the cocharacter lattice in the modified root datum $\tilde \Psi$, $Y$ is the cocharacter lattice in the original root datum $\Psi$, and $\ZZ[\tilde \Delta^\vee]$ is the coroot lattice in the modified root datum $\tilde \Psi$.  The ``tetractor'' condition is the relation
$$\kappa(y_1 + y_2) - (\kappa(y_1) + \kappa(y_2)) = [\delta \circ C](y_1, y_2),$$
for all $y_1, y_2 \in \tilde Y$, where $\delta\colon \ZZ / 2 \ZZ \rightarrow \ZZ / 4 \ZZ$ is the unique injective group homomorphism.  The set of tetractors of $C$ is a torsor for the elementary abelian 2-group $\Hom(\bar Y, \ZZ / 2 \ZZ)$.

The fourth section verifies that our construction of an L-group associated to a metaplectic group is compatible with some important examples in the literature.  First, we check compatibility with previous work on metaplectic tori \cite{Wei1}; we are able to parameterize irreducible genuine representations of split metaplectic tori by ${}^L \tilde T^\vee$-valued parameters, where ${}^L \tilde T^\vee = {}^L \alg{\tilde T}\LVee(\CC)$, and ``parameters'' refer to continuous homomorphisms from the Weil group of $F$.  Importantly, this parameterization relies on no choices (no additive characters, no tetractors, etc.), indicating that our L-group, constructed by a complicated double-twist, is the ``natural'' one (at least for split tori).  Of course, any natural parameterization for reductive groups must be compatible with that for tori, so this indicates that the double-twist is required for reductive groups as well.

Next, we check compatibility with the ``classical'' metaplectic correspondences.  Beginning with the monumental work of Waldspurger \cite{Wal1}, continuing with Adams-Barbasch \cite{A-B} (in the archimedean case), and most recently with Gan-Savin \cite{GS1} (in the $p$-adic case), one now knows the following:  contingent on the local Langlands conjectures for special orthogonal groups, there exists a finite-to-one map from the set of isomorphism classes of irreducible genuine admissible representations of the metaplectic group $Mp_{2n}$ (a central extension of $\alg{Sp}_{2n}(F)$ by $\mu_2$) to the set of conjugacy classes of Weil-Deligne representations with values in $\alg{Sp}_{2n}(\CC)$.  This is compatible with our results in the following way: the finite-to-one map described above depends (only) on a choice of a nontrivial additive character $\psi$ of $F$ (really, the $F^{\times 2}$-orbit of such a character), and our isomorphism from the L-group ${}^L Mp_{2n}^\vee$ to $\Gamma_F \times Sp_{2n}(\CC)$ also depends precisely on the $F^{\times 2}$-orbit of a nontrivial additive character $\psi$.  The choice of a ``tetractor'' here disappears -- there are only two tetractors in this setting, and switching between them has the same effect as replacing an additive character $\psi$ by its complex conjugate; in other words, the choice of tetractor becomes nothing more than a single choice of a square root of $-1$ in $\CC$.

Finally, we check compatibility with Savin's unramified Shimura correspondence for simply-laced groups from \cite{Sav}.  There Savin demonstrates that when $F$ is a $p$-adic field, and $(p,n) = 1$, and $\alg{G}$ is a simply-laced simply-connected split simple group over $F$, there is an isomorphism from the Iwahori Hecke algebra of $\tilde G_F$ to the Iwahori Hecke algebra of a linear group $G_{F,n}$.  This is compatible with our construction in two ways:  First, the L-group of $G_{F,n}$ is precisely what we call $\alg{\Gamma}_F \times \alg{\tilde G}\LVee$.  Second, the choices made by Savin in the process of identifying these Hecke algebras are precisely the choices needed in our construction of an isomorphism (over $\ZZ[i]$) from ${}^L \alg{\tilde G}\LVee$ to $\alg{\Gamma}_F \times \alg{\tilde G}\LVee$.  In this setting, the choice of additive character essentially disappears (it becomes a choice of square root of $-1$, in the unramified setting), but the choice of tetractor becomes more evident.  Our results are compatible with similar ``unramified correspondences'', such as McNamara's metaplectic Satake correspondence \cite{McN} and the recent thesis of Reich \cite{Rei} and work of Finkelberg-Lysenko \cite{F-L} in the geometric setting.  But the compatibility with Savin's work is most striking, since it illustrates not only an isomorphism between groups but also compatibility in choices of isomorphism.

The appendix provides foundations for the theory of Hopf algebras in the symmetric monoidal category of complete locally convex topological vector spces over $\CC$, with tensor product given by the completed injective tensor product.  Since our putative L-groups were defined by their Hopf algebras, computations with Langlands parameters involve Hopf algebras instead of groups.  In order to mix locally compact groups (like the Weil group of a local field) with complex algebraic groups, we exhibit a class of Hopf algebras that can be used for both sorts of groups.  Such Hopf algebras may have other applications as well, e.g. to locally compact quantum groups.

\subsection*{Directions}

We expect (or at least hope) the constructions of this paper to be eventually subsumed by a better construction, perhaps one that proceeds functorially from a central extension $\alg{G}'$ of $\alg{G}$ by $\alg{K}_2$ to an L-group, rather than have both $\alg{G}'$ and the L-group arise as functors from a common incarnating source.  In particular, the whole ``category of bisectors'' is a useful gadget, good for computations, but a bit of a kludge.  It is certainly not sufficient for carrying out descent, necessary to construct L-groups for nonsplit metaplectic groups.

It is worrying, or at least suspicious, that some of the data used to incarnate a metaplectic group -- the compatibility map $\phi: E_{Q_{sc}} \rightarrow f_{sc}^\ast E_C$ of Section \ref{IncSC} -- is not used in our construction of an L-group.  We hope that incorporating this data appropriately might eliminate some technicalities with fair bisectors as well.  But, at the moment, there is scant evidence for L-groups of metaplectic groups except for simply-connected groups and tori, and any speculation about an L-group would be hard to support.  We hope to rectify this in a later paper, by studying some depth zero representations.

We expect a similar construction, using double-twists of Hopf algebras, to work in more generality; this will probably require an L-group defined over $\ZZ[\zeta]$ (with $\zeta$ a primitive $n^{th}$ root of unity).  If one is willing to sacrifice a bit of arithmetic content, one could work over a field instead and attempt to describe the L-group from the Tannakian perspective.  This should be an exercise in translating the double-twist into the category of representations.  Work of \cite{F-L} and \cite{Rei} suggests that the ``second twist'' (written $\chi$ in this paper) should correspond to a twist in the commutativity constraint of the Tannakian category; but our twist involves not only $\alg{\tilde G}\LVee$ but also the Galois group $\alg{\Gamma}_F$.  The ``first twist'' is already of group-theoretic origin, and presumably simple to translate to the Tannakian category.

We hope to approach the nonsplit case in a future paper, and find evidence in the form of a depth zero parameterization for unramified metaplectic groups in the spirit of DeBacker-Reeder \cite{D-R}.

By working with group schemes over $\ZZ$ whenever possible, we have left the door open to studying a mod $\ell$ local Langlands correspondence for metaplectic groups, in the spirit of Vigneras \cite{Vig}.  It is hoped that the constructions of this paper will globalize, and working with schemes over $\ZZ$ will pay off in arithmetic applications.

\subsection*{Acknowledgments}

During the preparation of an earlier paper on metaplectic tori \cite{Wei1}, we benefited greatly from a correspondence with P. Deligne.  His advice on that paper was very helpful, and suggested many reasons why the results there are not entirely satisfactory.  Deligne's recommendation to consider the geometric work of Finkelberg-Lysenko \cite{F-L}, seconded by a later recommendation of B. Gross to consider the geometric setting in the recent thesis of R. Reich \cite{Rei}, was crucial to making a reasonable guess at an L-group.  We also thank B. Gross for pointing us towards the recent work of Gan, Gross, and Prasad \cite[\S 11]{GGP}, in which there are precise conjectures on a local Langlands parameterization for the metaplectic group $Mp_{2n}$.  

We are also thankful for the advice and personal communication with Wee Teck Gan and Gordan Savin.  Both have thought deeply about metaplectic groups; Savin's earlier paper \cite{Sav} and their recent work together \cite{GS1} were very helpful, as was their encouragement.

\section{Dual groups over $\ZZ$}

Let $F$ be a field.  $\alg{G}$ will be a split connected reductive group over $F$ with split maximal torus $\alg{T}$.  In this section, we describe how a Weyl-invariant quadratic form on the cocharacter lattice of $\alg{T}$ can be used to modify the root datum of $(\alg{G}, \alg{T})$; this construction can also be found in \cite[\S 11]{McN}.  A ``modified dual'' root datum can then be used to construct a group scheme $\alg{\tilde G}\LVee$ over $\ZZ$.

\subsection{Cartan and root data}

We recall the notion of Cartan datum and root datum, using Lusztig's definitions from \cite[\S 2.1]{Lus1} throughout.  
\subsubsection{Cartan data}
A {\em Cartan datum} is a pair $(I, \cdot)$ where $I$ is a finite set and $\cdot$ is a symmetric $\ZZ$-valued bilinear form on $\ZZ[I]$, such that 
\begin{enumerate}
\item
$0 < i \cdot i \in 2 \ZZ$ for all $i \in I$;
\item
$0 \geq 2 (i \cdot j) / (i \cdot i) \in \ZZ$ for all $i \neq j \in I$.
\end{enumerate}
We exclusively work with Cartan data {\em of finite type}, meaning that the symmetric bilinear form is positive definite.  The {\em Cartan matrix} associated to $(I, \cdot)$ is the matrix $(a_{ij})$ with rows and columns indexed by $I$, with entries $a_{ij} = 2 (i \cdot j) / (i \cdot i)$.

Let $h(i,j)$ be $2$, $3$, $4$, or $6$, according to whether $4 (i \cdot j)(j \cdot i) / (i \cdot i)(j \cdot j)$ is $0$, $1$, $2$, or $3$, respectively.   Let $\Braid_{(I,\cdot)}$ be the braid group with generators $\{ s_i : i \in I \}$ and relations $s_i s_j s_i \cdots = s_j s_i s_j \cdots$ where each side has $h(i,j)$ factors, as defined in \cite[\S 2.1.1]{Lus1}.  Let $\Weyl_{(I,\cdot)}$ be the Weyl group -- the quotient of $\Braid_{(I,\cdot)}$ by the relations $s_i^2 = 1$ for all $i \in I$.  


When $(I, \cdot)$ is a Cartan datum, we define $m_I$ to be the smallest positive integer for which $m_I / 2 (i \cdot i) \in \ZZ$ for all $i \in I$.  In this case, define the {\em dual Cartan datum} $(I, \checkdot)$ with the same underlying set, but with symmetric bilinear form given by
$$i \checkdot j = m_I \cdot \frac{i \cdot j}{(i \cdot i )( j \cdot j)} = m_I a_{ij}/ 2 (j \cdot j) .$$
Note that $i \checkdot i = m_I / (i \cdot i) \in 2 \ZZ$, so that $i \cdot i > j \cdot j$ implies $i \checkdot i < j \checkdot j$.  Also, we have
$$2 \frac{i \checkdot j}{i \checkdot i} = 2 \frac{i \cdot j}{j \cdot j} \in \ZZ.$$
It follows that
$$\frac{ 4 (i \checkdot j)(j \checkdot i)}{(i \checkdot i)(j \checkdot j)} = \frac{4 (i \cdot j)(j \cdot i)}{(j \cdot j)(i \cdot i)}.$$
Hence $\Braid_{(I, \checkdot)} = \Braid_{(I, \cdot)}$ and $\Weyl_{(I, \checkdot)} = \Weyl_{(I, \cdot)}$ (the constants $h(i,j)$ are the same in $(I, \cdot)$ and in $(I, \checkdot)$).


\subsubsection{Root data}
Let $(I, \cdot)$ be a Cartan datum with Cartan matrix $(a_{ij})$.  Recall all Cartan datum are assumed to be of finite type.  A {\em root datum} of type $(I, \cdot)$ consists of
\begin{enumerate}
\item
A pair $(Y,X)$ of finitely-generated free $\ZZ$-modules and a perfect pairing $\langle \bullet, \bullet \rangle \from Y \times X \rightarrow \ZZ$;
\item
Embeddings $(i \mapsto \alpha_i^\vee)$ from $I$ into $Y$ and $(i \mapsto \alpha_i)$ from $I$ into $X$ such that
$$\langle \alpha_i^\vee, \alpha_j \rangle = a_{ij} , \mbox{ for all } i,j \in I.$$
\end{enumerate}

Given a root datum $(Y,X)$ of type $(I, \cdot)$, we write $\Delta^\vee$ for the image of $I$ in $Y$, and $\Delta$ for the image of $I$ in $X$; these are the sets of {\em simple coroots} and {\em simple roots}, respectively.  Since $(I, \cdot)$ is of finite type, $\ZZ[\Delta^\vee]$ is a subgroup of $Y$ (called the {\em coroot lattice}) and $\ZZ[\Delta]$ is a subgroup of $X$ (called the {\em root lattice}). 

The root data of type $(I,\cdot)$ are the objects of a category described in \cite[\S 2.2.2]{Lus1}, with initial object given by the simply-connected root datum (in which $\ZZ[\Delta^\vee] = Y$) and final object given by the adjoint root datum (in which $\ZZ[\Delta] = X$).  For this reason, we write $Y_{sc}$ for $\ZZ[\Delta^\vee]$ and $X_{ad}$ for $\ZZ[\Delta]$;  we also write $f_{sc}\colon Y_{sc} \rightarrow Y$ for the inclusion map. 

Let $(Y,X)$ be a root datum of type $(I,\cdot)$.  For each $i \in I$, let $s_i \from Y \rightarrow Y$ be the $\ZZ$-linear map given by $s_i(y) = y - \langle y, \alpha_i \rangle \alpha_i^\vee$.  Also, we write $s_i\colon X \rightarrow X$ for the $\ZZ$-linear map given by $s_i(x) = x - \langle \alpha_i^\vee, x \rangle \alpha_i$.  These define actions of $\Weyl = \Weyl_{(I, \cdot)}$ on $X$ and $Y$, for which the pairing is $\Weyl$-invariant:
$$\langle s_i(y), s_i(x) \rangle = \langle y, x \rangle, \mbox{ for all } y \in Y, x \in X, i \in I.$$
Let $\Phi^\vee = \Weyl \cdot \Delta^\vee$ and $\Phi = \Weyl \cdot \Delta$ be the orbits resulting from this action.  These are the sets of {\em coroots} and {\em roots} of the root datum, respectively. 

When $\Psi = (Y,X)$ is a root datum of type $(I, \cdot)$, the {\em dual root datum} $\Psi^\vee = (X,Y)$ of type $(I, \checkdot)$ is given by the same embeddings of $I$ into $X$ and $I$ into $Y$, and the pairing $\langle \bullet, \bullet \rangle^\vee$ given by $\langle x,y \rangle^\vee= \langle y,x \rangle$ for all $x \in X$, $y \in Y$.  Note that
$$\langle \alpha_i, \alpha_j^\vee \rangle^\vee = \langle \alpha_j^\vee, \alpha_i \rangle = \frac{2 j \cdot i}{j \cdot j} = \frac{2 i \checkdot j}{i \checkdot i}.$$
The roots of the root datum $\Psi^\vee$ are the coroots of $\Psi$, and the coroots of $\Psi^\vee$ are the roots of $\Psi$.  The Weyl groups of $\Psi$ and $\Psi^\vee$ are equal.    

\subsection{Reductive groups}
Let $\alg{G}$ be a split connected reductive group over $F$.  Let $\alg{T} \subset \alg{B} \subset \alg{G}$ be an $F$-split maximal torus contained in a Borel subgroup of $\alg{G}$, defined over $F$.  Let $X = X(\alg{T})$ and $Y = Y(\alg{T})$ be the groups of characters and cocharacters of $\alg{T}$, respectively.  These are finite-rank free $\ZZ$-modules, endowed with a canonical $\ZZ$-valued perfect pairing.  Let $\Phi \subset X$ be the set of roots for the action of $\alg{T}$ on the Lie algebra $\alg{Lie}(\alg{G})$, and let $\Delta \subset \Phi$ be the set of positive simple roots with respect to the choice of Borel subgroup $\alg{B}$.  

For each $\alpha \in \Phi$, let $\alg{G}_\alpha$ be the smallest Zariski-closed subgroup of $\alg{G}$ containing the root subgroups for $\pm \alpha$.  There exists a unique $\alpha^\vee \in Y$ such that $\alpha^\vee$ has image in $\alg{G}_\alpha$ and $\langle \alpha^\vee, \alpha \rangle = 2$.  Write $\Phi^\vee$ for the resulting set of coroots, and $\Delta^\vee$ for the subset of simple coroots.

There exists a Cartan datum $(I, \cdot)$ of finite type, and bijections $(i \mapsto \alpha_i^\vee)$, $(i \mapsto \alpha_i)$ from $I$ to $\Delta^\vee$ and $\Delta$ respectively, for which $(Y,X)$ becomes a root datum of type $(I, \cdot)$.  Thus we find a root datum $\Psi = \Psi(\alg{G}, \alg{B}, \alg{T})$ of type $(I, \cdot)$ associated to the group $\alg{G}$. 

\subsection{Metaplectic root datum}

The following structure occurs naturally in the study of metaplectic groups, when using the framework of Brylinski and Deligne \cite{B-D}.
\begin{definition}
A {\em metaplectic structure} on the root datum $(Y,X)$ (of type $(I, \cdot)$) is a Weyl-invariant quadratic form $Q\colon Y \rightarrow \ZZ$ and a positive integer $n$.
\end{definition}
Whenever $Q\colon Y \rightarrow \ZZ$ is a quadratic form, we write $B\colon Y \times Y \rightarrow \ZZ$ for the associated bilinear form:
$$B(y_1, y_2) = Q(y_1 + y_2) - (Q(y_1) + Q(y_2)).$$
\begin{lemma}
\label{BQLemma}
Let $Q$ be a Weyl-invariant quadratic form from $Y$ to $\ZZ$.  Then for all $i \in I$,and all $y \in Y$, we have
$$B(\alpha_i^\vee, y) = Q(\alpha_i^\vee) \cdot \langle y, \alpha_i \rangle.$$
\end{lemma}
\proof
The Weyl-invariance of $Q$ implies the Weyl-invariance of $B$, from which we find
\begin{eqnarray*}
B(\alpha_i^\vee, y) & = & B(s_i(\alpha_i^\vee), s_i(y)) \\
& = & B(-\alpha_i^\vee, y - \langle y, \alpha_i \rangle \alpha_i^\vee) \\
& = & B(\alpha_i^\vee, \alpha_i^\vee) \langle y, \alpha_i \rangle  - B(\alpha_i^\vee, y) \\
& = & 2 Q(\alpha_i^\vee) \cdot \langle y, \alpha_i \rangle - B(\alpha_i^\vee, y).
\end{eqnarray*}
Hence $2 B(\alpha_i^\vee, y) = 2 Q(\alpha_i^\vee) \cdot \langle y, \alpha_i \rangle$.  The result follows immediately.
\qed

Recently, a number of authors (\cite{McN}, \cite{Rei}, \cite{F-L}) have found a ``dual group'' to a metaplectic group, or at least the root datum thereof.  The combinatorics of this construction are below.
\begin{construction}[Cf.\ {\cite[Theorem 11.1]{McN}}]
\label{MetaConst}
Let $(Q,n)$ be a metaplectic structure on the root datum $(Y,X)$ of type $(I, \cdot)$.  For each $i \in I$, let $n_i$ be the smallest positive integer for which $n_i \cdot Q(\alpha_i^\vee) \in n \ZZ$.  Define a symmetric bilinear form $\odot$ on $\ZZ[I]$ by
$$i \odot j = \frac{n^2}{n_i n_j} (i \cdot j).$$
Then $(I, \odot)$ is a Cartan datum of finite type.  Define also
\begin{itemize}
\item
$\tilde Y = \{ y \in Y : B(y, y') \in n \ZZ \mbox{ for all } y' \in Y \}$.
\item
$\tilde \alpha_i^\vee = n_i \alpha_i^\vee$, for all $i \in I$.  
\item
$\tilde X = \{ x \in X \otimes \QQ : \langle y,x \rangle \in \ZZ \mbox{ for all } y \in \tilde Y$.
\item
$\tilde \alpha_i = n_i^{-1} \alpha_i$, for all $i \in I$.
\end{itemize}
This forms a root datum $\tilde \Psi = (\tilde Y, \tilde X)$ of type $(I, \odot)$, via the embeddings $(i \mapsto \tilde \alpha_i^\vee)$ and $(i \mapsto \tilde \alpha_i)$.  We write $\tilde \Delta$, $\tilde \Delta^\vee$ for the set of simple roots and simple coroots in this new root datum, and $\tilde \Phi$, $\tilde \Phi^\vee$ for the roots and coroots.  We write $\tilde \Psi^\vee$ for the root datum dual to $\tilde \Psi$.
\end{construction}
\proof
By Lemma \ref{BQLemma} we find that 
$$B(\tilde \alpha_i^\vee, y) = n_i B(\alpha_i^\vee, y) = n_i Q(\alpha_i^\vee) \langle y, \alpha_i \rangle \in n \ZZ,$$
for all $y \in Y$.  Hence $\tilde \alpha_i^\vee \in \tilde Y$.  Again by Lemma \ref{BQLemma}, we find that
$$\langle \tilde y, \alpha_i \rangle Q(\alpha_i^\vee) = B(\alpha_i^\vee, \tilde y) \in n \ZZ,$$
for all $\tilde y \in \tilde Y$.  Since $n_i$ is the smallest positive integer for which $n_i Q(\alpha_i^\vee) \in n \ZZ$, it follows that $\langle \tilde y, \alpha_i \rangle \in n_i \ZZ$.  Hence $\langle \tilde y, \tilde \alpha_i \rangle \in \ZZ$ for all $\tilde y \in \tilde Y$.  Therefore $\tilde \alpha_i \in \tilde X$ as required.

We find that 
\begin{equation}
\label{Scale}
\langle \tilde \alpha_i^\vee, \tilde \alpha_j \rangle = n_i n_j^{-1} \langle \alpha_i^\vee, \alpha_j \rangle = \frac{2 n_i ( i \cdot j) }{n_j (i \cdot i)} = \frac{ 2 n^2 n_i n_i n_j (i \odot j) }{n^2 n_i n_i n_j (i \odot i)} = \frac{2 i \odot j}{i \odot i}.
\end{equation}

To see that $(I, \odot)$ is a Cartan datum of finite type, we find that $\odot$ gives a positive definite symmetric bilinear form on $\ZZ[I]$ since the bilinear form $\odot$ may be obtained from $\cdot$ by a $\QQ$-linear change of basis.  Moreover $\odot$ is an integer-valued bilinear form and $(i \odot i) \in 2 \ZZ$ since $i \odot j$ is an integer multiple of $i \cdot j$ for all $i,j \in I$.  By Equation \ref{Scale} and the facts that $\tilde \alpha_i^\vee \in \tilde Y, \tilde \alpha_j \in \tilde X$, we find that $0 \geq (2 i \odot j)/(i \odot i) \in \ZZ$, so $(I, \odot)$ is a Cartan datum.  Equation \ref{Scale} also implies that $\tilde \Psi$ is a root datum of type $(I, \odot)$.
\qed

\begin{remark}
This construction bears more than a small resemblence to the construction of \cite[\S 2.2.4, 2.2.5]{Lus1}.
\end{remark}

\subsection{Dual groups over $\ZZ$}

We keep the root datum $\tilde \Psi = (\tilde Y, \tilde X)$ of type $(I, \odot)$, from the previous section.  Let $\tilde \Psi^\vee = (\tilde X, \tilde Y)$ be the dual root datum (of type $(I, \checkodot)$).  We review Lusztig's integral form of a group scheme associated to $\tilde \Psi^\vee$ from \cite{Lus}, which completes and generalizes earlier work by Kostant \cite{Kos}.  Existence of such a group scheme follows from \cite[Exp. XXV, Th\'eor\`eme 1.1]{SGA3}, but we will heavily rely on Lusztig's construction in later sections. 

For any $\ZZ$-module $M$, let $M^\dual = \Hom_\ZZ(M, \ZZ)$ denote the full algebraic dual $\ZZ$-module.  If $M$ is a free $\ZZ$-module endowed with a $\ZZ$-basis ${\mathcal M}$, and $m \in {\mathcal M}$, we write $m^\dual$ for the unique element of $M^\dual$ satisfying $m^\dual(m) = 1$ and $m^\dual(m') = 0$ whenever $m \neq m' \in {\mathcal M}$.

Let $\Theta_\QQ$ be the associative unital $\QQ$-algebra, with generators $\{ \theta_i : i \in I \}$ and the Serre relations:
$$\sum_{ {p, q \in \NN} \atop {p+q=1 - \langle \tilde \alpha_i^\vee, \tilde \alpha_j \rangle} } (-1)^q \frac{\theta_i^p}{p!} \theta_j \frac{\theta_i^q}{q!} = 0,$$
for all $i \neq j \in I$.  This is the construction described in \cite[\S 1.3]{Lus}, though we use the letter $\Theta$ instead of Lusztig's $\alg{f}$.   For all $p \in \NN$, let $\theta_i^{(p)}$ be the divided power element $\theta_i^p / p! \in \Theta_\QQ$.  We write $\Theta$ for the (unital) $\ZZ$-subalgebra of $\Theta_\QQ$ generated by $\theta_i^{(p)}$ for all $i \in I$ and all $0 \leq p \in \ZZ$.

Let $\BB_\Theta$ be Lusztig's canonical basis of the $\ZZ$-module $\Theta$, discussed in \cite[\S 14.4]{Lus1}.  Let $\OO(\alg{\tilde N}\LVee)$ be the $\ZZ$-submodule of $\Theta^\dual$ spanned by $b^\dual$ for $b \in \BB_\Theta$ (the notation will be justified soon). 

Let $\UU$ be the (non-unital, in general) associative $\ZZ$-algebra generated by symbols $f_1^+ 1_y f_2^-$ and $f_1^- 1_yf_2^+$, for all $f_1, f_2 \in \Theta$ and $y \in \tilde Y$, subject to the relations of \cite[\S 1.4]{Lus} (N.B.\ our $\tilde Y$ is Lusztig's $X$, since we work with the modified dual root datum $\tilde \Psi^\vee = (\tilde X,\tilde Y)$).  As a $\ZZ$-module, there is an isomorphism
$$\Theta \otimes \ZZ[\tilde Y] \otimes \Theta \rightarrow \UU, \quad f_1 \otimes y \otimes f_2 \mapsto f_1^- 1_y f_2^+.$$

Let $\BB$ be Lusztig's canonical $\ZZ$-basis of $\UU$, described in \cite[\S 25]{Lus1} and \cite[\S 1.4]{Lus}.  What Lusztig writes as $\dot \BB$ and $\dot \UU$, we write unadorned as $\BB$ and $\UU$.  For all $y \in \tilde Y$, and $b \in \BB_\Theta$, the elements $b^- 1_y$ and $1_y b^+$ of $\UU$ are contained in the canonical basis $\BB$.  In particular, $1_y \in \BB$ for all $y \in \tilde Y$.  These elements satisfy the relationship $1_{y_1} \cdot 1_{y_2} = 0$ if $y_1 \neq y_2$, and $1_y \cdot 1_y = 1_y$, for all $y, y_1, y_2 \in Y$.  Let $\OO(\alg{\tilde G}\LVee)$ be the $\ZZ$-submodule of $\UU^\dual$ spanned by $b^\dual$ for all $b \in \BB$.  

In \cite[\S 3.8, 3.9]{Lus}, Lusztig describes commutative Hopf algebra structures on $\OO(\alg{\tilde N}\LVee)$ and on $\OO(\alg{\tilde G}\LVee)$. 
Justifying the notation for these Hopf algebras, we define group schemes over $\ZZ$ by:
$$\alg{\tilde G}\LVee = \Spec(\OO(\alg{\tilde G}\LVee)), \quad \alg{\tilde N}\LVee = \Spec(\OO(\alg{\tilde N}\LVee)), \quad \alg{\tilde T}\LVee = \Spec(\ZZ[\tilde Y]).$$
From \cite[Theorem 4.11]{Lus}, it follows that $\alg{\tilde G}\LVee$ is a split reductive group scheme over $\ZZ$, with root datum $\tilde \Psi^\vee$.  In other words, $\alg{\tilde G}\LVee$ is a smooth affine group scheme over $\ZZ$, with connected reductive geometric fibres, and $\alg{\tilde T}\LVee$ is a split maximal torus in $\alg{\tilde G}\LVee$ over $\ZZ$.  

Moreover, the construction yields a maximal torus $\alg{\tilde T}\LVee$ contained in a Borel subgroup $\alg{\tilde T}\LVee \alg{\tilde N}\LVee$ of $\alg{\tilde G}\LVee$.

For each $y \in \tilde Y$, there is a basis element $1_y \in \BB$ and we write $\BB_\alg{\tilde T} = \{ 1_y : y \in \tilde Y \}$ for the resulting subset of $\BB$.  The unit element of $\OO(\alg{\tilde G}\LVee)$ is $1_0^\dual$ (where $0$ denotes the zero element of the abelian group $\tilde Y$).  The counit for $\OO(\alg{\tilde G}\LVee)$ is given by $\epsilon(a^\dual) = 0$ if $a \not \in \BB_{\alg{\tilde T}}$ and $\epsilon(a^\dual) = 1$ if $a \in \BB_{\alg{\tilde T}}$.

When $M$ is a $\ZZ$-module, we say ``define a $\ZZ$-linear map $f\colon \OO(\alg{\tilde G}\LVee) \rightarrow M$ by $f(b^\dual) = $...'' to mean ``let $f$ be the unique $\ZZ$-linear map whose restriction to basis elements $b^\dual \in \BB$ is...''.
\begin{definition}
A $\ZZ$-linear map $f\colon \OO(\alg{\tilde G}\LVee) \rightarrow M$ is called {\em toral} if $f(b^\dual) = 0$ for all $b \not \in \BB_{\alg{\tilde T}}$.  To define a toral linear map $f$, it suffices to give $f(1_y^\dual)$ for all $y \in \tilde Y$.
\end{definition}

The algebra $\UU$ has a direct sum decomposition into free $\ZZ$-submodules:
$$\UU = \bigoplus_{y_1, y_2 \in \tilde Y} 1_{y_1} \UU 1_{y_2}.$$
We frequently use the fact (observed in \cite[\S 23.1.2]{Lus1}, and a simple consequence of the relations of \cite[\S 23.1.3]{Lus1}) that
\begin{equation}
\label{USC}
1_{y_1} \UU 1_{y_2} = 0 \mbox{ unless } (y_1 - y_2) \in \ZZ[\tilde \Delta^\vee].
\end{equation}
Here we recall that $\ZZ[\tilde \Delta^\vee]$ is the coroot lattice of the root datum $\tilde \Psi$, or equivalently the image of the coroot lattice of the simply-connected root datum $\tilde \Psi^{sc}$ under the canonical morphism from $\tilde \Psi^{sc}$ to $\tilde \Psi$.

The canonical basis is partitioned compatibly with the direct sum decomposition of $\UU$:
$$\BB = \bigsqcup_{y_1, y_2 \in \tilde Y} \left( 1_{y_1} \UU 1_{y_2} \cap \BB \right).$$
If $b \in \BB \cap 1_{y_1} \UU 1_{y_2}$, then we write $\lt{b} = y_1$ and $\rt{b} = y_2$.  We write $\nt{b} = \rt{b} - \lt{b} \in \ZZ[\tilde \Delta^\vee]$.


We use structure constants to capture the Hopf algebra structure on $\OO(\alg{\tilde G}\LVee)$ using Lusztig's notation.  We write $\hat m_c^{a,b}$ and $m_{a,b}^c$ for the (integer) structure constants of multiplication and comultiplication as in \cite[\S 3.1]{Lus}, so that
$$a^\dual \cdot b^\dual = \sum_{c \in \BB} \hat m_c^{a,b} c^\dual, \quad \Delta(c^\dual) = \sum_{a,b \in \BB} m_{a,b}^c a^\dual \otimes b^\dual.$$
We freely use relationships among these structure constants listed in \cite[\S 1.5]{Lus}.  In particular, we use the fact that for all $y \in \tilde Y$,
$$
\hat m_{1_y}^{a,b} = \begin{cases} 1 & \mbox{ if } \exists y_1, y_2 \in \tilde Y, \mbox{ for which } a = 1_{y_1}, b = 1_{y_2}, \mbox{ and } y_1 + y_2 = y, \\ 0& \mbox{ otherwise}. \end{cases}
$$

We also use the fact that if $y \in \tilde Y$, and $b \in \BB$, then
$$m_{1_y, b}^c = \begin{cases}  1 & \mbox{ if } b = c \mbox{ and }  \lt{b} = \lt{c} = y, \\ 0 & \mbox{ otherwise.} \end{cases}$$
Symmetrically,
$$m_{b, 1_y}^c = \begin{cases}  1 & \mbox{ if } b = c \mbox{ and }  \rt{b} = \rt{c} = y, \\ 0 & \mbox{ otherwise.} \end{cases}$$

By associativity and coassociativity, there exist unique integer structure constants with more indices, satisfying the following:
$$(a_1^\dual \cdot a_2^\dual) \cdot a_3^\dual = a_1^\dual \cdot (a_2^\dual \cdot a_3^\dual) = \sum_{a \in \BB} \hat m_{a}^{a_1, a_2, a_3} a^\dual,$$
$$(\Id \otimes \Delta)\Delta(a^\dual) = (\Delta \otimes \Id)\Delta(a^\dual) = \sum_{a_1, a_2, a_3 \in \BB} m_{a_1, a_2, a_3}^a (a_1^\dual \otimes a_2^\dual \otimes a_3^\dual).$$
In what follows, we utilize structure constants with three, four, or five indices, keeping associativity and coassociativity in mind throughout.

To summarize our constructions up to this point, we have the following.
\begin{construction}
Let $\alg{G}$ be a split connected reductive group over a field $F$.  Let $\Psi = (Y,X)$ denote its root datum, of type $(I, \cdot)$.

To each metaplectic structure $(Q,n)$ on $(Y,X)$, we may associate a root datum $\tilde \Psi = (\tilde Y, \tilde X)$.  To the dual root datum $\tilde \Psi^\vee = (\tilde X, \tilde Y)$, Lusztig associates a commutative Hopf algebra $\OO(\alg{\tilde G}\LVee)$ over $\ZZ$, whence we have a group scheme $\alg{\tilde G}\LVee$ over $\ZZ$.

We say that $\alg{\tilde G}\LVee$ is the {\em dual group} for $(\alg{G}, Q,n)$.
\end{construction}

When $F$ is a local field, this dual group is closely related to spherical unramified represenations of a ``metaplectic group'' $\tilde G_F$, under some additional hypotheses and after making many additional choices.

\section{Metaplectic groups}

We keep the notation of the previous section, where $\alg{G}$ is a split connected reductive group over a field $F$.  As before, $\alg{T}$ is a split maximal torus, contained in a Borel subgroup $\alg{B}$ of $\alg{G}$, and $(Y,X)$ the associated root datum of type $(I, \cdot)$.

In \cite{B-D}, Brylinski and Deligne study and classify central extensions of $\alg{G}$ by $\alg{K}_2$, in the category of sheaves of groups on the big Zariski site of schemes of finite type over $F$.  We write $\Cat{CExt}(\alg{G}, \alg{K}_2)$ for their category of central extensions.  In practice, objects of this category are constructed via cocycles; we describe this practice with a series of ``incarnation functors''.    

\subsection{Incarnating $\Cat{CExt}(Y, F^\times)$}

\begin{definition}
Let $\Cat{CExt}(Y, F^\times)$ denote the groupoid of central extensions 
$$1 \rightarrow F^\times \rightarrow E \rightarrow Y \rightarrow 1$$ of $Y$ by $F^\times$ (in the category of groups).  A morphism from $F^\times \rightarrow E \rightarrow Y$ to $F^\times \rightarrow E' \rightarrow Y$ in $\Cat{CExt}(Y, F^\times)$ is a group homomorphism $\phi\colon E \rightarrow E'$, which makes the following diagram commute.
$$\xymatrix{F^\times \ar[r] \ar[d]^= & E \ar[r] \ar[d]^{\phi} & Y \ar[d]^= \\ F^\times \ar[r] & E' \ar[r] & Y}$$
\end{definition}

Each central extension $F^\times \rightarrow E \rightarrow Y$ determines a commutator pairing $[\bullet, \bullet]\colon \bigwedge^2 Y \rightarrow F^\times$.  The universal coefficient theorem gives an exact sequence
$$0 = Ext(Y, F^\times) \rightarrow H^2(Y, F^\times) \rightarrow \Hom(\bigwedge^2 Y, F^\times).$$
It follows that the isomorphism class of a central extension $F^\times \rightarrow E \rightarrow Y$ is uniquely determined by its commutator pairing.

Let $Q\colon Y \rightarrow \ZZ$ be a quadratic form, and $B$ the associated bilinear form:
$$B(y_1, y_2) = Q(y_1 + y_2) - (Q(y_1) + Q(y_2)).$$
Let $\Cat{CExt}_Q(Y, F^\times)$ be the full subcategory of $\Cat{CExt}(Y, F^\times)$ whose objects are central extensions with commutator
$$[y_1, y_2] = (-1)^{B(y_1, y_2)}, \mbox{ for all } y_1, y_2 \in Y.$$

\begin{definition}
Let $Q\colon Y \rightarrow \ZZ$ be a quadratic form.  A {\em bisector} of $Q$ is a bilinear form $C\colon Y \times Y \rightarrow \ZZ / 2 \ZZ$ such that $Q(y) \equiv C(y,y)$, mod $2$, for all $y \in Y$.
\end{definition}

Bisectors of $Q$ are the objects of a category defined here.
\begin{definition}
Let $Q\colon Y \rightarrow \ZZ$ be a quadratic form.  Let $\Cat{Bis}_Q$ denote the category whose objects are bisectors of $Q$, and where the set of morphisms from $C_1$ to $C_2$ is the set of functions $H\colon Y / 2 Y \rightarrow \ZZ / 2 \ZZ$ such that
$$C_2(y_1, y_2) - C_1(y_1, y_2) = H(y_1 + y_2) - H(y_1) - H(y_2).$$
Composition of morphisms is given by pointwise addition:  $[H_1 \circ H_2](y) = H_1(y) + H_2(y)$.
\end{definition}

\begin{proposition}
The category $\Cat{Bis}_Q$ is a groupoid with one isomorphism class.  
\end{proposition}
\proof
The category $\Cat{Bis}_Q$ is nonempty -- every quadratic form $Q(y)$, mod $2$, can be expressed as $C(y,y)$ for {\em some} (not necessarily symmetric) bilinear form $C$.

Now let $C_1$ and $C_2$ be objects of $\Cat{Bis}_Q$.  The function $[C_2 - C_1]$ can be seen as a group cocycle; $[C_2 - C_1] \in Z^2(Y / 2 Y, \ZZ / 2 \ZZ)$.  The cohomology class of such a cocycle is determined by its associated quadratic form, so since $[C_2 - C_1](y,y) = 0$ we find that $C_2 - C_1$ is a coboundary: $C_2 - C_1 = \partial H$ (for some 1-chain $H$).  The function $H \colon Y / 2 Y \rightarrow  \ZZ / 2 \ZZ$ is a morphism from $C_1$ to $C_2$ in the category $\Cat{Bis}_Q$.  The inverse of the morphism $H$ is $-H$, so $H$ is an isomorphism.
\qed  

\begin{construction}
The following data determines an essentially surjective functor from $\Cat{Bis}_Q$ to $\Cat{CExt}_Q(Y, F^\times)$:  to each object $C$ of $\Cat{Bis}_Q$ we associate the central extension $E_C$ of $Y$ by $F^\times$ which is $Y \times F^\times$ as a set and in which multiplication is given by
$$(y_1, f_1) \cdot_C (y_2, f_2) = (y_1 + y_2, (-1)^{C(y_1, y_2)} f_1 f_2).$$

For each morphism $H\colon C_1 \rightarrow C_2$, we associate the isomorphism $\phi_H\colon E_{C_1} \rightarrow E_{C_2}$ given by
$$\phi_H(y, f) = (y, (-1)^{H(y)} f).$$
\end{construction}
\proof
The commutator pairing of the central extension $E_C$ is 
$$[y_1, y_2] = (-1)^{C(y_1, y_2) + C(y_2, y_1)}.$$
If $C(y,y) \ident Q(y)$, mod $2$, then $C(y_1, y_2) + C(y_2, y_1) \ident B(y_1, y_2)$, mod $2$, and so $E_C$ is an object of $\Cat{CExt}_Q(Y, F^\times)$.

For morphisms, note that $\phi_H$ is invertible and 
\begin{eqnarray*}
\phi_H( (y_1, f_1) \cdot_{C_1} (y_2, f_2) ) & = & \phi_H \left( y_1 + y_2, f_1 f_2 (-1)^{C_1(y_1, y_2)} \right) \\
& = & \left( y_1 + y_2, f_1 f_2 (-1)^{C_1(y_1, y_2) + H(y_1 + y_2)} \right) \\
& = & \left( y_1 + y_2, f_1 f_2 (-1)^{C_2(y_1, y_2) + H(y_1) + H(y_2)} \right) \\ 
& = & \phi_H(y_1, f_1) \cdot_{C_2} \phi_H(y_2, f_2).
\end{eqnarray*}
Hence $\phi_H$ is an isomorphism from $E_{C_1}$ to $E_{C_2}$ whenever $H\colon C_1 \rightarrow C_2$ is a morphism in $\Cat{Bis}_Q$.  Since every quadratic form $Q$ has a bisector, the functor is essentially surjective.
\qed

\subsection{Incarnating split metaplectic tori}

In \cite[Proposition 3.11]{B-D}, Brylinski and Deligne define an equivalence of categories from $\Cat{CExt}(\alg{T}, \alg{K}_2)$ to $\bigsqcup_Q \Cat{CExt}_Q(Y, F^\times)$ -- the groupoid obtained as the disjoint union of the groupoids $\Cat{CExt}_Q(Y, F^\times)$ as $Q$ ranges over all quadratic forms on $Y$.  The inverse of this functor is well-defined up to natural isomorphism; consider the resulting composition of functors:
$$\Cat{Inc}_{\alg{T}}:  \Cat{Bis}_Q \rightarrow \Cat{CExt}(Y, F^\times) \rightarrow \Cat{CExt}(\alg{T}, \alg{K}_2).$$

We call this composite functor the {\em incarnation} functor; if $C$ is an object of $\Cat{Bis}_Q$ we write $\alg{T}_C' = \Cat{Inc}(C)$ for the associated central extension of $\alg{T}$ by $\alg{K}_2$, and we say that $\alg{T}_C'$ is incarnated by $(Q,C)$.  

Upon taking $F$-points, $T_C' = \alg{T}_C'(F)$ is a central extension of groups:
$$1 \rightarrow \alg{K}_2(F) \rightarrow T_C' \rightarrow T_C \rightarrow 1.$$
Identifying $T_C = \alg{T}_C(F)$ with $Y \otimes F^\times$, we can realize $T_C'$ explicitly by generators and relations.  $T_C'$ is generated by symbols $y'(u)$ for all $y \in Y$ and $u \in F^\times$, and symbols $z$ for all $z \in \alg{K}_2(F)$, together with the following relations:
\begin{enumerate}
\item
$[y_1'(u_1), y_2'(u_2)] = \{ u_1, u_2 \}^{B(y_1, y_2)}$;
\item
$y_1'(u) \cdot y_2'(u) = (y_1 + y_2)'(u) \cdot \{ u,u \}^{C(y_1, y_2)}$;
\item
$y'(u_1) \cdot y'(u_2) = y'(u_1 u_2) \cdot \{ u_1, u_2 \}^{Q(y)}$;
\item
those relations which imply that $\alg{K}_2(F)$ is embedded as a central subgroup.
\end{enumerate}
These relations follow from \cite[\S 3.9]{B-D}.  

\begin{remark}
Since the universal symbol $\{ \bullet, \bullet \}: F^\times \times F^\times \rightarrow \alg{K}_2(F)$ is skew-symmetric, the element $\{ u,u \}$ has order one or two in $\alg{K}_2(F)$.  Hence $\{ u,u \}^{C(y_1, y_2)}$ is well-defined for $C(y_1, y_2) \in \ZZ / 2 \ZZ$.
\end{remark}

A morphism $H\colon C_1 \rightarrow C_2$ in $\Cat{Bis}_Q$ yields a group homomorphism $\tau_H\colon T_{C_1}' \rightarrow T_{C_2}'$ given by $\tau_H(z) = z$ for all $z \in \alg{K}_2(F)$ and $\tau_H(y'(u)) = y'(u) \cdot \{ u,u \}^{H(y)}$ for all $y \in Y$, $u \in F^\times$.  This describes the incarnation functor (at the level of $F$-points):
$$\Cat{Inc}_{\alg{T}(F)}\colon \Cat{Bis}_Q \rightarrow \Cat{CExt}(\alg{T}, \alg{K}_2) \rightarrow \Cat{CExt}(T, \alg{K}_2(F)).$$

\subsection{Split semisimple simply-connected metaplectic groups}

Recall that $Y_{sc} = \ZZ[\Delta^\vee] \subset Y$ is the coroot lattice.  Let $f_{sc}\colon Y_{sc} \rightarrow Y$ denote the inclusion.  Corresponding to $f_{sc}$, there is a homomorphism of split reductive groups $\alg{G}_{sc} \rightarrow \alg{G}$, where $\alg{G}_{sc}$ is the simply-connected cover of the derived subgroup of $\alg{G}$.

When $Q\colon Y \rightarrow \ZZ$ is a {\em Weyl-invariant} quadratic form, its restriction $Q_{sc}$ to $Y_{sc}$ is a Weyl-invariant quadratic form on the coroot lattice.  Corresponding to $Q_{sc}$, there is a canonical (uniquely determined up to unique isomorphism) central extension $\alg{G}_{sc}'$ of $\alg{G}_{sc}$ by $\alg{K}_2$.  The $F$-points of this central extension can be defined by relations of Steinberg \cite{Ste}, for example.  The group $\alg{G}_{sc}'(F)$ is generated by elements $e_\alpha(x)$ for all $\alpha \in \Phi^\vee$ and $x \in F$.  

Defining $e_\alpha = e_\alpha(1)$, we obtain (via \cite[\S 11]{B-D}) a central extension $E_{Q_{sc}}$ of $Y_{sc}$ by $F^\times$, canonically associated to the quadratic form $Q$.  
\begin{construction}
The central extension $E_{Q_{sc}}$ admits the following presentation:  the generators are symbols $\{ e_i : i \in I \}$ and $\{ f : f \in F^\times \}$, and the relations are:
\begin{enumerate}
\item
$F^\times$ is contained in the center of $E_{Q_{sc}}$;
\item
For all $i,j \in I$, $[e_i, e_j] = (-1)^{B(\alpha_i^\vee, \alpha_j^\vee)}$.
\end{enumerate}
The sequence $F^\times \rightarrow E_{Q_{sc}} \rightarrow Y_{sc}$ is uniquely determined by the conditions that $f \in F^\times$ maps to the corresponding generator of $E_{Q_{sc}}$, and for all $i \in I$, $e_i \in E_{Q_{sc}}$ maps to $\alpha_i^\vee \in Y_{sc}$.
\end{construction}

\subsection{Incarnating split metaplectic groups}
\label{IncSC}
A bisector $C \in \Cat{Bis}_Q$ incarnates a central extension of $\alg{T}$ by $\alg{K}_2$, via the work of \cite{B-D}.  More data is required to incarnate a central extension of the split reductive group $\alg{G}$ by $\alg{K}_2$; this data includes a compatibility condition with the simply-connected group which we describe briefly below.    

On one hand, to each Weyl-invariant quadratic form $Q\colon Y \rightarrow \ZZ$, the restriction $Q_{sc}$ of $Q$ to $Y_{sc}$ is a Weyl-invariant quadratic form on $Y_{sc}$.  From $Q_{sc}$ we have obtained a canonical central extension $E_{Q_{sc}} \in \Cat{CExt}_{Q_{sc}}(Y_{sc}, F^\times)$.  

On the other hand, to each $E \in \Cat{CExt}_Q(Y, F^\times)$, one may pull back $E$ to get a central extension $f_{sc}^\ast E \in \Cat{CExt}_{Q_{sc}}(Y_{sc}, F^\times)$.  The two central extensions $E_{Q_{sc}}$ and $f_{sc}^\ast E$ are isomorphic (in the category $\Cat{CExt}(Y_{sc}, F^\times)$), since the isomorphism class of a central extension is determined by its commutator in this setting.  But this isomorphism is not unique.

\begin{definition}
Let $\Cat{CExt}_Q^{sc}(Y, F^\times)$ be the category whose objects are pairs $(E, \phi)$, where $E$ is a central extension of $Y$ by $F^\times$, and $\phi\colon E_{Q_{sc}} \rightarrow f_{sc}^\ast E$ is an isomorphism in the category $\Cat{CExt}_{Q_{sc}}(Y_{sc}, F^\times)$.  The morphisms of $\Cat{CExt}_Q^{sc}(Y, F^\times)$, from $(E_1, \phi_1)$ to $(E_2, \phi_2)$, are morphisms $\tau\colon E_1 \rightarrow E_2$ in $\Cat{CExt}_Q(Y, F^\times)$ such that $f_{sc}^\ast \tau \circ \phi_1 = \phi_2$.
\end{definition}

Just as we defined an essentially surjective functor from $\Cat{Bis}_Q$ to $\Cat{CExt}_Q(Y, F^\times)$, we can add a bit of extra data to define a category $\Cat{Bis}_Q^{sc}$ with an essentially surjective functor to $\Cat{CExt}_Q^{sc}(Y, F^\times)$.
\begin{definition}
Let $\Cat{Bis}_Q^{sc}$ be the category of pairs $(C, \eta)$, where $C$ is a bisector of $Q$, and $\eta\colon Y_{sc} \rightarrow F^\times$ is a group homomorphism.  Such a homomorphism $\eta$ is equivalent to a collection $\eta(i) = \eta(\alpha_i^\vee)$ of elements of $F^\times$, since $Y_{sc}$ is the free $\ZZ$-module spanned by $\Delta^\vee$.  A morphism of pairs, from $(C_1, \eta_1)$ to $(C_2, \eta_2)$ will mean a morphism $H\colon C_1 \rightarrow C_2$ in the category $\Cat{Bis}_Q$ satisfying the additional condition $\eta_2(i) / \eta_1(i) = (-1)^{H(\alpha_i^\vee)}$ for all $i \in I$.
\end{definition}

\begin{construction}
Define a functor from $\Cat{Bis}_Q^{sc}$ to $\Cat{CExt}_Q^{sc}(Y, F^\times)$ as follows:  to each object $(C,\eta)$ of $\Cat{Bis}_Q^{sc}$ we associate the extension $E_C$ of $Y$ by $F^\times$ as before.  We also associate the unique homomorphism $\phi_\eta$ from $E_{Q_{sc}}$ to $f_{sc}^\ast E_C$ sending $e_i$ to $(\alpha_i^\vee, \eta(i))$.

To a morphism $H\colon (C_1, \eta_1) \rightarrow (C_2, \eta_2)$, we associate the homomorphism $\phi_H\colon E_{C_1} \rightarrow E_{C_2}$ as before; this homomorphism intertwines the homomorphisms $\phi_{\eta_1}$ and $\phi_{\eta_2}$ from $E_{Q_{sc}}$ to $E_{C_1}$ and $E_{C_2}$.
\end{construction}
\proof
First, we check that there exists a unique homomorphism from $E_{Q_{sc}}$ to $f_{sc}^\ast E_C$ sending $e_i$ to $(\alpha_i^\vee, \eta(i))$.  This follows from the facts that $Y_{sc}$ is generated by the coroots, and the commutator relation in $E_{Q_{sc}}$ coincides with that of $E_C$, when pulled back to $Y_{sc}$.

Now, observe that for a morphism $H\colon (C_1, \eta_1) \rightarrow (C_2, \eta_2)$, we have
$$\phi_H(\alpha_i^\vee, \eta_1(i)) = (\alpha_i^\vee, \eta_1(i) (-1)^{H(\alpha_i^\vee)} ) = (\alpha_i^\vee, \eta_2(i)).$$
Hence $\phi_H$ intertwines the homomorphisms from $E_{Q_{sc}}$ as required.
\qed

The main theorem of \cite{B-D} provides an equivalence of categories from the groupoid $\Cat{CExt}(\alg{G}, \alg{K}_2)$ to $\bigsqcup_Q \Cat{CExt}_Q^{sc}(Y, F^\times)$ , where $Q$ ranges over all Weyl-invariant quadratic forms on $Y$.  Define
$$\Cat{BD}\colon \bigsqcup_Q \Cat{CExt}_Q^{sc}(Y, F^\times) \rightarrow \Cat{CExt}(\alg{G}, \alg{K}_2)$$
to be an inverse (well-defined up to natural equivalence) of the functor of Brylinski-Deligne.  Let $\Cat{Inc}_{\alg{G}}$ be the composite of $\Cat{BD}$ with the functor described in the previous construction:
$$\Cat{Inc}_{\alg{G}}\colon \Cat{Bis}_Q^{sc} \rightarrow \Cat{CExt}_Q^{sc}(Y, F^\times) \rightarrow \Cat{CExt}(\alg{G}, \alg{K}_2).$$

This functor will be called ``incarnation''.  It is the functorial construction of a central extension of $\alg{G}$ by $\alg{K}_2$, using the data $(C, \eta)$ of a bisector of a Weyl-invariant quadratic form, together with a homomorphism $\eta\colon Y_{sc} \rightarrow F^\times$.

\subsection{Fair bisectors}

Though we can incarnate central extensions of $\alg{G}$ by $\alg{K}_2$ using pairs $(C, \eta)$ as above, it will become useful to place another condition on bisectors.  

\begin{definition}
Let $Q$ be a Weyl-invariant quadratic form on $Y$.  A bisector $C$ of $Q$ is called {\em fair} if 
$$\forall y \in Y, \forall \alpha^\vee \in \Delta^\vee, Q(\alpha^\vee) \in 2 \ZZ \Rightarrow C(\alpha^\vee, y) = C(y, \alpha^\vee) = 0.$$
Let $\Cat{Bis}_Q^{f}$ denote the full subcategory of $\Cat{Bis}_Q$ consisting of fair bisectors.
\end{definition}

\begin{proposition}
Let $Q$ be a Weyl-invariant quadratic form on $Y$.  Then $Q$ possesses a fair bisector.
\end{proposition}
\proof
By Lemma \ref{BQLemma}, we find that whenever $Q(\alpha^\vee)$ is even, $B(\alpha^\vee, y)$ is even.

Now, let $\Delta_{odd}^\vee = \{ v_1, \ldots, v_\ell \}$ denote the set of simple coroots not contained in $2 Y$.  Then the set $\bar \Delta_{odd}^\vee$ of reductions (mod $2Y$) of these coroots is a $(\ZZ / 2 \ZZ)$-linearly independent subset of $Y / 2 Y$.  Thus it may be completed to an ordered basis $\{ v_1, \ldots, v_m \}$ of the $(\ZZ / 2 \ZZ)$-vector space $Y / 2 Y$.

Define a bilinear form $\bar C\colon Y / 2 Y \times Y / 2 Y \rightarrow \ZZ / 2 \ZZ$ by
$$\bar C(v_i, v_j) = \begin{cases} B(v_i, v_j) \mbox{ mod } 2, & \mbox{ if } i < j \\ Q(v_i) \mbox{ mod } 2, & \mbox{ if } i = j \\ 0, & \mbox{ if } i > j. \end{cases}$$
Then the lift $C\colon Y \times Y \rightarrow \ZZ / 2 \ZZ$ is a fair bisector of $Q$.
\qed

\begin{remark}
This construction is based on the cocycle used to construct covering groups at the beginning of \cite{Sav}.  Indeed, when the root datum is simply-connected and simply-laced, this construction coincides with the construction found there.
\end{remark}

\begin{remark}
The importance of fairness will become clear later, when constructing the L-group of a metaplectic group.
\end{remark}

\begin{definition}
Let $\Cat{Bis}_Q^{f,sc}$ denote the full subcategory of $\Cat{Bis}_Q^{sc}$, whose objects are those pairs $(C, \eta)$ in which $C$ is a fair bisector of $Q$.
\end{definition}

\begin{proposition}
The incarnation functor $\Cat{Inc}_{\alg{G}}\colon \Cat{Bis}_Q^{f, sc} \rightarrow \Cat{CExt}(\alg{G}, \alg{K}_2)$ is essentially surjective.
\end{proposition}
\proof
Since every Weyl-invariant quadratic form has a fair bisector, the inclusion $\Cat{Bis}_Q^{f, sc} \rightarrow \Cat{Bis}_Q^{sc}$ is essentially surjective.    

By \cite{B-D}, it suffices to prove that the functor from $\Cat{Bis}_Q^{sc}$ to $\Cat{CExt}_Q^{sc}(Y, F^\times)$ is essentially surjective.  We have already seen that the functor from $\Cat{Bis}_Q$ to $\Cat{CExt}_Q(Y, F^\times)$ is essentially surjective.  

Finally, the data of a morphism from $E_{Q_{sc}}$ to $f_{sc}^\ast E_C$ is precisely given by the data of a homomorphism $\eta\colon Y_{sc} \rightarrow F^\times$.
\qed

\begin{remark}
The previous proposition indicates that we can construct all of the central extensions studied by Brylinski-Deligne in \cite{B-D}, {\em up to isomorphism}, by these incarnation functors.  For working with nonsplit groups via descent, it will become important to include all morphisms in the category $\Cat{CExt}_Q(\alg{G}, \alg{K}_2)$, but this is beyond the scope of the present article.
\end{remark}

\section{Cocycle double-twists and the L-group}

\subsection{Setup}
Let $F$ now denote a {\em local} field, with $F \not \isom \CC$.  

\subsubsection{The metaplectic group and its dual group}

Let $\alg{G}$ be a split connected reductive group over $F$.  Let $\alg{T} \subset \alg{B} \subset \alg{G}$ be a split maximal torus contained in a Borel subgroup of $\alg{G}$, all defined over $F$.  Let $(Y,X)$ be the associated root datum, of type $(I, \cdot)$.  

Let $Q\colon Y \rightarrow \ZZ$ be a Weyl-invariant quadratic form on $Y$.  Let $n$ be a positive integer which divides $\# \mu(F)$.  Associated to $(Q,n)$, let $(\tilde Y, \tilde X)$ denote the resulting root datum of type $(I, \odot)$, and $\tilde \Psi^\vee = (\tilde X, \tilde Y)$ the dual root datum of type $(I, \checkodot)$.  Let $\alg{\tilde G}\LVee$ denote Lusztig's reductive group scheme over $\ZZ$ associated to the root datum $(I, \checkodot)$.  The group $\alg{\tilde G}\LVee$ comes equipped with a maximal torus $\alg{\tilde T}\LVee = \Spec(\ZZ[\tilde Y])$.  The center of $\alg{\tilde G}\LVee$ may be identified explicitly as the group scheme:
$$\alg{Z}(\alg{\tilde G}\LVee) = \Spec \left(  \frac{ \ZZ[\tilde Y] }{\langle \tilde \alpha_i^\vee : i \in I \rangle } \right) .$$

Let $(C,\eta)$ be an object of $\Cat{Bis}_Q^{f, sc}$, so $C\colon Y \times Y \rightarrow \ZZ / 2 \ZZ$ is a $\ZZ$-bilinear map satisfying the following conditions:
\begin{enumerate}
\item
$Q(y) \ident C(y,y)$, mod $2$.
\item
For all $y \in Y$ and $i \in I$, if $Q(\alpha_i^\vee) \in 2 \ZZ$ then $C(\alpha_i^\vee, y) = C(y, \alpha_i^\vee) = 0$.
\end{enumerate}

Let $\alg{G}'$ be the central extension of $\alg{G}$ by $\alg{K}_2$ over $F$, incarnated by $(C, \eta)$.  Let $\Hilb_n\colon \alg{K}_2(F) \rightarrow \mu_n(F)$ be the surjective homomorphism related to the (full) Hilbert symbol by $\Hilb_n = \Hilb^{\# \mu(F) / n}$

Pushing $\alg{K}_2(F) \rightarrow \alg{G}'(F) \rightarrow \alg{G}(F)$ forward via $\Hilb_n$ yields a central extension
$$1 \rightarrow \mu_n(F) \rightarrow \tilde G_F \rightarrow G_F \rightarrow 1.$$
We call $\tilde G_F$ a metaplectic cover of $G_F$.  As we assume $F$ is a local field, the group $\tilde G_F$ has a natural topology which makes it a locally compact group, and the above sequence is an exact sequence of topological groups.

\begin{remark}
Note that the ``dual group'' $\alg{\tilde G}\LVee$ depends only on $Q$ and $n$ and the root datum $(Y,X)$.  On the other hand, the metaplectic group $\tilde G_F$ depends also on the data $(C, \eta)$.  The additional data of the bisector $C$ of $Q$ does not appear in the construction of the dual group, but it will appear in the construction of the {\em L-group}.
\end{remark}

\subsubsection{Galois groups}

Let $\bar F$ be a separable algebraic closure of $F$, and let $\Gamma_F = \Gal(\bar F / F)$.  Then $\Gamma_F = \liminv Gal(L/F)$, as $L$ ranges over the finite Galois extensions of $F$.  We write $\OO(\Gal(L/F))$ for the ring of functions from $\Gal(L/F)$ to $\ZZ$.  We write $\OO(\Gamma_F)$ for $\limdir \OO(\Gal(L/F))$; it is a commutative Hopf algebra over $\ZZ$, and we write $\alg{\Gamma}_F = \Spec(\OO(\Gamma_F))$ for the resulting affine group scheme over $\ZZ$.  We identify the $\ZZ$-module $\OO(\Gamma_F) \otimes_\ZZ \OO(\Gamma_F)$ with $\OO(\Gamma_F \times \Gamma_F)$ in the unique way for which $[f \otimes g](\gamma_1, \gamma_2) = f(\gamma_1) \cdot g(\gamma_2)$; the comultiplication in the Hopf algebra $\OO(\Gamma_F)$ can then be given by
$$[\Delta f](\gamma_1, \gamma_2) = f(\gamma_1 \gamma_2)$$
for all $\gamma_1, \gamma_2 \in \Gamma_F$ and all $f \in \OO(\Gamma_F)$.  

Let $\widehat{F^\times}$ denote the inverse limit $\liminv F^\times / U$, where $U$ ranges over all open subgroups of $F^\times$ of finite index.  Let $\rec\colon \Gamma_F \rightarrow \widehat{F^\times}$ denote the reciprocity isomorphism of local class field theory, normalized in the nonarchimedean case so that a geometric Frobenius element maps to a uniformizing element.  Define a locally constant function $h\colon \Gamma_F \times \Gamma_F \rightarrow \mu_n(F)$ by
$$h(\gamma_1, \gamma_2) = \Hilb_n(\rec(\gamma_1), \rec(\gamma_2)).$$
We identify the 2-torsion of $\mu_n(F)$ with a subgroup of $\ZZ^\times$ in the unique way.  Thus $h(\gamma, \gamma) \in \ZZ^\times$ for all $\gamma \in \Gamma_F$.  Moreover, if $2m \in n \ZZ$, then $h(\gamma_1, \gamma_2)^m \in \ZZ^\times$ for all $\gamma_1, \gamma_2 \in \Gamma_F$.

\subsubsection{The na\"ive L-group}

The product $\alg{\Gamma}_F \times_\ZZ \alg{\tilde G}\LVee$ of group schemes is a group scheme over $\ZZ$ which we might call the {\em na\"ive L-group}, since we will find a more intricate construction necessary later.  In both Hopf algebras $\OO(\alg{\tilde G}\LVee)$ and $\OO(\Gamma_F)$, we write $1$ for the unit and $\epsilon$ for the counit -- confusion should not arise.

The na\"ive L-group $\alg{\Gamma}_F \times \alg{\tilde G}\LVee$ is a group scheme over $\ZZ$, and its Hopf algebra is $\OO(\Gamma_F) \otimes \OO(\alg{\tilde G}\LVee)$.  All tensor products will be understood over $\ZZ$.  In this Hopf algebra, multiplication and comultiplication are given componentwise, the unit is $1 \otimes 1$, and the counit is $\epsilon \otimes \epsilon$.  We apply the natural commutativity and associativity isomorphisms, and the natural isomorphism from $\OO(\Gamma_F) \otimes \OO(\Gamma_F)$ to $\OO(\Gamma_F \times \Gamma_F)$ whenever convenient -- for example, we consider the comultiplication in $\OO(\Gamma_F) \otimes \OO(\alg{\tilde G}\LVee)$ to take values in $\OO(\Gamma_F \times \Gamma_F) \otimes \OO(\alg{\tilde G}\LVee) \otimes \OO(\alg{\tilde G}\LVee)$.  Accordingly, we can write 
$$\Delta(f \otimes a^\dual) = \Delta(f) \otimes \Delta(a^\dual) = \sum_{a_1, a_2 \in \BB} m_{a_1, a_2}^a \Delta(f) \otimes (a_1^\dual \otimes a_2^\dual),$$ 
for all $f \in \OO(\Gamma_F)$ and all $a \in \BB$, to describe the comultiplication in $\OO(\Gamma_F) \otimes \OO(\alg{\tilde G}\LVee)$.

The na\"ive L-group $\alg{\Gamma}_F \times \alg{\tilde G}\LVee$ does not seem to suffice, when attempting to generalize Langlands conjectures to metaplectic groups.  In order to construct an appropriate L-group, we apply a {\em double-twist}.  Using a pair of cocycles $(\tau, \chi)$, we ``double-twist'' the Hopf algebra $\OO(\Gamma_F) \otimes \OO(\alg{\tilde G}\LVee)$ to define a Hopf algebra
$$\OO({}^L \alg{\tilde G}\LVee) = \OO(\Gamma_F) \otimes_{\tau, \chi} \OO(\alg{\tilde G}\LVee).$$
The double-twist is usually called the {\em cocycle bicrossproduct} in the literature on Hopf algebras and quantum groups; we refer to Majid's book \cite{Maj} for a treatment from this perspective.

The algebra $\OO({}^L \alg{\tilde G}\LVee)$ (constructed later in this section) will be a commutative Hopf algebra over $\ZZ$, with associated group scheme ${}^L \alg{\tilde G}\LVee$ over $\ZZ$.  When $n$ is odd, ${}^L \alg{\tilde G}\LVee$ is {\em equal} to the na\"ive L-group.  But interestingly, when $n = 2$ the group scheme ${}^L \alg{\tilde G}\LVee$ is isomorphic to the na\"ive L-group $\alg{\Gamma}_F \times \alg{\tilde G}\LVee$, but only after extending the base ring from $\ZZ$ to $\ZZ[i]$; moreover, the isomorphism is not unique, but relies often on a choice of additive character of $F$ and other (more combinatorial) choices.

\subsection{The first twist}

For any $\ZZ$-linear map $\tau\from \OO(\alg{\tilde G}\LVee) \rightarrow \OO(\Gamma_F \times \Gamma_F)$, define a $\ZZ$-linear map
$\Delta_\tau \from \OO(\Gamma_F) \otimes \OO(\alg{\tilde G}\LVee) \rightarrow \OO(\Gamma_F \times \Gamma_F) \otimes \left( \OO(\alg{\tilde G}\LVee) \otimes \OO(\alg{\tilde G}\LVee) \right)$ by
$$\Delta_\tau(h \otimes b^\dual) = \sum_{b_1, b_2, b_3 \in \BB} m_{b_1, b_2, b_3}^b \left( \Delta(h) \cdot \tau(b_3^\dual) \right) \otimes \left( b_1^\dual \otimes b_2^\dual \right).$$
We say that $\Delta_\tau$ is obtained from twisting the comultiplication by $\tau$.

\begin{definition}[Cf.\ {\cite[Proposition 6.3.8]{Maj} }]
Let $\tau$ be a $\ZZ$-linear map as above.  We say that $\tau$ is a {\em (1,2)-cocycle} if for all $b \in \BB$, the following three conditions hold:
\begin{enumerate}
\item[(OT1)]
$[\epsilon \otimes \Id](\tau(b^\dual)) = 1_{\OO(\Gamma_F)} \cdot \epsilon(b^\dual) = [\Id \otimes \epsilon](\tau(b^\dual))$;
\item[(OT2)]
$\sum_{b_1, b_2 \in \BB} m_{b_1, b_2}^b \tau(b_2^\dual ) \otimes b_1^\dual = \sum_{b_1, b_2 \in \BB} m_{b_1, b_2}^b \tau(b_1^\dual) \otimes b_2^\dual$;
\item[(OT3)]
$\sum_{b_1, b_2 \in \BB} m_{b_1, b_2}^b (1 \otimes \tau(b_1^\dual)) \cdot [\Id \otimes \Delta](\tau(b_2^\dual)) = \\ \sum_{b_1, b_2 \in \BB} m_{b_1, b_2}^b (\tau(b_1^\dual) \otimes 1) \cdot [\Delta \otimes \Id](\tau(b_2^\dual))$.
\end{enumerate}
\end{definition}

\begin{remark}  We use the term ``(1,2)-cocycle'' to remind the reader that the domain has one factor $\OO(\alg{\tilde G}\LVee)$, and the range is the tensor product of two factors $\OO(\Gamma_F) \otimes \OO(\Gamma_F)$.  Later we will encounter (2,1)-cocycles in which this is reversed.  The conditions labelled ``OT'' should remind the reader of ``One-Two''.  Under additional commutativity and cocommutativity assumptions on Hopf algebras, the (1,2)-cocycles and (2,1)-cocycles are special cases of $(p,q)$-cocycles in a double-complex (cf.\ \cite[\S 2]{Hof}).    
\end{remark}

\begin{proposition}
If $\tau$ is a (1,2)-cocycle, then $\Delta_\tau$ makes $\OO(\Gamma_F) \otimes \OO(\alg{\tilde G}\LVee)$ into a coassociative coalgebra over $\ZZ$ with the same counit $\epsilon \otimes \epsilon$ as before.
\end{proposition}
\proof
The proof from \cite[Proposition 6.3.8]{Maj} carries through.
\qed

\begin{proposition}
Suppose that $\sigma\colon \tilde Y \times \Gamma_F \times \Gamma_F \rightarrow A$ is a function satisfying the following five properties for all $\gamma_1, \gamma_2, \gamma_3, \gamma \in \Gamma_F$ and all $y, y_1, y_2 \in \tilde Y$:
\begin{enumerate}
\item
$\sigma(\tilde \alpha_i^\vee; \gamma_1, \gamma_2) = 1$ for all $i \in I$;
\item
$\sigma(y_1 + y_2; \gamma_1, \gamma_2) = \sigma(y_1; \gamma_1, \gamma_2) \cdot \sigma(y_2; \gamma_1, \gamma_2)$ and $\tau(0; \gamma_1, \gamma_2) = 1$;
\item
There exists a subgroup $\Gamma_L$ of finite index in $\Gamma_F$ such that $\sigma(y; \gamma_1 \delta_1, \gamma_2 \delta_2) = \sigma(y; \gamma_1, \gamma_2)$ for all $\delta_1, \delta_2 \in \Gamma_L$;
\item
$\sigma(y; \gamma, 1) = \sigma(y; 1, \gamma) = 1$;
\item
$\sigma(y; \gamma_1, \gamma_2 \gamma_3) \sigma(y; \gamma_2, \gamma_3) = \sigma(y; \gamma_1 \gamma_2, \gamma_3) \sigma(y; \gamma_1, \gamma_2)$.
\end{enumerate}
Define a toral $\ZZ$-linear map $\tau\colon \OO(\alg{\tilde G}\LVee) \rightarrow \OO(\Gamma_F \times \Gamma_F)$ by
$$\tau(1_y^\circ)(\gamma_1, \gamma_2) = \sigma(y; \gamma_1, \gamma_2),$$
Then $\tau$ is a (1,2)-cocycle.
\end{proposition}
\proof
First, observe that by (2), for any fixed $\gamma_1, \gamma_2$, the function $y \mapsto \sigma(y; \gamma_1, \gamma_2)$ is a homomorphism from $\tilde Y$ to $\ZZ^\times$.  On the other hand, for any fixed $y \in \tilde Y$, we write $\sigma_y(\gamma_1, \gamma_2)$ for $\sigma(y; \gamma_1, \gamma_2)$ and observe that the function $\sigma_y$ is a 2-cocycle from $\Gamma_F$ to $\ZZ^\times$.
  
We now check the (1,2)-cocycle conditions of $\tau$ directly.  For condition (OT1), we find that $\tau(b^\dual) = \epsilon(b^\dual) = 0$ unless $b^\dual = 1_y^\dual$ for some $y \in \tilde Y$ (since $\tau$ is toral).  When $b^\dual = 1_y^\dual$, we find that for all $\gamma \in \Gamma_F$, $1_{\OO(\Gamma_F)}(\gamma) \cdot \epsilon(b^\dual) = 1$ and
$$[\epsilon \otimes \Id](\tau(b^\dual))(\gamma) = \sigma(y; 1,\gamma) = 1 = \sigma(y; \gamma, 1) = [\Id \otimes \epsilon](\tau(b^\dual))(\gamma).$$
Hence (OT1) holds.

For (OT2), suppose that $b \in \BB$, and $\tau(b_2^\dual) \neq 0$.  Then $b_2 = 1_{y_2}$ for some $y_2 \in \tilde Y$, and $m_{b_1, b_2}^b = 0$ unless $b = b_1$ and $y_2 = \rt{b}$.  Thus we find that
$$\sum_{b_1, b_2 \in \BB} m_{b_1, b_2}^b \tau(b_2^\dual) \otimes b_1^\dual = \sigma_{y_2} \otimes b^\dual.$$
Similarly, if $b \in \BB$ and $\tau(b_1^\dual) \neq 0$, then $b_1 = 1_{y_1}$ for some $y_1 \in \tilde Y$.  $m_{b_1, b_2}^b = 0$ unless $b = b_2$ and $y_1 = \lt{b}$.  In this case, we find that
$$\sum_{b_1, b_2 \in \BB} m_{b_1, b_2}^b \tau(b_1^\dual) \otimes b_2^\dual = \sigma_{y_1} \otimes b^\dual.$$
We find that $\sigma_{y_1} \otimes b^\dual = \sigma_{y_2} \otimes b^\dual$, since $y_1 - y_2 = \nt{b}$ is a sum of coroots, and so $\sigma_{y_1 - y_2}$ is trivial by (1) and (2).  This proves (OT2).

Third, suppose that $b \in \BB$, and $\tau(b_1^\dual) \neq 0$ and $\tau(b_2^\dual) \neq 0$.  Then $b_1 = 1_{y_1}$ and $b_2 = 1_{y_2}$ for some $y_1, y_2 \in Y$.  Then $m_{b_1, b_2}^b = 0$ unless $b = b_1 = b_2 = 1_y$ for $y = y_1 = y_2 \in \tilde Y$.  We find that
$$\sum_{b_1, b_2 \in \BB} m_{b_1, b_2}^b (1 \otimes \tau(b_1^\dual)) \cdot [\Id \otimes \Delta](\tau(b_2^\dual)) = (1 \otimes \sigma_y) \cdot [\Id \otimes \Delta](\sigma_y).$$
This is an element of $\OO(\Gamma_F \times \Gamma_F \times \Gamma_F)$ whose evaluation at $(\gamma_1, \gamma_2, \gamma_3)$ is given by
$$\left[  (1 \otimes \sigma_y) \cdot [\Id \otimes \Delta](\sigma_y) \right] (\gamma_1, \gamma_2, \gamma_3) = \sigma(y; \gamma_2, \gamma_3) \cdot \sigma(y; \gamma_1, \gamma_2 \gamma_3).$$
Similarly, we find that
$$\sum_{b_1, b_2 \in \BB} m_{b_1, b_2}^b (\tau(b_1^\dual) \otimes 1) \cdot [\Delta \otimes \Id](\tau(b_2^\dual)) = (\sigma_y \otimes 1) \cdot [\Delta \otimes \Id](\sigma_y),$$
whose evaluation at $(\gamma_1, \gamma_2, \gamma_3)$ is given by
$$\left[ (\sigma_y \otimes 1) \cdot [\Delta \otimes \Id](\sigma_y) \right] (\gamma_1, \gamma_2, \gamma_3) = \sigma(y; \gamma_1, \gamma_2) \cdot  \sigma(y; \gamma_1 \gamma_2, \gamma_3).$$
Now (5) implies (OT3).
\qed

\begin{remark}
There is a natural bijection between functions $\sigma$ satisfying the five conditions of the above proposition and morphisms of schemes over $\ZZ$ from $(\alg{\Gamma}_F \times \alg{\Gamma}_F)$ to $\alg{Z}(\alg{\tilde G}\LVee)$, satisfying the usual conditions for a (normalized) 2-cocycle.  Thus the ``first twist'' of this section is nothing other than the usual twisting of a direct product of groups by a group-theoretic 2-cocycle.  Our reasons for working with Hopf algebras instead of groups will become clear after the next section.
\end{remark}

Our (1,2)-cocycles arise from the following construction.
\begin{construction}
Recall that $Q\colon Y \rightarrow \ZZ$ is an integer-valued, Weyl-invariant quadratic form.  Moreover, for all $y \in \tilde Y$, we have $2 Q(y) = B(y,y) \in n \ZZ$.  Hence $h(\gamma_1, \gamma_2)^{Q(y)} \in \{ \pm 1 \}$ for all $\gamma_1, \gamma_2 \in \Gamma_F$.  

The function
$$\sigma(y; \gamma_1, \gamma_2) = h(\gamma_1, \gamma_2)^{Q(y)}$$
satisfies the five conditions of the previous proposition, and hence defines a (1,2)-cocycle by
$$\tau(1_y^\dual)(\gamma_1, \gamma_2) = h(\gamma_1, \gamma_2)^{Q(y)}.$$
\end{construction}
\proof
Since $Q(\tilde \alpha_i^\vee) \in n \ZZ$ for all $i \in I$, we find that $\sigma$ gives a homomorphism from $\tilde Y$ to $\ZZ^\times$ which is trivial on coroots, for any fixed $\gamma_1, \gamma_2$.  This demonstrates Properties 1 and 2.  Local constancy of $h$ implies Property 3.  Since $h$ is bimultiplicative, Properties 4 and 5 hold as well.
\qed

In what follows, $\tau$ will always denote the (1,2)-cocycle constructed above.  Let $\OO(\Gamma_F) \otimes_\tau \OO(\alg{\tilde G}\LVee)$ denote the coalgebra over $\ZZ$ obtained by replacing the comultiplication of $\OO(\Gamma_F) \otimes \OO(\alg{\tilde G}\LVee)$ by $\Delta_\tau$.  

\subsection{The second twist}
In the previous section, we discussed twisting the coalgebra structure on the Hopf algebra $\OO(\Gamma_F) \otimes \OO(\alg{G}\LVee)$ by a (1,2)-cocycle.  Such a twist could also have been seen from the group-theoretic standpoint, with a $\alg{Z}(\alg{\tilde G}\LVee)$-valued 2-cocycle on $\alg{\Gamma}_F$.    

In this section, we study a twist which can {\em only} be seen properly from the Hopf algebraic perspective.  We twist the algebra structure of $\OO(\Gamma_F) \otimes \OO(\alg{\tilde G}\LVee)$ instead of the coalgebra structure.
\begin{remark}
If we worked over a field instead of over $\ZZ$, this new twist corresponds to a twist of the Tannakian category of representations of the na\"ive L-group.  Closely related methods of twisting the commutativity constraint in the Tannakian category are described in Reich's thesis \cite{Rei} and in another context by Lusztig in \cite{Lus1}; but the twists found there do not account for ramification -- they twist the dual group $\alg{\tilde G}\LVee$ without incorporating the Galois group $\alg{\Gamma}_F$.
\end{remark}

For any $\ZZ$-linear map $\chi\colon \OO(\alg{\tilde G}\LVee) \otimes \OO(\alg{\tilde G}\LVee) \rightarrow \OO(\Gamma_F)$, define a $\ZZ$-linear map 
$$\cdot_\chi\colon (\OO(\Gamma_F) \otimes \OO(\alg{\tilde G}\LVee)) \otimes (\OO(\Gamma_F) \otimes \OO(\alg{\tilde G}\LVee)) \rightarrow \OO(\Gamma_F) \otimes \OO(\alg{\tilde G}\LVee), \mbox{ by }$$
$$(f \otimes a^\dual) \cdot_\chi (h \otimes b^\dual) = \sum_{a_1, a_2, b_1, b_2 \in \BB} m_{a_1, a_2}^a m_{b_1, b_2}^b \left( f \cdot h \cdot \chi(a_1^\dual, b_1^\dual) \right) \otimes a_2^\dual b_2^\dual.$$

\begin{definition}
Let $\chi$ be a $\ZZ$-linear map as above.  We say that $\chi$ is a (2,1)-cocycle if for all $a, b, c \in \BB$, the following three conditions hold.
\begin{enumerate}
\item[(TO1)]
$\chi(a^\dual, 1) = \epsilon(a^\dual) \cdot 1 = \chi(1, a^\dual)$.
\item[(TO2)]
$\sum m_{a_1, a_2}^a m_{b_1, b_2}^b \epsilon(a_1^\dual b_1^\dual) \chi(a_2^\dual, b_2^\dual) = \sum m_{a_1, a_2}^a m_{b_1, b_2}^b \epsilon(a_2^\dual b_2^\dual) \chi(a_1^\dual, b_1^\dual) $
\item[(TO3)]
$\sum m_{b_1, b_2}^b m_{c_1, c_2}^c \chi(b_1^\dual, c_1^\dual) \cdot \chi(a^\dual, b_2^\dual c_2^\dual) = \sum m_{a_1, a_2}^a m_{b_1, b_2}^b \chi(a_1^\dual, b_1^\dual) \chi(a_2^\dual b_2^\dual, c^\dual).$
\end{enumerate}
We say that $\chi$ is {\em symmetric} if $\chi(a^\dual, b^\dual) = \chi(b^\dual, a^\dual)$ for all $a,b \in \BB$.
\end{definition}

\begin{proposition}[Cf.\ {\cite[Lemma 4.5 and Example 4.10]{BCM} }]
If $\chi$ is a (2,1)-cocycle, the $\cdot_\chi$ makes $\OO(\Gamma_F) \otimes \OO(\alg{\tilde G}\LVee)$ an associative unital algebra with unit element $1 \otimes 1$.  If, in addition, $\chi$ is symmetric, then this algebra is commutative.
\end{proposition}
\proof
The proof from \cite{BCM} carries through.
\qed

Our (2,1)-cocycles will come from the following construction.
\begin{construction}
Recall that $C\colon Y \times Y \rightarrow \ZZ / 2 \ZZ$ is a bilinear form (it is a fair bisector of $Q$).  Define a toral $\ZZ$-linear map $\chi\colon  \OO(\alg{\tilde G}\LVee) \otimes \OO(\alg{\tilde G}\LVee) \rightarrow \OO(\Gamma_F)$ by
$$\chi(1_{y_1}^\dual, 1_{y_2}^\dual)(\gamma) = h(\gamma, \gamma)^{C(y_1, y_2)},$$
for all $y_1, y_2 \in \tilde Y$ and $\gamma \in \Gamma$.  Then $\chi$ is a symmetric (2,1)-cocycle.  
\end{construction}
\proof
Define $\eta\colon \Gamma_F \rightarrow \{ \pm 1 \}$ by $\eta(\gamma) = h(\gamma, \gamma)$, in what follows.

For (TO1), observe that $\chi(a^\dual, 1) = \epsilon(a^\dual) \cdot 1 = \chi(1, a^\dual) = 0$ unless $a^\dual = 1_y^\dual$ for some $y \in Y$.  On the other hand, if $a^\dual = 1_y^\dual$, then $\chi(a^\dual, 1) = \eta^{C(y, 0)} = 1$ and $\chi(1, a^\dual) = \eta^{C(0,y)} = 1$, and $\epsilon(1_y^\dual) \cdot 1 = 1$.

For (TO2), note that $\epsilon(a_1 b_1) = \sum_{c_1} \hat m_{c_1}^{a_1, b_1} \epsilon(c_1^\dual)$.  Furthermore, $\epsilon(c_1^\dual) = 0$ unless $c_1^\dual = 1_{y_1}^\dual$ for some $y_1 \in \tilde Y$ and $\hat m_{c_1}^{a_1, b_1} = 0$ unless $a_1^\dual = 1_{u_1}^\dual$ and $b_1^\dual = 1_{v_1}^\dual$ for some $u_1, v_1 \in \tilde Y$ such that $u_1 + v_1 = y_1$.  Also, $\chi(a_2^\dual, b_2^\dual) = 0$ unless $a_2^\dual = 1_{u_2}^\dual$ and $b_2^\dual = 1_{v_2}^\dual$ for some $u_2, v_2 \in \tilde Y$.  But in this case, $m_{a_1, a_2}^a$ and $m_{b_1, b_2}^b$ vanish unless $u_1 = u_2$ and $v_1 = v_2$.  In this way, we find that the terms on the left and right sides of (TO2) vanish except when $a_1^\dual = a_2^\dual = a^\dual = 1_u^\dual$, and $b_1^\dual = b_2^\dual = b^\dual = 1_v^\dual$ for elements $u,v \in \tilde Y$.  The only nonvanishing term of (TO2), on the left and on the right side, is $\eta^{C(u,v)}$.

For (TO3), let $L(a,b,c)$ and $R(a,b,c)$ denote the left and right sides:
\begin{eqnarray*}
L(a,b,c) & = & \sum m_{b_1, b_2}^b m_{c_1, c_2}^c \chi(b_1^\dual, c_1^\dual) \cdot \chi(a^\dual, b_2^\dual c_2^\dual), \\
R(a,b,c) & = & \sum m_{a_1, a_2}^a m_{b_1, b_2}^b  \chi(a_1^\dual, b_1^\dual) \cdot \chi(a_2^\dual b_2^\dual, c^\dual).
\end{eqnarray*}
An argument as before demonstrates that all terms vanish except those where $a_1^\dual = a_2^\dual = a^\dual = 1_u^\dual$, $b_1^\dual = b_2^\dual = b^\dual = 1_v^\dual$, and $c_1^\dual = c_2^\dual = c^\dual = 1_w^\dual$, for some $u,v,w \in \tilde Y$.  For such $u,v,w$, we find that
$$L(a,b,c) = \eta^{C(v,w)} \cdot \eta^{C(u, v + w)}, \quad R(a,b,c) = \eta^{C(u,v)} \cdot \eta^{C(u + v, w)}.$$
The condition (TO3) follows immediately from $\ZZ$-bilinearity of $C$.  We have verified that $\chi$ is a (2,1)-cocycle.

For symmetry, we observe that when $y_1, y_2 \in \tilde Y$, we have
$$C(y_1, y_2) + C(y_2, y_1) \ident B(y_1, y_2), \mbox{ mod } 2.$$
Either $n$ is odd, in which case $h(\gamma, \gamma)$ is identically $1$ (and $\chi = \epsilon \otimes \epsilon$ is trivial), or $n$ is even, and $B(y_1, y_2) \in n \ZZ$ is even and thus $C(y_1, y_2) = C(y_2, y_1) \in \ZZ / 2 \ZZ$.  In either case, $h(\gamma, \gamma)^{C(y_1, y_2)} = h(\gamma, \gamma)^{C(y_2, y_1)}$ so $\chi$ is a symmetric (2,1)-cocycle.
\qed

\subsection{Compatibility}
We have discussed two ways to modify the Hopf algebra $\OO(\Gamma_F) \otimes \OO(\alg{\tilde G}\LVee)$.  We may modify the coalgebra structure using a (1,2)-cocycle $\tau\colon \OO(\alg{\tilde G}\LVee) \rightarrow \OO(\Gamma_F) \otimes \OO(\Gamma_F)$.  We may modify the algebra structure using a (2,1)-cocycle $\chi\colon \OO(\alg{\tilde G}\LVee) \otimes \OO(\alg{\tilde G}\LVee) \rightarrow \OO(\Gamma_F)$.  We write
$$\OO(\Gamma_F) \otimes_{\tau, \chi} \OO(\alg{\tilde G}\LVee)$$
for the simultaneous algebra and coalgebra obtained by applying the twist $\tau$ to the coalgebra and $\chi$ to the algebra structure.  Additional conditions are required for this to be a Hopf algebra.

\begin{definition}
We say that the pair $(\tau, \chi)$ is {\em compatible} if the following two conditions hold.
\begin{enumerate}
\item[(Com1)]
For all $a,b \in \BB$, the following equality holds in $\OO(\Gamma_F) \otimes \OO(\alg{\tilde G}\LVee)$:
$$\sum_{a_1, a_2} \sum_{b_1, b_2} m_{a_1, a_2}^a m_{b_1, b_2}^b \chi(a_1^\dual, b_1^\dual) \otimes a_2^\dual b_2^\dual = \sum_{a_1, a_2} \sum_{b_1, b_2} m_{a_1, a_2}^a m_{b_1, b_2}^b \chi(a_2^\dual, b_2^\dual) \otimes a_1^\dual b_1^\dual.$$
\item[(Com2)]
For all $a,b \in \BB$, the following equality holds in $\OO(\Gamma_F \times \Gamma_F)$:
\begin{eqnarray*}
\sum_{a_1, a_2, a_3} \sum_{b_1, b_2, b_3} m_{a_1, a_2, a_3}^a m_{b_1, b_2, b_3}^b \tau(a_3^\dual) \cdot \tau(b_3^\dual) \cdot \left( \chi(a_1^\dual, b_1^\dual) \otimes \chi(a_2^\dual, b_2^\dual) \right) = \\ \sum_{a_1, a_2} \sum_{b_1, b_2} m_{a_1, a_2}^a m_{b_1, b_2}^b  \tau(a_1^\dual b_1^\dual) \cdot \Delta \left( \chi(a_2^\dual, b_2^\dual) \right).
\end{eqnarray*}
\end{enumerate}
\end{definition}

\begin{definition}
We say that $\chi\colon \OO(\alg{\tilde G}\LVee) \otimes \OO(\alg{\tilde G}\LVee) \rightarrow \OO(\Gamma_F)$ is {\em convolution-invertible} if there exists $\lambda\colon \OO(\alg{\tilde G}\LVee) \otimes \OO(\alg{\tilde G}\LVee) \rightarrow \OO(\Gamma_F)$ such that for all $a,b \in \BB$,
$$\sum_{a_1, a_2} \sum_{b_1, b_2} m_{a_1, a_2}^a m_{b_1, b_2}^b \chi(a_1^\dual, b_1^\dual) \cdot \lambda(a_2^\dual, b_2^\dual) = \epsilon(a^\dual) \epsilon(b^\dual).$$
We say that $\chi$ is {\em involutive}, if the above equality holds with $\lambda = \chi$.
\end{definition}

\begin{thm}[Cf.\ {\cite[Theorem 2.9]{Maj3}}]
If $(\tau, \chi)$ is a compatible pair, then $\OO(\Gamma_F) \otimes_{\tau, \chi} \OO(\alg{\tilde G}\LVee)$ is a $\ZZ$-bialgebra with respect to the multiplication $\cdot_\chi$, unit $1 \otimes 1$, comultiplication $\Delta_\tau$, and counit $\epsilon \otimes \epsilon$.  If, moreover, the cocycle $\chi$ is convolution-invertible, then this bialgebra is a Hopf algebra over $\ZZ$.
\end{thm}
\proof
Knowing from the previous two sections that $\OO(\Gamma_F) \otimes_{\tau, \chi} \OO(\alg{\tilde G}\LVee)$ is an associative unital algebra, and coassociative counital coalgebra, we must check that the comultiplication is an algebra homomorphism; in other words, we must verify that for all $f,g \in \OO(\Gamma_F)$ and all $a,b \in \BB$,
\begin{equation}
\label{Compat}
\Delta_\tau \left( (f \otimes a^\dual) \cdot_\chi (g \otimes b^\dual) \right) = \Delta_\tau(f \otimes a^\dual) \cdot_\chi \Delta_\tau(g \otimes b^\dual).
\end{equation}
The left side of \ref{Compat} is the following function of $a,b \in \BB$, $f,g \in \OO(\Gamma_F)$:
$$L(a,b) = \sum m_{a_1, a_2, a_3, a_4}^a m_{b_1, b_2, b_3, b_4}^b \Delta f \cdot \Delta g \cdot \Delta(\chi(a_4^\dual, b_4^\dual)) \cdot \tau(a_3^\dual \cdot b_3^\dual) \otimes (a_1^\dual b_1^\dual \otimes a_2^\dual b_2^\dual).$$
The right side of \ref{Compat} is the following function of $a,b \in \BB$, $f,g \in \OO(\Gamma_F)$:
\begin{eqnarray*} R(a,b) & = & \sum m_{a_1, a_2, a_3, a_4, a_5}^a m_{b_1, b_2, b_3, b_4, b_5}^b \Delta f \cdot \Delta g \cdot \tau(a_5^\dual) \cdot \tau(b_5^\dual) \\ & & \cdot \left( \chi(a_2^\dual, b_2^\dual) \otimes \chi(a_4^\dual,b_4^\dual) \right) \otimes (a_1^\dual b_1^\dual \otimes a_3^\dual b_3^\dual).
\end{eqnarray*} 
\qed

Using (Com1) and coassociativity in the Hopf algebra $\OO(\alg{\tilde G}\LVee)$, the expression $\chi(a_2^\dual,b_2^\dual) \otimes a_3^\dual b_3^\dual$ may be replaced by $\chi(a_3^\dual,b_3^\dual) \otimes a_2^\dual b_2^\dual$.  The right side becomes:
\begin{eqnarray*} R(a,b) & = & \sum m_{a_1, a_2, a_3, a_4, a_5}^a m_{b_1, b_2, b_3, b_4 b_5}^b \Delta f \cdot \Delta g \cdot \tau(a_5^\dual) \cdot \tau(b_5^\dual) \\ & & \cdot \left( \chi(a_3^\dual, b_3^\dual) \otimes \chi(a_4^\dual,b_4^\dual) \right) \otimes (a_1^\dual b_1^\dual \otimes a_2^\dual b_2^\dual).
\end{eqnarray*} 
Using (Com2) and coassociativity, the expression $\tau(a_5^\dual) \cdot \tau(b_5^\dual) \cdot (\chi(a_3^\dual, b_3^\dual) \otimes \chi(a_4^\dual, b_4^\dual))$ may be replaced by $\tau(a_3^\dual b_3^\dual) \cdot \Delta( \chi(a_4^\dual, b_4^\dual) )$ (changing summation and structure constants as appropriate).  We find that 
\begin{eqnarray*} R(a,b) & = & \sum m_{a_1, a_2, a_3, a_4}^a m_{b_1, b_2, b_3, b_4}^b \Delta f \cdot \Delta g \cdot \tau(a_3^\dual b_3^\dual) \\ & & \cdot \Delta( \chi(a_4^\dual, b_4^\dual) ) \otimes (a_1^\dual b_1^\dual \otimes a_2^\dual b_2^\dual).
\end{eqnarray*} 
The commutativity of multiplication in $\OO(\Gamma_F)$ yields $L(a,b) = R(a,b)$.  

To check that $\OO(\Gamma_F) \otimes_{\tau, \chi} \OO(\alg{\tilde G}\LVee )$ is a bialgebra, we must furthermore check compatibility among unit and comultiplication, counit and comultiplication, unit and counit.  The unit and counit are the same as the untwisted algebra, and so their compatibility is automatic.

For unit and comultiplication, observe that $\Delta_\tau(1 \otimes 1) = \Delta(1) \cdot \tau(1) \otimes (1 \otimes 1) = \Delta(1) \otimes (1 \otimes 1) = (1 \otimes 1) \otimes 1 \otimes 1$.  

For counit and multiplication, observe that 
\begin{eqnarray*}
[\epsilon \otimes \epsilon]( (f \otimes a^\dual) \cdot_\chi (g \otimes b^\dual) ) & = & [\epsilon \otimes \epsilon] \left( \sum_{a_1, a_2} \sum_{b_1, b_2}  f g \chi(a_1^\dual, b_1^\dual) \otimes a_2^\dual b_2^\dual \right) \\
& = & \sum_{a_1, a_2} \sum_{b_1, b_2} \epsilon(f) \epsilon(g) \epsilon(\chi(a_1^\dual, b_1^\dual)) \otimes \epsilon(a_2^\dual b_2^\dual ) \\
& = & \epsilon(f) \epsilon(g) \otimes \epsilon(a^\dual) \epsilon(b^\dual),
\end{eqnarray*}
where we observe that $\epsilon(\chi(a_1^\dual, b_1^\dual)) = 1$ if $a_1^\dual, b_1^\dual \in \BB_{\alg{\tilde T}}$ and $\epsilon(\chi(a_1^\dual, b_1^\dual)) = 0$ otherwise.  This finishes the proof that $\OO(\Gamma_F) \otimes_{\tau, \chi} \OO(\alg{\tilde G}\LVee)$ is a bialgebra.

For the antipode, we note that Schauenberg (in \cite[Theorem A.2]{Sch}) constructs an antipode $S$ when we work with Hopf algebras over the fraction field $\QQ$ of $\ZZ$.  His construction involves only the existing Hopf algebra structures and the convolution-inverse of $\chi$.  Thus Schauenburg's antipode is $\ZZ$-linear in our setting.  Since all our $\ZZ$-modules are free, the commutativity of the anitpode diagrams after tensoring with $\QQ$ implies the commutativity of the diagrams in the original modules over $\ZZ$.  Hence, when $\chi$ is convolution-invertible, $\OO(\Gamma_F) \otimes_{\tau, \chi} \OO(\alg{\tilde G}\LVee)$ is a Hopf algebra.
\qed

We have discussed constructions of a (1,2)-cocycle $\tau$ and (2,1)-cocycle $\chi$.  Now we demonstrate compatibility.
\begin{construction}
Recall that $Q\colon Y \rightarrow \ZZ$ is a Weyl-invariant integer-valued quadratic form, and $C$ is a fair bisector of $Q$.  Define toral $\ZZ$-linear maps $\tau$ and $\chi$ as before by: 
$$\tau(1_y^\dual)(\gamma_1, \gamma_2) = h(\gamma_1, \gamma_2)^{Q(y)}, \quad \chi(1_{y_1}^\dual, 1_{y_2}^\dual) = h(\gamma, \gamma)^{C(y_1, y_2)}.$$
We have seen that $\tau$ is a (1,2)-cocycle, and $\chi$ is a symmetric (2,1)-cocycle.  Moreover, the pair $(\tau, \chi)$ is compatible and $\chi$ is involutive.
\end{construction}
\proof
We begin by proving the first compatibility condition (Com1).  The left side of this condition is
$$L(a,b) = \sum m_{a_1, a_2}^a m_{b_1, b_2}^b \chi(a_1^\dual, b_1^\dual) \otimes a_2^\dual b_2^\dual.$$
The right side of (Com2) is
$$R(a,b) = \sum m_{a_1, a_2}^a m_{b_1, b_2}^b \chi(a_2^\dual, b_2^\dual) \otimes a_1^\dual b_1^\dual.$$
The terms in $L(a,b)$ vanish if $a_1, b_1 \not \in \BB_{\alg{\tilde T}}$; nonvanishing implies that $a = a_2$ and $a = 1_{u_1}$, and $b = b_2$ and $b_1 = 1_{v_1}$, where $u_1 = \lt{a}$, $v_1 = \lt{b}$.

In this way, given $a,b \in \BB$, with $u_1 = \lt{a}$, $u_2 = \rt{a}$, $v_1 = \lt{b}$, $v_2 = \rt{b}$, we find that
$$L(a,b)(\gamma) = h(\gamma, \gamma)^{C(u_1, v_1)} \otimes a^\dual b^\dual, \quad R(a,b)(\gamma) = h(\gamma, \gamma)^{C(u_2, v_2)} \otimes a^\dual b^\dual.$$
Now, recall $C$ is bilinear and {\em fair}.  If $n$ is odd, then $h(\gamma, \gamma) = 1$ and there is nothing left to show.  If $n$ is even and $Q(\alpha_i^\vee)$ is odd, then $n_i$ is even and so $C(\tilde \alpha_i^\vee, y) = n_i C(\alpha_i^\vee, y) = 0$ for all $y \in \tilde Y$.  If $n$ is even and $Q(\alpha_i^\vee)$ is even, then fairness of $C$ implies that $C(\tilde \alpha_i^\vee, y) = C(\alpha_i^\vee, y) = 0$ for all $y \in \tilde Y$.  

Thus, as a consequence of fairness of $C$, and the facts that $u_2 - u_1 = \nt{a} \in \ZZ[\tilde \Delta^\vee]$, $v_2 - v_1 = \nt{b} \in \ZZ[\tilde \Delta^\vee]$, we find that $C(u_1, v_1) = C(u_2, v_2)$.  Hence $L(a,b) = R(a,b)$.

Next we prove the compatibility condition (Com2).  The left side of this condition is
$$L(a,b) = \sum_{a_1, a_2, a_3} \sum_{b_1, b_2, b_3} m_{a_1, a_2, a_3}^a m_{b_1, b_2, b_3}^b \tau(a_3^\dual) \cdot \tau(b_3^\dual) \cdot \left( \chi(a_1^\dual, b_1^\dual) \otimes \chi(a_2^\dual, b_2^\dual) \right).$$
The terms of $L(a,b)$ vanish except when $a_1, a_2, a_3$ and $b_1, b_2, b_3$ are in $\BB_{\alg{\tilde T}}$; this nonvanishing implies that $a,b \in \BB_{\alg{\tilde T}}$ as well.  Thus the nonvanishing terms occur when $a = a_1 = a_2 = a_3 = 1_u$ and $b = b_1 = b_2 = b_3 = 1_v$ for some $u,v \in \tilde Y$.  The corresponding term $L(a,b)$, evaluated at $\gamma_1, \gamma_2$ is
\begin{eqnarray*}
L(1_u, 1_v)(\gamma_1, \gamma_2) & = & h(\gamma_1, \gamma_2)^{Q(u)} h(\gamma_1, \gamma_2)^{Q(v)} \cdot h(\gamma_1, \gamma_1)^{C(u,v)} \cdot h(\gamma_2, \gamma_2)^{C(u,v)}, \\
& = & h(\gamma_1, \gamma_2)^{Q(u) + Q(v)} \cdot h(\gamma_1, \gamma_1)^{C(u,v)} \cdot h(\gamma_2, \gamma_2)^{C(u,v)}.
\end{eqnarray*}

The right side of (Com2) is
$$R(a,b) = \sum_{a_1, a_2} \sum_{b_1, b_2} m_{a_1, a_2}^a m_{b_1, b_2}^b  \tau(a_1^\dual b_1^\dual) \cdot \Delta \left( \chi(a_2^\dual, b_2^\dual) \right).$$
The terms vanish, except when $a = a_1 = a_2 = 1_u$ and $b = b_1 = b_2 = 1_v$, for some $u,v \in \tilde Y$.  In this case, $a_1^\dual b_1^\dual = 1_{u + v}^\dual$.  We find that
$$R(1_u, 1_v)(\gamma_1, \gamma_2) = h(\gamma_1, \gamma_2)^{Q(u + v)} \cdot h(\gamma_1 \gamma_2, \gamma_1 \gamma_2)^{C(u,v)}.$$

Observe that $h(\gamma_1 \gamma_2, \gamma_1 \gamma_2) = h(\gamma_1, \gamma_1) h(\gamma_2, \gamma_2)$, by skew-symmetry of the Hilbert symbol.  Also, note that $Q(u+v) \ident Q(u) + Q(v)$, mod $n \ZZ$, when $u,v \in \tilde Y$.  Hence $L(a,b) = R(a,b)$ and (Com2) is satisfied.

Finally, we observe that $\chi$ is involutive; we have
$$\sum m_{a_1, a_2}^a m_{b_1, b_2} \chi(a_1^\dual, b_1^\dual) \cdot \chi(a_2^\dual, b_2^\dual) = \epsilon(a^\dual) \epsilon(b^\dual).$$
Indeed, the terms on the left vanish unless $a = a_1 = a_2 = 1_u$ and $b = b_1 = b_2 = 1_v$ for some $u,v \in \tilde Y$.  In this case, the left side is $\chi(1_u^\dual, 1_v^\dual) \chi(1_u^\dual, 1_v^\dual) = 1$.
\qed

The above construction finishes the proof of the following theorem.
\begin{thm}
The double-twist $\OO(\Gamma_F) \otimes_{\tau, \chi} \OO(\alg{\tilde G}\LVee)$ is a commutative Hopf algebra over $\ZZ$.
\end{thm}

Define ${}^L \alg{\tilde G}\LVee = \Spec \left( \OO(\Gamma_F) \otimes_{\tau, \chi} \OO(\alg{\tilde G}\LVee) \right)$ to be the resulting group scheme over $\ZZ$.  There is a sequence of Hopf algebra homomorphisms over $\ZZ$:
$$\OO(\Gamma_F) \rightarrow \OO(\Gamma_F) \otimes_{\tau, \chi} \OO(\alg{\tilde G}\LVee) \rightarrow \OO(\alg{\tilde G}\LVee)$$
given by the $\ZZ$-linear maps satisfying $f \mapsto (f \otimes 1)$ and $(f \otimes a^\dual) \mapsto \epsilon(f) \cdot a^\dual$ for all $f \in \OO(\Gamma_F)$ and all $a \in \BB$.  

Thus there are corresponding homomorphisms of group schemes over $\ZZ$:
$$\alg{\tilde G}\LVee \rightarrow {}^L \alg{\tilde G}\LVee \rightarrow \alg{\Gamma}_F.$$
\begin{proposition}
The map $\alg{\tilde G} \LVee \rightarrow {}^L \alg{\tilde G}\LVee$ is a closed embedding of group schemes over $\ZZ$, and identifies $\alg{\tilde G}\LVee$ with $\Ker \left( {}^L \alg{\tilde G}\LVee \rightarrow \alg{\Gamma}_F \right)$.  The map $ {}^L \alg{\tilde G}\LVee \rightarrow \alg{\Gamma}_F$ is faithfully flat.
\end{proposition}
\proof
The ring homomorphism $\OO(\Gamma_F) \otimes_{\tau, \chi} \OO(\alg{\tilde G}\LVee) \rightarrow \OO(\alg{\tilde G}\LVee)$ satisfying $(f \otimes a^\dual) \mapsto \epsilon(f) \cdot a^\dual$ is surjective; its kernel is $\Ker(\epsilon) \otimes \OO(\alg{\tilde G}\LVee)$.  Hence we find that the map of schemes $\alg{\tilde G}\LVee \rightarrow {}^L \alg{\tilde G}\LVee$ is a closed embedding, and identifies
$$\alg{\tilde G}\LVee \ident \Ker \left( {}^L \alg{\tilde G}\LVee \rightarrow \alg{\Gamma}_F \right).$$  

As a $\OO(\Gamma_F)$-module, via $f \mapsto (f \otimes 1)$, we find that $\OO(\Gamma_F) \otimes_{\tau, \chi} \OO(\alg{\tilde G}\LVee)$ is free with basis $\{ 1 \otimes a^\dual : a \in \BB \}$.  Hence the map ${}^L \alg{\tilde G}\LVee \rightarrow \alg{\Gamma}_F$ is faithfully flat.
\qed

We see in the next section that this sequence splits, at least when $n$ is odd or $n = 2$, with the middle term isomorphic to the direct product of the outer terms, but not without some other choices and a possible extension of scalars from $\ZZ$ to $\ZZ[i]$.

For what follows, and for general considerations of functoriality, it is useful to define a category of which ${}^L \alg{\tilde G}\LVee$ is a prototypical object.
\begin{definition} 
Let $\Cat{EGp} = \Cat{EGp}_{\tilde \Psi, F}$ (it depends on the root datum $\tilde \Psi$ and the field $F$) denote the category whose objects are sequences of group schemes over $\ZZ$ and group scheme morphisms over $\ZZ$:  
$$\alg{\tilde G}\LVee \rightarrow \alg{E} \rightarrow \alg{\Gamma}_F.$$
A morphism in $\Cat{EGp}$ from $\alg{\tilde G}\LVee \rightarrow \alg{E}_1 \rightarrow \alg{\Gamma}_F$ to $\alg{\tilde G}\LVee \rightarrow \alg{E}_2 \rightarrow \alg{\Gamma}_F$ will mean a group scheme homomorphism $e\colon \alg{E}_1 \rightarrow \alg{E}_2$ making the following diagram commute:
$$\xymatrix{\alg{\tilde G}\LVee \ar[r] \ar[d]^{=} &  \alg{E}_1 \ar[r] \ar[d]^{e} & \alg{\Gamma}_F \ar[d]^{=}  \\
\alg{\tilde G}\LVee \ar[r] &  \alg{E}_2 \ar[r] & \alg{\Gamma}_F.}$$
\end{definition}

\subsection{Change of bisector}

We have twisted the Hopf algebra $\OO(\Gamma_F) \otimes  \OO(\alg{\tilde G}\LVee)$ by a compatible pair $(\tau, \chi)$ (with $\chi$ involutive) in order to obtain a new commutative Hopf algebra $\OO(\Gamma_F) \otimes_{\tau, \chi} \OO(\alg{\tilde G}\LVee)$ over $\ZZ$.  The data required for this construction included a quadratic form $Q$, and a fair bisector $C$ of $Q$.  

Recalling that the metaplectic group $\tilde G$ depends functorially on a pair $(C, \eta)$ (via the functor $\Cat{Inc}_\alg{G}$), it is important to see that the Hopf algebra construction also depends functorially on $(C, \eta)$.  We demonstrate that a morphism from $(C_1, \eta)$ to $(C_2, \eta)$ (in the category $\Cat{Bis}_Q^{f, sc}$) yields a morphism of Hopf algebras. 

\begin{remark}
We only consider morphisms $(C_1, \eta) \rightarrow (C_2, \eta)$, with the same component $\eta\colon Y_{sc} \rightarrow F^\times$ in source and target.  Without this assumption, it seems that another construction is required.  This indicates that our construction of an L-group may be insufficient, especially if one wants to work with nonsplit groups via descent.
\end{remark}

Consider a morphism $H$ in the category $\Cat{Bis}_Q^{f,sc}$ from $(C_1, \eta)$ to $(C_2, \eta)$.  This means that $H: Y / 2 Y \rightarrow \ZZ / 2 \ZZ$ is a function satisfying
\begin{enumerate}
\item
$C_2(y_1, y_2) - C_1(y_1, y_2) = H(y_1 + y_2) - H(y_1) - H(y_2)$;
\item
$0 = C_2(\tilde \alpha_i^\vee, y) - C_1(\tilde \alpha_i^\vee, y) = H(y + \tilde \alpha_i^\vee) - H(y) - H(\tilde \alpha_i^\vee)$;
\item
$1 = \eta(\alpha_i^\vee) / \eta(\alpha_i^\vee) = (-1)^{H(\alpha_i^\vee)}$.
\end{enumerate}

As a consequence of these identities, we find that
\begin{equation}
\label{HC}
H(y + \tilde \alpha_i^\vee) = H(y) \mbox{ for all } i \in I, y \in Y.
\end{equation}

Hence $H$ gives a well defined function from $Y / (2Y + \ZZ[\tilde \Delta^\vee])$ to $\ZZ / 2 \ZZ$.

Define a toral $\ZZ$-linear map from $\OO(\alg{\tilde G}\LVee)$ to $\ZZ$ by $\omega_H(1_y^\dual) = (-1)^{H(y)}$.  If $\chi_1, \chi_2$ are the (2,1)-cocycles constructed from $C_1, C_2$, respectively, then
$$\chi_1(1_u^\dual, 1_v^\dual) \omega_H(1_{u+v}^\dual) = \chi_2(1_u^\dual, 1_v^\dual) \omega_H(1_u^\dual) \omega_H(1_v^\dual),$$
for all $u,v \in \tilde Y$.  

Define a $\ZZ$-linear endomorphism: $\alpha_H\colon \OO(\Gamma_F) \otimes \OO(\alg{\tilde G}\LVee) \rightarrow \OO(\Gamma_F) \otimes \OO(\alg{\tilde G}\LVee)$, by
$$\alpha_H(f \otimes a^\dual) = \sum_{a_1, a_2} m_{a_1, a_2}^a f \cdot \omega_H(a_1^\dual) \otimes a_2^\dual.$$
\begin{proposition}
Let $\tau$ be the (1,2)-cocycle constructed before.  Let $\chi_1$ and $\chi_2$ be (2,1)-cocycles constructed via fair bisectors $C_1$ and $C_2$ respectively.  Then for any morphism $H\colon (C_1, \eta) \rightarrow (C_2, \eta)$ of bisectors, the $\ZZ$-linear endomorphism $\alpha_H$ is a Hopf algebra isomorphism from $\OO(\Gamma_F) \otimes_{\tau, \chi_1} \OO(\alg{\tilde G}\LVee)$ to $ \OO(\Gamma_F) \otimes_{\tau, \chi_2} \OO(\alg{\tilde G}\LVee)$.
\end{proposition}
\proof
It is straightforward to check that the endomorphism $\alpha_H$ preserves unit and counit.  To check that $\alpha_H$ is an algebra homomorphism, we compute for any $a,b \in \BB$:
\begin{eqnarray*}
\alpha_H( (f \otimes a^\dual) \cdot_{\chi_1} (g \otimes b^\dual)) & = & \sum m_{a_1, a_2}^a m_{b_1, b_2}^b \alpha_H(fg \chi_1(a_1^\dual, b_1^\dual) \otimes a_2^\dual b_2^\dual) \\
& = & \sum m_{a_1, a_2, a_3}^a m_{b_1, b_2, b_3}^b fg \chi_1(a_1^\dual, b_1^\dual) \omega_H(a_2^\dual b_2^\dual) \otimes a_3^\dual b_3^\dual \\
& = & fg \chi_1(1_u^\dual, 1_v^\dual) \omega_H(1_{u+v}^\dual) \otimes a^\dual b^\dual \\
& = & fg \chi_2(1_u^\dual, 1_v^\dual) \omega_H(1_u^\dual) \omega_H(1_v^\dual) \otimes a^\dual b^\dual \\
& = & \alpha_H(f \otimes a^\dual) \cdot_{\chi_2} \alpha_H(g \otimes b^\dual).
\end{eqnarray*}
In the middle, we use the fact that all terms vanish except those where $a_1 = a_2 = 1_u$ and $b_1 = b_2 = 1_v$, where $u = \lt{a}$, $v = \lt{b}$.

To check that $\alpha_H$ is a coalgebra homomorphism, we compute
\begin{eqnarray*}
[\alpha_H \otimes \alpha_H](\Delta_\tau(f \otimes a^\dual)) & = & [\alpha_H \otimes \alpha_H] \left( \sum m_{a_1, a_2, a_3}^a \Delta(f) \tau(a_3^\dual) \otimes (a_1^\dual \otimes a_2^\dual) \right) \\
& = & \sum m_{a_1, a_2, a_3, a_4, a_5}^a \Delta(f) \tau(a_5^\dual) \\ 
& & \left( \omega_H(a_{1}^\dual) \otimes \omega_H(a_{3}^\dual) \right) \otimes (a_{2}^\dual \otimes a_{4}^\dual) \\
& = & \sum m_{a_1, a_2}^a \Delta(f) \tau(1_{v_2}^\dual) \\
& &  (\omega_H(1_{u_1}^\dual) \otimes \omega_H(1_{v_1}^\dual)) \otimes (a_1^\dual \otimes a_2^\dual) \\
& = & \sum m_{a_1, a_2}^a \Delta(f) \Delta(\omega(1_{u_1}^\dual)) \tau(1_{v_2}^\dual) \otimes (a_1^\dual \otimes a_2^\dual) \\
& = & \sum m_{a_1, a_2}^a \Delta(f) \Delta(\omega(a_1^\dual)) \tau(a_4^\dual) \otimes (a_2^\dual \otimes a_3^\dual) \\
& = & \sum m_{a_1, a_2}^a \Delta_\tau(f \omega(a_1^\dual) \otimes a_2^\dual) \\
& = & \Delta_\tau( \alpha_H(f \otimes a^\dual))
\end{eqnarray*}
Here, we write $u_1 = \lt{a}_1$, $u_2 = \rt{a}_1$, $v_1 = \lt{a}_2$, $v_2 = \rt{a}_2$.  All steps use techniques from before, except perhaps in the middle, where we use the equality
$$\omega_H(1_{u_1}^\dual) \otimes \omega_H(1_{v_1}^\dual) = \Delta(\omega_H(1_{u_1}^\dual)).$$
This follows from Equation \ref{HC}, using the fact that $u_1 - v_1 = \nt{a} \in \ZZ[\tilde \Delta^\vee]$ and so $H(u_1) = H(v_1)$.

\begin{corollary}
Let $\eta\colon Y_{sc} \rightarrow F^\times$ be a homomorphism, and let $Q\colon Y \rightarrow \ZZ$ be a quadratic form.  Let $\Cat{Bis}_{Q, \eta}^{f,sc}$ be the full subcategory of $\Cat{Bis}_{Q}^{f,sc}$ consisting of objects of the form $(C,\eta)$.  

Then our construction of the group scheme ${}^L \alg{\tilde G}\LVee$ over $\ZZ$ defines a functor from $\Cat{Bis}_{Q,\eta}^{f,sc}$ to the category $\Cat{EGp}$.  In other words, the group scheme ${}^L \alg{\tilde G}\LVee$ depends functorially on the choice of bisector $C$ (with $Q, \eta$ fixed). 
\end{corollary}

\begin{remark}
Our construction is not entirely satisfactory, since it does not take $\eta$ into account, nor does it allow for morphisms in $\Cat{CExt}(Y, F^\times)$ which intertwine different maps $\eta_1, \eta_2: Y_{sc} \rightarrow F^\times$.  But our construction is satisfactory when $\eta$ is trivial, and takes into account automorphisms of bisectors.
\end{remark}

\subsection{Trivialization}

We have discussed a method to twist the Hopf algebra $\OO(\Gamma_F) \otimes \OO(\alg{\tilde G}\LVee)$ by a compatible pair $(\tau, \chi)$, with $\chi$ convolution-invertible, to obtain a new Hopf algebra $\OO(\Gamma_F) \otimes_{\tau, \chi} \OO(\alg{\tilde G}\LVee)$.  Here we discuss the ``trivialization'' of this double-twist after extension of scalars. 

Let $A$ be a flat commutative $\ZZ$-algebra.  Let $\omega$ be an $A$-linear map from $\OO_A(\alg{\tilde G}\LVee)$ to $\OO_A(\Gamma_F)$, where subscripts denote extension of scalars ($\bullet_A = \bullet \otimes_\ZZ A$).  Let $\alpha_\omega$ be the $A$-linear endomophism of $\OO_A(\Gamma_F) \otimes_A \OO_A(\alg{\tilde G}\LVee)$ given by
$$\alpha_\omega(f \otimes a^\dual) = \sum_{a_1, a_2} m_{a_1, a_2}^a f \cdot \omega(a_1^\dual) \otimes a_2^\dual.$$

\begin{definition}
We say that $\omega$ is convolution-invertible if there exists an $A$-linear map $\rho\colon \OO(\alg{\tilde G}\LVee) \rightarrow \OO(\Gamma_F)$ such that $\sum_{a_1, a_2} m_{a_1, a_2}^a \omega(a_1^\dual) \rho(a_2^\dual) = \epsilon(a^\dual)$ for all $a \in \BB$.  In this case, we say that $\rho$ is a convolution-inverse of $\omega$.
\end{definition}

\begin{proposition}
When $\rho$ is a convolution-inverse of $\omega$, $\alpha_\rho \circ \alpha_\omega = \Id$.
\end{proposition}
\proof
We check directly:
\begin{eqnarray*}
[\alpha_\rho \circ \alpha_\omega](f \otimes a^\dual) & = & \alpha_\rho \left( \sum_{a_1, a_2} m_{a_1, a_2}^a f \cdot \omega(a_1^\dual) \otimes a_2^\dual \right) \\
& = & \sum_{a_1, a_2, a_3} m_{a_1, a_2, a_3}^a f \cdot \omega(a_1^\dual) \rho(a_2^\dual) \otimes a_3^\dual \\
& = & \sum_{a_1, a_2} m_{a_1, a_2}^a f \epsilon(a_1^\dual)  \otimes a_2^\dual =  f \otimes a^\dual.
\end{eqnarray*}
\qed

\begin{proposition}
Let $(\tau, \chi)$ be a compatible pair as defined in the previous sections.  Suppose that $\omega\colon \OO_A(\alg{\tilde G}\LVee) \rightarrow \OO_A(\Gamma_F)$ satisfies the following axioms.
\begin{enumerate}
\item[(Triv1)] $\omega(1) = 1$, and $\epsilon \circ \omega = \epsilon$;
\item[(Triv2)] For all $a,b \in \BB$,
$$\omega(a^\dual b^\dual) = \sum m_{a_1, a_2}^a m_{b_1, b_2}^b \omega(a_1^\dual) \omega(b_1^\dual) \chi(a_2^\dual, b_2^\dual).$$
\item[(Triv3)]  For all $a \in \BB$,
$$\sum m_{a_1, a_2, a_3, a_4}^a \omega(a_1^\dual) \otimes \omega(a_3^\dual) \otimes a_2^\dual \otimes a_4^\dual = \sum  m_{a_1, a_2, a_3, a_4}^a \Delta(\omega(a_1^\dual)) \tau(a_4^\dual) \otimes a_2^\dual \otimes a_3^\dual.$$
(Both sides are elements of $\OO_A(\Gamma_F \times \Gamma_F) \otimes  \OO_A(\alg{\tilde G}\LVee) \otimes  \OO_A(\alg{\tilde G}\LVee)$.)
\end{enumerate}
Then $\alpha_\omega$ is a homomorphism of Hopf algebras over $A$, from the untwisted Hopf algebra $\OO_A(\Gamma_F) \otimes \OO_A(\alg{\tilde G}\LVee)$ to the double-twisted Hopf algebra $\OO_A(\Gamma_F) \otimes_{\tau,\chi} \OO_A(\alg{\tilde G}\LVee)$ (with tensor products understood over $A$).  If, moreover, $\omega$ is convolution-invertible, then $\alpha_\omega$ is an isomorphism of Hopf algebras over $A$.
\end{proposition}
\proof
First we use (Triv1) to demonstrate that $\alpha_\omega$ preserves unit and counit.  Observe that $\alpha_\omega(1 \otimes 1) = 1 \otimes 1$, since $\Delta(1) = 1 \otimes 1$ and $\omega(1) = 1$.  Next, observe that 
$$\epsilon(\alpha_\omega(f \otimes a^\dual)) = \sum m_{a_1, a_2}^a \epsilon(f) \epsilon(\omega(a_1^\dual)) \cdot \epsilon(a_2^\dual) = \sum m_{a_1, a_2}^a \epsilon(f) \epsilon(a_1^\dual) \epsilon(a_2^\dual) = \epsilon(f) \cdot \epsilon(a^\dual).$$
Thus $\alpha_\omega$ preserves unit and counit.

Now we check that $\alpha_\omega$ is an algebra homomorphism, using (Triv2):
\begin{eqnarray*}
\alpha_\omega( (f \otimes a^\dual) \cdot (g \otimes b^\dual) ) & = & \alpha_\omega(fg \otimes a^\dual b^\dual) \\
& = & \sum_{a_1, a_2} \sum_{b_1, b_2} m_{a_1, a_2}^a m_{b_1, b_2}^b f g \omega(a_1^\dual b_1^\dual) \otimes a_2^\dual b_2^\dual \\
& = &  \sum_{{a_1, a_2, a_3} \atop {b_1, b_2, b_3} } m_{a_1, a_2, a_3}^a m_{b_1, b_2, b_3}^b f g \omega(a_1^\dual) \omega(b_1^\dual) \chi(a_2^\dual, b_2^\dual) \otimes a_3^\dual b_3^\dual \\
& = &  \sum_{a_1, a_2} \sum_{b_1, b_2} m_{a_1, a_2}^a m_{b_1, b_2}^b \left( f \omega(a_1^\dual) \otimes a_2^\dual \right) \cdot_\chi \left( g \omega(b_1^\dual) \otimes b_2^\dual \right) \\
& = & \alpha_\omega(f \otimes a^\dual) \cdot_\chi \alpha_\omega(g \otimes b^\dual).
\end{eqnarray*}

Finally, we check that $\alpha_\omega$ is a coalgebra homomorphism, using (Triv3):
\begin{eqnarray*}
[\alpha_\omega \otimes \alpha_\omega](\Delta(f \otimes a^\dual)) & = & [\alpha_\omega \otimes \alpha_\omega] \left( \sum_{a_1, a_2} m_{a_1, a_2}^a \Delta(f)  \otimes a_1^\dual \otimes a_2^\dual \right) \\
& = & \sum_{a_1, a_2, a_3, a_4} m_{a_1, a_2, a_3, a_4}^a \Delta(f) (\omega(a_1^\dual) \otimes \omega(a_3^\dual) ) \otimes (a_2^\dual \otimes a_4^\dual) \\
& = & \sum_{a_1, a_2, a_3, a_4} m_{a_1, a_2, a_3, a_4}^a\Delta(f) \Delta(\omega(a_1^\dual)) \tau(a_4^\dual) \otimes (a_2^\dual \otimes a_3^\dual) \\
& = & \sum_{a_1, a_2} m_{a_1, a_2}^a \Delta_\tau(f \omega(a_1^\dual) \otimes a_2^\dual) \\
& = & \Delta_\tau(\alpha_\omega(f \otimes a^\dual)).
\end{eqnarray*}
\qed

\subsection{The L-group of a double-cover}
Recall that $F$ is a local field and $F \not \isom \CC$.  Using the Weyl-invariant quadratic form $Q\colon Y \rightarrow \ZZ$, the positive integer $n$, and the fair bisector $C$ of $Q$, we have constructed a compatible pair $(\tau, \chi)$, from which we have a Hopf algebra over $\ZZ$:
$$\OO({}^L \alg{\tilde G}\LVee) = \OO(\Gamma_F) \otimes_{\tau, \chi} \OO(\alg{\tilde G}\LVee).$$

If $n$ is odd, then both $\chi$ and $\tau$ are trivial, i.e.\ $\chi(a^\dual, b^\dual) = \epsilon(a^\dual) \cdot \epsilon(b^\dual)$ and $\tau(a^\dual) = \epsilon(a^\dual) \otimes \epsilon(a^\dual)$.  In this case,
$$\OO({}^L \alg{\tilde G}\LVee) = \OO(\Gamma_F) \otimes_{\tau, \chi} \OO(\alg{\tilde G}\LVee) = \OO(\Gamma_F) \otimes \OO(\alg{\tilde G}\LVee)$$
as commutative Hopf algebras over $\ZZ$.  In other words, the L-group ${}^L \alg{\tilde G}\LVee$ is equal to the na\"ive L-group $\alg{\Gamma}_F \times \alg{\tilde G}\LVee$, as group schemes over $\ZZ$.

Assume for the remainder of this section that $n = 2$ and $\Char(F) \neq 2$.  Interesting arithmetic problems arise when $n$ is even, and we have not yet tackled the case when $n$ is a multiple of $4$.

\subsubsection{The Hilbert symbol as a coboundary}
Recall that $\Hilb\colon F^\times \times F^\times \rightarrow \mu(F)$ denotes the Hilbert symbol.  Let $\Hilb_2\colon F^\times \times F^\times \rightarrow \mu_2$ be the quadratic Hilbert symbol, so $\Hilb_2 = \Hilb^{\# \mu(F) / 2}$.  For any continuous nontrivial additive character $\psi\colon F \rightarrow \CC^\times$, and any nonzero element $x \in F^\times$, let $\weil(x, \psi)$ denote the {\em Weil index} \cite[\S A.3]{RR} (based on the original treatment of \cite{Weil}).  When $\psi$ is such a character, and $c \in F^\times$, we write ${}^c \psi$ for the character given by ${}^c \psi(x) = \psi(cx)$ for all $x \in F$.  This gives a transitive action of $F^\times$ on the set of continuous nontrivial additive characters of $F$.  

The Weil index $\weil(x, \psi)$ is a fourth root of unity in $\CC$, which satisfies
\begin{enumerate}
\item[(W1)]
$\weil(x,\psi)^2 = \Hilb_2(x,x)$ for all $x \in F^\times$;
\item[(W2)]
$\Hilb_2(x,y) = \weil(xy, \psi) \weil(x,\psi)^{-1} \weil(y,\psi)^{-1}$ for all $x,y \in F^\times$;
\item[(W3)]
$\weil(x c^2, \psi) = \weil(x, \psi)$ for all $c,x \in F^\times$;
\item[(W4)]
For all $x,c \in F^\times$, $\weil(x, {}^c \psi) = \Hilb_2(x,c) \weil(x,\psi)$.
\end{enumerate}

\begin{remark}
Using (W2), with $x = y = 1$, we find that $\weil(1, \psi) =1$ for all $\psi$.
\end{remark}

Let $\ZZ[i]$ be the smallest subring of $\CC$ containing $\ZZ$ and $i$.  Thus we identify $\weil(x,\psi)$ as an element of $\mu_4 = \ZZ[i]^\times$.  The above properties of the Weil index yield the following.
\begin{proposition}
The Weil index gives a bijection between the following two sets:
\begin{enumerate}
\item
The set of orbits $\Omega$ for the action of $F^{\times 2}$ on the set of nontrivial continuous additive characters $F \rightarrow \CC^\times$.
\item
The set $W$ of functions $w\colon F^\times / F^{\times 2} \rightarrow \mu_4$ satisfying 
$$\Hilb_2(x,y) = w(xy) w(x)^{-1} w(y)^{-1}.$$
\end{enumerate}
\end{proposition}
\proof
Each nontrivial continuous additive character $\psi$ yields a function $w = \weil(\bullet, \psi)$, which depends only on the $F^{\times 2}$-orbit of $\psi$ by (W4).  The map from $\Omega$ to $W$ being well-defined, it remains to check that the map is a bijection.  Property (W4) and the nondegeneracy of the quadratic Hilbert symbol implies that the map $\Omega \rightarrow W$ is injective -- the function $\weil(\bullet, \psi)$ determines the $F^{\times 2}$-orbit of the additive character $\psi$.

The set $W$ is a torsor for $\Hom(F^\times / F^{\times 2}, \mu_4) = \Hom(F^\times / F^{\times 2}, \mu_2)$, which has cardinality $\# (F^\times / F^{\times 2})$.  Similarly, $F^\times$ acts transitively on the nontrivial additive characters of $F$, so $\Omega$ has cardinality $\# (F^\times / F^{\times 2})$.  Hence $\# \Omega = \# W$, so the injective map from $\Omega$ to $W$ is bijective.
\qed

\subsubsection{The bisector as a coboundary}
We consider the following elementary abelian 2-group:
$$\bar Y = \frac{\tilde Y}{2Y + \ZZ[\tilde \Delta^\vee]}.$$
Note that $2Y \subset \tilde Y$, since $n = 2$.  For any $y \in \tilde Y$, write $\bar y$ for its image in $\bar Y$.  Recall that $C\colon Y \times Y \rightarrow \ZZ / 2 \ZZ$ is a fair bisector of $Q$, so $C$ restricts to $\tilde Y \times \tilde Y$ and from there it factors through a symmetric bilinear form
$$\bar C\colon \bar Y \times \bar Y \rightarrow \ZZ / 2 \ZZ.$$

\begin{definition}
Let $\delta\colon \ZZ / 2 \ZZ \rightarrow \ZZ / 4 \ZZ$ be the unique injective group homomorphism:   $\delta(0) = 0$, $\delta(1) = 2$.
For any symmetric bilinear form $\bar C\colon \bar Y \times \bar Y \rightarrow \ZZ / 2 \ZZ$ as above, a {\em tetractor} of $\bar C$ is a function $\kappa\colon \bar Y \rightarrow \ZZ / 4 \ZZ$ satisfying
$$\kappa(\bar y_1 + \bar y_2) - (\kappa(\bar y_1) + \kappa(\bar y_2)) = [\delta \circ \bar C](\bar y_1, \bar y_2).$$
\end{definition}
Note that if $\kappa$ is a tetractor of $\bar C$, then $-\kappa$ is also a tetractor of $\bar C$.  Moreover, if $\bar C$ is nonzero, then $\kappa \neq - \kappa$.  

\begin{proposition}
\label{Kappa}
Any symmetric bilinear form $\bar C\colon \bar Y \times \bar Y \rightarrow \ZZ / 2 \ZZ$ has at least one tetractor.  Every tetractor $\kappa$ satisfies $\kappa(0) = 0$ and for all $y \in \tilde Y$, $\kappa(\bar y) \ident Q(y)$ mod $2$.\end{proposition}
\proof
Consider the bilinear map $\delta \circ \bar C\colon \bar Y \times \bar Y \rightarrow \ZZ / 4 \ZZ$.  Note that $\delta$ induces the zero map in cohomology, when applied to the coefficients of cohomology of $\ZZ / 2 \ZZ$:
$$H^2(\delta) = 0\colon H^2( \ZZ / 2 \ZZ, \ZZ / 2 \ZZ) \rightarrow H^2( \ZZ / 2 \ZZ, \ZZ / 4 \ZZ).$$
Since $\bar Y$ is an elementary abelian 2-group, we find that $\delta \circ \bar C$ is a coboundary when viewed as an element of $Z^2(\bar Y, \ZZ / 4 \ZZ)$.  In other words, there exists a function $\kappa\colon \bar Y \rightarrow \ZZ / 4 \ZZ$ such that
\begin{equation}
\label{KapId}
[\delta \circ \bar C](\bar y_1,\bar y_2) = \kappa(\bar y_1 + \bar y_2) - \kappa(\bar y_1) - \kappa(\bar y_2).
\end{equation}
Hence $\bar C$ has a tetractor.

Now let $\kappa$ be a tetractor.  Observe that if $\bar C(\bar y, \bar y) = 0$ then $0 = [\delta \circ \bar C](\bar y, \bar y) = \kappa(\bar y + \bar y) - \kappa(\bar y) - \kappa(\bar y) = 2 \kappa(\bar y)$.  Hence if $\bar C(\bar y, \bar y) = 0$, then $\kappa(\bar y) \in \{ 0, 2 \}$.  On the other hand, if $\bar C(\bar y, \bar y) = 1$, then $2 = [\delta \circ \bar C](\bar y, \bar y) = 2 \kappa(\bar y)$, so we find that $\kappa(\bar y) \in \{1, 3 \}$.  In any case, we find that 
$$\bar C(\bar y, \bar y) \ident \kappa(\bar y), \mbox{ mod } 2.$$
Since $Q(y) \ident C(y,y)$, mod $2$, the identity is proven.

Finally, note that
$$\kappa(0) = -(\kappa(0 + 0) - \kappa(0) - \kappa(0)) = - [d \circ \bar c](0,0) = 0,$$
so $\kappa(0) = 0$.  
\qed

\subsubsection{The L-group and its trivialization}

Recall that the compatible pair $(\tau, \chi)$ is defined by
$$\tau(1_y^\dual)(\gamma_1, \gamma_2) = h(\gamma_1, \gamma_2)^{Q(y)}, \quad \chi(1_{y_1}^\dual, 1_{y_2}^\dual) = h(\gamma, \gamma)^{C(y_1, y_2)}.$$
These are used to define the Hopf algebra over $\ZZ$,
$$\OO({}^L \alg{\tilde G}\LVee) = \OO(\Gamma_F) \otimes_{\tau, \chi} \OO(\alg{\tilde G}\LVee),$$
which fits into a natural sequence of Hopf algebras over $\ZZ$:
$$\OO(\Gamma_F) \rightarrow \OO({}^L \alg{\tilde G}\LVee) \rightarrow \OO(\alg{\tilde G}\LVee).$$
This corresponds to a sequence of group schemes over $\ZZ$:
$$\alg{\tilde G}\LVee \rightarrow {}^L \alg{\tilde G}\LVee \rightarrow \alg{\Gamma}_F.$$

Since we assume $n = 2$ here, we have
$$h(\gamma_1, \gamma_2) = \Hilb_2(\rec(\gamma_1), \rec(\gamma_2)).$$

\begin{construction}
\label{TrivConstruct}
Let $\psi\colon F \rightarrow \CC^\times$ be a nontrivial continuous additive character.  Let $\kappa\colon \bar Y \rightarrow \ZZ / 4 \ZZ$ be a tetractor of $\bar C$ (as exists by Proposition \ref{Kappa}): 
$$[\delta \circ \bar C](\bar y_1,\bar y_2) = \kappa(\bar y_1 + \bar y_2) - \kappa(\bar y_1) - \kappa(\bar y_2).$$
Define a toral $\ZZ[i]$-linear map $\omega = \omega_{\psi, \kappa}\from \OO_{\ZZ[i]}(\alg{\tilde G}\LVee) \rightarrow \OO_{\ZZ[i]}(\Gamma_F)$ by
$$\omega(1_y^\dual)(\gamma) = \weil(\rec(\gamma), \psi)^{\kappa(\bar y)}.$$
Then $\omega$ satisfies conditions (Triv1), (Triv2), and (Triv3), is convolution-invertible, and therefore defines an isomorphism of Hopf algebras over $\ZZ[i]$:
$$\alpha_\omega\colon \OO_{\ZZ[i]}(\Gamma_F) \otimes \OO_{\ZZ[i]}(\alg{\tilde G}\LVee) \isom \OO_{\ZZ[i]}({}^L \alg{\tilde G}\LVee).$$
Thus corresponding to $\psi$ and $\kappa$, there is an isomorphism of group schemes over $\ZZ[i]$:
$$\lambda_{\psi, \kappa}\colon  {}^L \alg{\tilde G}\LVee \isom \alg{\Gamma}_F \times \alg{\tilde G}\LVee.$$
\end{construction}
\proof
A convolution-inverse of $\omega$ is given by replacing $\kappa$ by $-\kappa$.  Now we check directly that $\omega$ satisfies conditions (Triv1), (Triv2), and (Triv3).  As a temporary abbreviation, we write $[\gamma]$ for $\weil(\rec(\gamma), \psi)$ throughout this proof.  From the properties of the Weil index and Hilbert symbol, we have
$$[\gamma_1 \gamma_2] [\gamma_1]^{-1} [\gamma_2]^{-1} = h(\gamma_1, \gamma_2), \quad [\gamma]^2 = h(\gamma, \gamma).$$
In this notation, the toral $\ZZ[i]$-linear map $\omega$ is given by
$$\omega(1_y^\dual)(\gamma) = [\gamma]^{\kappa(\bar y)}.$$

For (Triv1), note that $\omega(1) = 1$ since $\omega(1_0^\dual)(\gamma) = [\gamma]^{\kappa(0)} = 1$, and $\epsilon \circ \omega = \epsilon$ since $\epsilon \circ \omega(1_y^\dual) = [1]^{\kappa(\bar y)} = 1$.

For (Triv2), suppose that $a,b \in \BB$; if $\omega(a^\dual b^\dual) \neq 0$ then $a^\dual b^\dual = 1_y^\dual$ for some $y \in \tilde Y$.  From this it follows that $a^\dual = 1_u^\dual$ and $b^\dual = 1_v^\dual$ for some $u,v \in \tilde Y$ satisfying $u + v = y$.  This implies that $a_1^\dual = a_2^\dual = 1_u^\dual$ and $b_1^\dual = b_2^\dual = 1_v^\dual$ in every nonzero term on the right side of (Triv2).  Moreover, for any $a,b \in \BB$, the only nonzero terms on the right side of (Triv2) are those where $a_1^\dual = a_2^\dual = 1_u^\dual$ and $b_1^\dual = b_2^\dual = 1_v^\dual$ for some $u,v \in \tilde Y$.  Thus for (Triv2) we are left to check that
$$\omega(1_u^\dual 1_v^\dual) = \omega(1_u^\dual) \cdot \omega(1_v^\dual) \cdot \chi(1_u^\dual, 1_v^\dual).$$
Rephrasing, we must check that
$$[\gamma]^{\kappa(\bar u + \bar v)} = [\gamma]^{\kappa(\bar u)} [\gamma]^{\kappa(\bar v)} h(\gamma, \gamma)^{C(u,v)}.$$
Since $\kappa$ is a tetractor of $\bar C$, and $[\gamma]^2 = h(\gamma, \gamma)$, we have
$$[\gamma]^{\kappa(\bar u + \bar v)} = [\gamma]^{\kappa(\bar u)} [\gamma]^{\kappa(\bar v)}  [\gamma]^{[\delta \circ \bar C](\bar u, \bar v)} =  [\gamma]^{\kappa(\bar u)} [\gamma]^{\kappa(\bar v)}  h(\gamma, \gamma)^{\bar C(\bar u, \bar v)}.$$
Hence (Triv2) holds.

For (Triv3), suppose that $a \in \BB$.  Consider the result of applying comultiplication in the following order:
$$[\Delta \otimes \Delta](\Delta(a)) = \sum_{a_1, a_2, a_3, a_4} m_{a_1, a_2, a_3, a_4}^a (a_1 \otimes a_2) \otimes (a_3 \otimes a_4).$$
If $a_1^\dual = 1_{u_1}^\dual$ and $a_3^\dual = 1_{v_1}^\dual$ then $m_{a_1, a_2, a_3, a_4} = 0$ unless $a_2 \in 1_{u_1} \UU 1_{u_2}$ and $a_4 \in 1_{v_1} \UU 1_{v_2}$ for some $u_2, v_2 \in \tilde Y$, and $u_2 = v_1$.  From this it follows that 
$$u_1 \ident u_2 \ident v_1 \ident v_2, \mbox{ mod } \ZZ[\tilde \Delta^\vee].$$
Thus we write $\bar y $ for the common image of $u_1, u_2, v_1, v_2$ in $\bar Y$.

We find that the left side of (Triv3) is
$$L(a) = \sum_{a_2, a_4} m_{a_2, a_4}^a \omega(1_{u_1}^\dual) \otimes \omega(1_{v_1}^\dual) \otimes a_2 \otimes a_4.$$
Evaluation of $L(a)$ at $\gamma_1, \gamma_2 \in \Gamma_F$ yields (after a bit of reindexing)
$$L(a)(\gamma_1, \gamma_2) = \sum_{a_1, a_2} m_{a_1, a_2}^a [\gamma_1]^{\kappa(\bar y)} [\gamma_2]^{\kappa(\bar y)} (a_1 \otimes a_2).$$

A similar approach demonstrates that the right side of (Triv3) is
$$R(a) = \sum_{a_2, a_3} m_{a_2, a_3}^a \Delta(\omega(1_{u_1}^\dual)) \tau(1_{v_2}^\dual) \otimes a_2 \otimes a_3,$$
where $a_2 \in 1_{u_1} \UU 1_{u_2}$ and $a_3 \in 1_{v_1} \UU 1_{v_2}$.  Evaluation of $R(a)$ at $\gamma_1, \gamma_2 \in \Gamma_F$ yields (after a bit of reindexing)
$$R(a)(\gamma_1, \gamma_2) = \sum_{a_1, a_2} m_{a_1, a_2}^a [\gamma_1 \gamma_2]^{\kappa(\bar y)} h(\gamma_1, \gamma_2)^{Q(y)}(a_1 \otimes a_2).$$

To show that $L(a) = R(a)$, and verify (Triv3), it suffices to prove that
$$[\gamma_1]^{\kappa(\bar y)} [\gamma_2]^{\kappa(\bar y)} = [\gamma_1 \gamma_2]^{\kappa(\bar y)} \cdot h(\gamma_1, \gamma_2)^{Q(y)}.$$
This follows directly from the facts that $[\gamma_1 \gamma_2] / [\gamma_1] [\gamma_2] = h(\gamma_1, \gamma_2)$, and $Q(y) \ident \kappa(\bar y)$, mod $2$.
\qed

\section{Conjectures and Compatibilities}

We maintain the setup of the previous section.  To remind the reader, we have
\begin{itemize}
\item
A split connected reductive group $\alg{G}$ over a local field $F \not \isom \CC$; $\alg{T} \subset \alg{B} \subset \alg{G}$; the associated root system $\Psi = (Y,X)$ of type $(I, \cdot)$.  
\item
A quadratic form $Q\colon Y \rightarrow \ZZ$, positive integer $n \vert \# \mu(F)$;  the modified root datum $\tilde \Psi = (\tilde Y, \tilde X)$ and dual $\tilde \Psi\LVee = (\tilde X, \tilde Y)$ of type $(I, \checkodot)$; Lusztig's group scheme $\alg{\tilde G}\LVee$ over $\ZZ$ associated to the root datum $\tilde \Psi\LVee$, called the {\em metaplectic dual group}.
\item
A metaplectic group $\tilde G$, incarnated by $(C,\eta)$, where $C$ is a fair bisector of $Q$.  The metaplectic group fits into a short exact sequence $1 \rightarrow \mu_n(F) \rightarrow \tilde G \rightarrow G \rightarrow 1$.
\item
The L-group ${}^L \alg{\tilde G}\LVee$, depending functorially on $C$ (with $Q$,$\eta$ fixed), fitting into a short exact sequence of affine group schemes over $\ZZ$:
$$\alg{\tilde G}\LVee \rightarrow {}^L \alg{\tilde G}\LVee \rightarrow \alg{\Gamma}_F.$$
The Hopf algebra of ${}^L \alg{\tilde G}\LVee$ is obtained by a double-twist:
$$\OO({}^L \alg{\tilde G}\LVee) = \OO(\Gamma_F) \otimes_{\tau, \chi} \OO(\alg{\tilde G}\LVee).$$
\end{itemize}

Let $\epsilon\colon \mu_n(F) \rightarrow \CC^\times$ be an injective homomorphism.  An admissible representation (always on a $\CC$-vector space) $(\pi, V)$ of $\tilde G$ is called $\epsilon$-{\em genuine} if $\pi(\zeta) v = \epsilon(\zeta) \cdot v$, for all $\zeta \in \mu_n(F)$; as $\epsilon$ will be fixed, we simply say ``genuine'' to mean $\epsilon$-genuine.  We write $\Irr_\epsilon(\tilde G)$ for the set of isomorphism classes of irreducible genuine admissible representations of $\tilde G$.

\subsection{Parameters}

Let $\Weil_F$ denote the Weil group of $F$; this is a locally compact topological group, endowed with a continuous homomorphism $\iota\colon \Weil_F \rightarrow \Gamma_F$ with dense image.  We refer to Tate's excellent exposition \cite{Tat} for details.  In what follows we work with complex parameters and representations; we define ${}^L \tilde G\LVee = {}^L \alg{\tilde G}\LVee(\CC)$.  Then ${}^L \tilde G\LVee$ is a locally compact topological group whose neutral component is the complex Lie group $\tilde G\LVee = \alg{\tilde G}\LVee(\CC)$.  

\begin{definition}
A {\em Weil parameter} in the L-group ${}^L \tilde G\LVee$ is a continuous homomorphism $\rho\colon \Weil_F \rightarrow {}^L \tilde G\LVee$ such that $\rho(\Fr)$ is semisimple for a Frobenius element $\Fr \in \Weil_F$ (this condition is vacuous if $F \isom \RR$) and the composite map
$$\xymatrix{ \Weil_F \ar[r]^\rho & {}^L \tilde G\LVee \ar[r] & \Gamma_F}$$
is the map $\iota\colon\Weil_F \rightarrow \Gamma_F$.
\end{definition}
\begin{remark}
Weil parameters are called {\em infinitesimal characters} in Vogan's discussion of the local Langlands conjectures \cite{Vog}.
\end{remark}

Since ${}^L \tilde G\LVee$ is constructed through the complicated process of double-twisting the Hopf algebra $\OO(\Gamma_F) \otimes \OO(\alg{\tilde G}\LVee)$, it is difficult to work with parameters group-theoretically.  Instead, it is useful to work with parameters from the standpoint of Hopf algebras.

The Weil group is a locally compact group, and the ring $\Cont(\Weil_F)$ of continuous complex-valued functions on $\Weil_F$, with the compact open topology, is naturally a {\em complete locally convex} Hopf algebra.  The Galois group $\Gamma_F$ is also a locally compact group, and $\Cont(\Gamma_F) = \OO_\CC(\Gamma_F)$ is also a complete locally convex Hopf algebra, endowed with a continuous injective Hopf algebra homomorphism $\iota^\ast: \OO_\CC(\Gamma_F) \rightarrow \Cont(\Weil_F)$.  We refer to the Appendix for details on this category of topological Hopf algebras.  Similarly, $\OO_\CC(\alg{\tilde G}\LVee)$ can be viewed as complete locally convex Hopf algebras.  

From Proposition \ref{App}, we find that giving a continuous homomorphism $\rho$ from $\Weil_F$ to ${}^L \tilde G\LVee$ is equivalent to giving a continuous Hopf algebra homomorphism
$$\rho^\ast\colon \OO_\CC(\Gamma_F) \otimes_{\tau, \chi} \OO_\CC(\alg{\tilde G}\LVee) \rightarrow \Cont(\Weil_F).$$
For $\rho$ to give a parameter, we must have $\rho^\ast(f \otimes 1) = \iota^\ast(f)$ for all $f \in \OO_\CC(\Gamma_F) \subset \Cont(\Weil_F)$.  Since $(f \otimes a^\dual) = (f \otimes 1) \cdot_\chi (1 \otimes a^\dual)$ for all $f \in \OO(\Gamma_F)$ and all $a \in \BB$, we find that to give $\rho^\ast$, it suffices to give $\rho^\ast(1 \otimes a^\dual)$ for all $a^\dual \in \BB$.
\begin{proposition}
\label{HopfPar}
The set of Weil parameters $\rho\colon \Weil_F \rightarrow {}^L \tilde G\LVee$ is in natural bijection with the set of $\CC$-linear functions $\phi\colon \OO_\CC(\alg{\tilde G}\LVee) \rightarrow \Cont(\Weil_F)$ such that
\begin{enumerate}
\item
$\phi(1_0^\dual) = 1$ and $\epsilon(a^\dual) = \epsilon(\phi(a^\dual))$ for all $a \in \BB$.
\item
$\phi(a^\dual \cdot b^\dual) = \phi(a^\dual) \cdot \phi(b^\dual) \cdot \chi(a^\dual, b^\dual)$ for all $a,b \in \BB$.
\item
$\sum_{a_1, a_2, a_3} m_{a_1, a_2,a_3}^a \psi(a_3^\dual) \cdot [\phi(a_1^\dual) \otimes \phi(a_2^\dual)] = \Delta \phi(a^\dual)$ for all $a \in \BB$.
\item
$\phi$ factors through the restriction map from $\OO_\CC(\alg{\tilde G}\LVee)$ to the semisimple locus $\OO_\CC((\alg{\tilde G}\LVee)^{ss})$.
\end{enumerate}
\end{proposition}
\proof
Given a Weil parameter $\rho$, we write $\phi(a^\dual) = \rho^\ast(1 \otimes a^\dual)$.  Conversely, given a function $\phi$ as above, we write $\rho^\ast$ for the unique $\CC$-linear map satisfying $\rho^\ast(f \otimes a^\dual) = \iota^\ast(f) \cdot \phi(a^\dual)$ for all $f \in \OO(\Gamma_F)$ and all $a \in \BB$.  Such a $\CC$-linear map $\rho^\ast$ is automatically continuous, since $\OO_\CC(\alg{\tilde G}\LVee)$ has the inductive limit topology from its finite-dimensional subspaces (it is a DF-space), so there is an isomorphism of topological vector spaces
$$\OO(\alg{\tilde G}\LVee) \isom \bigoplus_{b \in \BB} \left( \OO(\Gamma_F) \otimes \CC b^\dual \right).$$

Conditions (1), (2), and (3) are equivalent to $\rho^\ast$ being a unital and counital algebra and coalgebra homomorphism (a homomorphism of Hopf algebras).  Condition (4) is equivalent to the semisimplicity condition on $\rho$ by Corollary \ref{ImageMap} (and the openness of the semisimple locus).
\qed  

When $F \isom \RR$, the Weil-Deligne group of $F$ is taken to be just the Weil group $\Weil_F$.  We following the treatment of Gross and Reeder \cite{G-R} to define Weil-Deligne parameters when $F$ is nonarchimedean.  For all $w \in \Weil_F$, we write $\vert w \vert = \vert \rec(w) \vert$, where $\rec\colon \Weil_F \rightarrow F^\times$ is the reciprocity map of local class field theory.  Note that when $\Fr$ is a Frobenius eleement of $\Weil_F$, $\vert \Fr \vert = q^{-1}$, where $q$ is the cardinality of the residue field of $F$.

\begin{definition}
A {\em Weil-Deligne parameter} in the L-group ${}^L \tilde G\LVee$ is a pair $(\rho, N)$, where $\rho\colon \Weil_F \rightarrow {}^L \tilde G\LVee$ is a Weil parameter in ${}^L \tilde G\LVee$ and $N$ is a nilpotent element of the Lie algebra ${}^L \Lie{\tilde g}^\vee$ of $\tilde G\LVee$, such that $\rho(w) N \rho(w)^{-1} = \vert w \vert \cdot N$, for all $w \in \Weil_F$.
\end{definition}

We say that two Weil parameters, or two Weil-Deligne parameters, are equivalent, if they are in the same $\tilde G\LVee$-orbit, where this group acts by conjugation on ${}^L \tilde G\LVee$ and by the adjoint action on its Lie algebra.  Let $\Par({}^L \tilde G\LVee)$ be the set of equivalence classes of Weil-Deligne parameters in ${}^L \tilde G\LVee$.  

The goal of this section is to provide evidence for a {\em parameterization map} -- a finite-to-one (but not necessarily surjective) map from $\Irr_\epsilon(\tilde G)$ to $\Par({}^L \tilde G\LVee)$.  More work is required to conjecture the image of this parameterization map, and to guess at its fibres.  We do not develop the theory for nonsplit forms of metaplectic groups, which would be required in any full adaptation of the local Langlands conjectures.  Nor do we develop a theory of L-packets for metaplectic groups.

\subsection{Metaplectic split tori}

Here we consider the simplest case, where $\alg{G} = \alg{T}$ is a split torus over a local field $F$ and $F \not \isom \CC$.  Recall that 
$$\tilde Y = \{ y \in Y : B(y, y') \in n \ZZ \mbox{ for all } y' \in Y \}.$$
Let $\alg{T}^\sharp$ denote the split torus over $F$ with cocharacter lattice $\tilde Y$; let $e^\sharp\colon \alg{T}^\sharp \rightarrow \alg{T}$ be the $F$-isogeny corresponding to the inclusion $\tilde Y \rightarrow Y$.  Let $e^\sharp\colon T^\sharp \rightarrow T$ be the induced homomorphism of $F$-points.

The metaplectic torus $\tilde T = \tilde T_F$ fits into a short exact sequence
$$1 \rightarrow \mu_n \rightarrow \tilde T \rightarrow T \rightarrow 1.$$
In what follows, we write $(u,v)_n = \epsilon(\Hilb_n(u,v))$ for all $u,v \in F^\times$.

\subsubsection{The center}

For $y \in Y$ and $u \in F^\times$, we have $y(u) \in T = T_F$; such elements generate $T$.  The commutator pairing associated to the central extension $\tilde T$ is given by
$$[ y_1(u_1), y_2(u_2) ] = \left( u_1, u_2 \right)_n^{B(y_1, y_2) }.$$
The center of $\tilde T$ contains $\mu_n$, and hence fits into a short exact sequence of locally compact abelian groups
$$1 \rightarrow \mu_n \rightarrow Z(\tilde T) \rightarrow Z^\dag(T) \rightarrow 1.$$
A short calculation with elementary divisors, using the nondegeneracy of the Hilbert symbol, yields the following.
\begin{proposition}[{\cite[Proposition 4.1]{Wei1}}]
Let $F$ be a local field.  Then $Z^\dag(T)$ equals the image of $T^\sharp = \alg{T}^\sharp(F)$ under the isogeny $e^\sharp$.  
\end{proposition}

We arrive at the following commutative diagram, with exact rows, each row the pullback of the one beneath it.
$$\xymatrix{
	&					&					& \Ker(e^\sharp) \ar[d] \ar[dl] & \\
1 \ar[r] & \mu_n \ar[r] \ar[d]^= & \tilde T^\sharp \ar[r] \ar[d] &  T^\sharp \ar[r] \ar[d]^{e^\sharp}& 1 \\
1 \ar[r] & \mu_n \ar[r] \ar[d]^= & Z(\tilde T) \ar[r] \ar[d] & Z^\dag(T) \ar[r] \ar[d] & 1 \\
1 \ar[r] & \mu_n \ar[r] & \tilde T \ar[r] & T \ar[r] & 1. }$$
The cover $\tilde T_F^\sharp \rightarrow T^\sharp$ splits canonically over the finite subgroup $\Ker(e^\sharp) \subset T^\sharp$.
\begin{corollary}
The continuous homomorphisms from $Z(\tilde T)$ to $\CC^\times$ are in canonical bijection with the continuous homomorphisms from $\tilde T^\sharp$ to $\CC^\times$ which are trivial on $\Ker(e^\sharp)$ (via the canonical splitting).
\end{corollary}

By a variant of the Stone von-Neumann theorem \cite[Theorem 3.1, Proposition 3.2]{Wei1}, the irreducible representations of $\tilde T$ are classified up to isomorphism by the characters of its center $Z(\tilde T)$.  
\begin{corollary}
The set $\Irr_\epsilon(\tilde T)$ is in natural bijection with the continuous genuine homomorphisms from $Z(\tilde T)$ to $\CC^\times$, which in turn are in natural bijection the continuous genuine homomorphisms from $\tilde T^\sharp$ to $\CC^\times$ trivial on $\Ker(e^\sharp)$.
\end{corollary}

Thus to parameterize $\Irr_{\epsilon}(\tilde T)$, we must parameterize $\Hom_{\epsilon}(\tilde T^\sharp, \CC^\times)$ (continuous homomorphisms which restrict to $\epsilon$ on $\mu_n$).  The group $\tilde T^\sharp$ admits the following presentation.
\begin{proposition}
The group $\tilde T^\sharp$ is generated by elements $\zeta \in \mu_n$, and elements $\tilde y(u)$ for all $y \in \tilde Y$ and all $u \in F^\times$, subject to the relations:
\begin{enumerate}
\item
all elements of $\tilde T^\sharp$ commute;
\item
for $y_1, y_2 \in \tilde Y$ and $u \in F^\times$, $\tilde y_1(u) \tilde y_2(u) = \widetilde{y_1 y_2}(u) \cdot (u,u)_n^{C(y_1, y_2)}$;
\item
for $y \in \tilde Y$ and $u_1, u_2 \in F^\times$, $\tilde y(u_1) \tilde y(u_2) = \tilde y(u_1, u_2) \cdot (u_1, u_2)_n^{Q(y)}$.
\end{enumerate}
\end{proposition}
\proof
These relations are consequences of the construction of $\tilde T$ from $(Q,C)$ -- the incarnation of split metaplectic tori.  Choosing a basis $y_1, \ldots, y_r$ of the free $\ZZ$-module $\tilde Y$, every element of $\tilde T$ can be written uniquely in the form
$$\tilde y_1(u_1) \tilde y_2(u_2) \cdots \tilde y_r(u_r) \cdot \zeta,$$
for $u_1, \ldots, u_r \in F^\times$ and $\zeta \in \mu_n$ (using Relations 1 and 2).  Relations 1 and 3 suffice to multiply elements of $\tilde T$ in such a form.  Hence the relations suffice for a presentation of $\tilde T$.
\qed

\subsubsection{Parameters}
We have found a natural map from $\Irr_\epsilon(\tilde T)$ to $\Hom_{\epsilon}(\tilde T^\sharp, \CC^\times)$, whose fibres have cardinality zero or one.  Now we demonstrate a Langlands-style parameterization using our metaplectic L-group.
\begin{thm}
There is a natural bijection between $\Hom_{\epsilon}(\tilde T^\sharp, \CC^\times)$ and $\Par({}^L \tilde T\LVee)$.
\end{thm}
\proof
First, note that a Weil parameter is the same as a Weil-Deligne parameter, since the complex Lie group $\tilde T\LVee$ has no nilpotent elements.  Moreover, since every element of ${}^L \tilde T\LVee$ is semisimple, a Weil parameter is simply a continuous homomorphism $\rho\colon \Weil_F \rightarrow {}^L \tilde T\LVee$ which intertwines the maps to $\Gamma_F$.

Recall the L-group of the metaplectic torus is a group scheme over $\ZZ$:
$${}^L \alg{\tilde T}\LVee = \Spec \left( \OO(\Gamma_F) \otimes_{\tau,\chi} \OO(\alg{\tilde T}\LVee) \right),$$
where $(\tau, \chi)$ is the compatible pair:
$$\tau(y)(\gamma_1, \gamma_2) = h(\gamma_1, \gamma_2)^{Q(y)}, \quad \chi(y_1, y_2)(\gamma) = h(\gamma, \gamma)^{C(y_1, y_2)}.$$
Here the Hopf algebra $\OO(\alg{\tilde T}\LVee)$ is the commutative {\em and cocommutative} Hopf algebra $\ZZ[\tilde Y]$ (with comultiplication $\Delta(y) = y \otimes y$), so we write $y$ instead of $1_y^\dual$.  Also, recall that $h(\gamma_1, \gamma_2) = (\rec(\gamma_1), \rec(\gamma_2))_n$.

By Proposition \ref{HopfPar}, giving a Weil parameter $\rho$ is equivalent to giving a $\CC$-linear map $\phi\colon \CC[\tilde Y] \rightarrow \Cont(\Weil_F)$ satisfying the following, for all $y, y_1, y_2 \in \tilde Y$, and all $w, w_1, w_2 \in \Weil_F$:
\begin{enumerate}
\item
$\phi(0; w) = 1$ and $\phi(y; 1) = 1$;
\item
$\phi(y_1 y_2; w) = \phi(y_1; w) \phi(y_2; w) h(w,w)^{C(y_1, y_2)}$;
\item
$\phi(y; w_1 w_2) = \phi(y; w_1) \phi(y; w_2) h(w_1, w_2)^{Q(y)}$.
\end{enumerate}
For any such $\phi$, and any $y \in \tilde Y$, the function $\phi(y; \bullet)\colon \Weil_F \rightarrow \CC$ factors through the abelianization $\Weil_F^{ab}$ by (1) and (3).  Thus to give a parameter, it is necessary and sufficient (via the reciprocity isomorphism of local class field theory) to give functions $\phi\colon Y \times F^\times \rightarrow \CC$ satisfying the following, for all $y, y_1, y_2 \in \tilde Y$ and all $u, u_1, u_2 \in F^\times$:
\begin{enumerate}
\item
$\phi(0; u) = 1$ and $\phi(y; 1) = 1$.
\item
$\phi(y_1 y_2; u) = \phi(y_1; u) \cdot \phi(y_2; u) \cdot (u,u)_n^{C(y_1, y_2)}$.
\item
$\phi(y; u_1 u_2) = \phi(y; u_1) \cdot \phi(y; u_2) \cdot (u_1, u_2)_n^{Q(y)}$.
\end{enumerate}
By the presentation of $\tilde T^\sharp$, this data is precisely what is needed for $\tilde y(u) \mapsto \phi(y; u)$ to give a genuine homomorphism from $\tilde T^\sharp$ to $\CC^\times$.
\qed

In this way, we identify a one-to-one correspondence between $\Irr_\epsilon(\tilde T)$ -- the irreducible genuine representations of $\tilde T$, up to equivalence -- and a subset of $\Par({}^L \tilde T\LVee)$.  This subset can be identified as well, using the isogeny of $F$-tori corresonding to the inclusion $n Y \hookrightarrow Y^\sharp$, but we leave this for a later paper.
  
\subsection{Correspondences for $Mp_{2n}$}

Consider the classical metaplectic group over the local field $F$, with $F \neq \CC$ and $\Char(F) \neq 2$:
$$1 \rightarrow \mu_2 \rightarrow Mp_{2n} \rightarrow Sp_{2n} \rightarrow 1$$
In other words, we consider the case where $\alg{G} = \alg{Sp}_{2n}$ over $F$, $\alg{T} \subset \alg{B}$ an $F$-split maximal torus contained in a Borel subgroup over $F$.  Let $\{ \alpha_1, \ldots, \alpha_n \}$ denote the set of simple roots, with $\alpha_1$ the unique long simple root.  Thus $\alpha_1^\vee$ is the unique short simple coroot.  The quadratic form $Q$ corresponding to the metaplectic group $Mp_{2n}$ is the unique Weyl-invariant one for which $Q(\alpha_1^\vee) = 1$ and $Q(\alpha_i^\vee) = 2$ for $2 \leq i \leq n$.

Here $Y = \ZZ[\Delta^\vee]$ (since $\alg{G}$ is simply-connected) and $B(\alpha_i^\vee, \alpha_j^\vee) \in 2 \ZZ$ for all $1 \leq i,j \leq n$.  Hence $\tilde Y = Y$ and $\tilde X = X$.  Moreover, $\tilde \alpha_1^\vee = 2 \alpha_1^\vee$ and $\tilde \alpha_i^\vee = \alpha_i^\vee$ for $2 \leq i \leq n$.  It follows that $\alg{\tilde G}\LVee \isom \alg{Sp}_{2n}$ (the Chevalley group over $\ZZ$) in this setting.  Moreover, $\bar Y = \tilde Y / (2 Y + \ZZ[\tilde \Delta^\vee]) \isom \ZZ / 2 \ZZ$, generated by the image of $\alpha_1^\vee$.  Hence, for any nontrivial continuous additive character $\psi\colon F \rightarrow \CC^\times$ and any tetractor $\kappa\colon \bar Y \rightarrow \ZZ / 4 \ZZ$ as described before, we have an isomorphism of group schemes over $\ZZ[i]$:
$$\lambda_{\psi, \kappa}:  {}^L \alg{\tilde G}\LVee \isom \alg{\Gamma}_F \times \alg{Sp}_{2n}.$$
Moreover, as there are exactly two tetractors, they have the form $\pm \kappa$, and changing between tetractors is equivalent to changing $\psi$ to $\bar \psi$.  Hence the isomorphisms of group schemes over $\ZZ[i]$ for the two tetractors are related by the ``conjugation'' sending $i$ to $-i$.  In other words, the choice of tetractor is equivalent to a choice of primitive fourth root of unity in $\CC$.  

Thus our construction of an L-group for $Mp_{2n}$ is compatible with the following well-known conjecture.  
\begin{conj}
For each nontrivial continuous additive character $\psi$ of $F$, there is a finite-to-one parameterization map from $\Irr_\epsilon(Mp_{2n})$ to $\Par(\alg{\Gamma}_F \times \alg{Sp}_{2n})$.  
\end{conj}

Waldspurger proves this conjecture for $n = 1$, accepting the local Langlands correspondence for $G = PGL_2(F), G^\vee = SL_2(\CC)$, in his representation-theoretic Shimura correspondence \cite{Wal1}.  Indeed, Waldspurger constructs a local ``Shimura correspondence'', from $\Irr_\epsilon(Mp_{2n})$ to $\Irr(PGL_2)$, and the L-group of $PGL_2$ is $\alg{\Gamma}_F \times \alg{SL}_2$.  Moreover, this correspondence is finite-to-one (the fibres are related by the ``Waldspurger involution'', when working over a $p$-adic field with $p \neq 2$).

More generally, in the archimedean setting, Adams and Barbasch \cite{A-B} have exhibited a Shimura-like correspondence 
$$\Irr_\epsilon(Mp_{2n}) \rightarrow \bigcup_{p+q=2n+1} \Irr(SO(p,q)).$$
Again, this correspondence depends on a choice of nontrivial continuous additive character of $\RR$.  By the local Langlands parameterization for $SO(p,q)$ over $\RR$, this verifies the conjecture written above, in the archimedean case.  

Finally, when $F$ is a $p$-adic field with $p$ odd, the recent work of Gan and Savin \cite[Corollary 1.2]{GS1} (see also the fundamental work of Kudla-Rallis \cite{K-R}) gives a finite-to-one map from $\Irr_\epsilon(Mp_{2n})$ to $\Par(\alg{\Gamma}_F \times \alg{Sp}_{2n})$ as conjectured above, contingent upon the local Langlands correspondence for special orthogonal groups over $p$-adic fields.  

Our construction of an L-group of $Mp_{2n}$ is compatible with the above conjecture, not only in the parameterization by $\Par(\alg{\Gamma}_F \times \alg{Sp}_{2n})$, but also by the choices necessary for such a parameterization.  A fourth root of unity $i \in \CC$, and an additive character $\psi$, suffice in our conjectural parameterization as well as in the work of Waldspurger, Adams-Barbasch, and Gan-Savin.  As observed in \cite[Theorem 5.3]{GS1}, the choice of additive character only matters up to $F^{\times 2}$-orbit, compatible with our findings on isomorphisms of L-groups over $\ZZ[i]$.
\begin{remark}
It would be interesting to see if the $\epsilon$-dichotomy studied in \cite{GS1} can also be seen on the ``Galois side'' via our L-group.
\end{remark}  

\subsection{Unramified correspondence}

Let $F$ be a $p$-adic field, with valuation ring $R$ and residue field $\FF_q$ of cardinality $q$.   Let $\Inertia_F \subset \Gamma_F$ denote the inertia subgroup.  Let $\varpi$ be a uniformizing parameter for $F$ (the choice will not matter).  Let $\alg{G}$ be a split reductive group over $F$ as before.  Assume $n=2$ (the case with $n$ odd is less interesting).  Recall that to each choice of continuous nontrivial additive character $\psi$ and tetractor $\kappa$, we have found an isomorphism $\lambda_{\psi, \kappa}$ of group schemes making the following triangle commute.
$$\xymatrix{
{}^L \alg{\tilde G}\LVee \ar[rr]^{\lambda_{\psi, \kappa}} \ar[dr] & & \alg{\Gamma}_F \times \alg{\tilde G}\LVee \ar[dl] \\
& \alg{\Gamma}_F & }$$

Let $\rho\colon \Weil_F \rightarrow {}^L \tilde G\LVee$ be a Weil parameter.  Define $\rho_{\psi, \kappa}$ to be the composition:  $\rho_{\psi, \kappa} = \lambda_{\psi, \kappa} \circ \rho$.  Let $\phi\colon \OO_\CC(\alg{\tilde G}\LVee) \rightarrow \Cont(\Weil_F)$ be the $\CC$-linear map corresponding to the Weil parameter $\rho$ by Proposition \ref{HopfPar}.  Let $pr_2: \alg{\Gamma}_F \times \alg{\tilde G}\LVee \rightarrow \alg{\tilde G}\LVee$ denote the projection onto the second factor.  

The map of Hopf algebras, $[pr_2 \circ \rho_{\psi, \kappa}]^\ast: \OO_\CC(\alg{\tilde G}\LVee) \rightarrow \Cont(\Weil_F)$ is given, for all $a \in \BB$ with $y = \lt{a}$, by
\begin{eqnarray*}
(pr_2 \circ \rho_{\psi, \kappa})^\ast(a^\dual) & = & [\rho^\ast \circ \lambda_{\psi, \kappa}^\ast \circ pr_2^\ast](a^\dual) \\
& = &  [\rho^\ast \circ \lambda_{\psi, \kappa}^\ast](1 \otimes a^\dual) \\
& = & [\rho^\ast](\weil(\rec(\bullet), \psi)^{\kappa(\bar y)} \otimes a^\dual) \\
& = & \weil(\rec(\bullet), \psi)^{\kappa(\bar y)} \cdot \phi(a^\dual).
\end{eqnarray*}
Here we have used Construction \ref{TrivConstruct} to prove that
\begin{equation}
\label{EQunr}
(pr_2 \circ \rho_{\psi, \kappa})^\ast(a^\dual) = \weil(\rec(\bullet), \psi)^{\kappa(\bar y)} \cdot \phi(a^\dual).
\end{equation}

To keep track of these Hopf algebras and homomorphisms over $\CC$, we have the following commutative diagram.
\begin{equation}
\label{BigCD}
\xymatrix{
& \Cont(\Weil_F) & & \OO_\CC(\alg{\tilde G}\LVee) \ar[dl]_{pr_2^\ast} \ar[ll]_{\weil(\rec(\bullet), \psi)^{\kappa(\bar y)}  \cdot \phi^\ast} \\
\OO_\CC({}^L \alg{\tilde G}\LVee) \ar[ur]^{\rho^\ast}  & & \OO_\CC(\alg{\Gamma}_F) \otimes \OO_\CC(\alg{\tilde G}\LVee) \ar[ul]_{\rho_{\psi, \kappa}^\ast} \ar[ll]_{\lambda_{\psi, \kappa}^\ast}  \\
& \OO_\CC(\alg{\Gamma}_F) \ar[ur] \ar[ul]& }
\end{equation}

\begin{definition}
We say that $\rho$ is {\em unramified} with respect to $(\psi, \kappa)$, or simply that $\rho_{\psi, \kappa}$ is unramified, if $\rho_{\psi, \kappa}(w) = (w,1)$ for all $w \in \Inertia_F$.  A Weil-Deligne parameter $(\rho, N)$ will be called unramified with respect to $(\psi, \kappa)$ if its Weil parameter $\rho_{\psi, \kappa}$ is unramified.
\end{definition}

When $\rho_{\psi, \kappa}$ is unramified, we write $g_{\psi, \kappa} = g_{\psi, \kappa}(\rho)$ for the unique element of $\tilde G\LVee$ such that $\rho_{\psi, \kappa}(\Fr) = (\Fr, g)$ for every (geometric) Frobenius element $\Fr \in \Weil_F$.  Then, when $\rho_{\psi, \kappa}$ is unramified we find that it is completely determined by the composition 
$$pr_2 \circ \rho_{\psi, \kappa}: \Weil_F / \Inertia_F \rightarrow \tilde G\LVee, \quad \Fr^d \mapsto g_{\psi, \kappa}(\rho)^d.$$

Note that $g_{\psi, \kappa}(\rho) \in \alg{\tilde G}\LVee(\CC)$; evaluation of the regular function $a^\dual$ at this $\CC$-point of  $\alg{\tilde G}\LVee$ is given by
$$a^\dual(g_{\psi, \kappa}(\rho)) = a^\dual \left( [pr_2 \circ \rho_{\psi, \kappa}](\Fr) \right) =  \weil(\rec(\Fr), \psi)^{\kappa(\bar y)} \cdot \phi(a^\dual)(\Fr).$$

We can rephrase the ``unramified'' condition in terms of Hopf algebras.  This allows us to see a relationship between Hecke algebras and Hopf algebras, and to see how the unramified condition depends on the choice $(\psi, \kappa)$.
\begin{proposition}
The Weil parameter $\rho$ is unramified with respect to $(\psi, \kappa)$ if and only if
$$\weil(\rec(\bullet), \psi)^{\kappa(\bar y)} \cdot \phi(a^\dual) \in \Cont(\Weil_F / \Inertia_F),$$
for all $a \in \BB$, with $y = \lt{a} \in \tilde Y$.
\end{proposition}
\proof
We find that $\rho_{\psi, \kappa}$ is unramified if and only if $pr_2 \circ \rho_{\psi, \kappa}\colon \Weil_F \rightarrow \tilde G\LVee$ factors through $\Weil_F / \Inertia_F$.  This occurs if and only if, for all $h \in \Cont(\tilde G\LVee)$, $(pr_2 \circ \rho_{\psi, \kappa})^\ast(h) \in \Cont(\Weil_F / \Inertia_F)$.  By density of $\OO_\CC(\alg{\tilde G}\LVee)$ in $\Cont(\tilde G\LVee)$, this occurs if and only if $(pr_2 \circ \rho_{\psi, \kappa})^\ast(a^\dual) \in \Cont(\Weil_F / \Inertia_F)$ for all $a \in \BB$.  

Now the result follows from Equation \ref{EQunr} (or see the commutative diagram \ref{BigCD}).
\qed

\begin{thm}
Suppose that $p$ is odd.  Suppose that $\rho$ is a Weil parameter and $\rho_{\psi, \kappa}$ is unramified.  If $v \in R^\times$ and $\psi' = {}^v \psi$, then $\rho_{\psi', \kappa}$ is also unramified.

Let $\bar v$ be the image of $v \in R^\times$ in the residue field $\FF_q^\times$ and let $\Leg$ denotes the Legendre symbol for $\FF_q^\times$:
$$\Leg(\bar v) = \begin{cases} 1 & \mbox{ if } \bar v \in \FF_q^{\times 2}\\ -1 & \mbox{ if } \bar v \not \in \FF_q^{\times 2}. \end{cases} $$

Recall that $\alg{Z}(\alg{\tilde G}\LVee) = \Spec(\ZZ[\tilde Y] / \langle \tilde \alpha_i^\vee : i \in I \rangle$.  Let $z_Q(v)  \in  \alg{Z}(\alg{\tilde G}\LVee)(\ZZ)$ be the central element corresponding to the ring homomorphism
$$\frac{\ZZ[\tilde Y]}{\langle \tilde \alpha_i^\vee : i \in I \rangle} \rightarrow \ZZ, \quad y \mapsto \Leg(v)^{Q(y)}.$$

Then the elements $g_{\psi, \kappa}$ and $g_{\psi', \kappa}$ satisfy
$$g_{\psi', \kappa} = z_Q(v) \cdot g_{\psi, \kappa}.$$
\end{thm}
\proof
Since $p$ is odd, we find that for all $u,v \in R^\times$, $\Hilb_2(u,v) = 1$.   Therefore $\Hilb_2(\rec(\bullet), v) \in \Cont(\Weil_F / \Inertia_F)$.  For all continuous nontrivial additive characters $\psi$, and all $v \in R^\times$,
$$\weil(\rec(\bullet), {}^v \psi)^{\kappa(\bar y)} \cdot \phi(a^\dual) = \Hilb_2(\rec(\bullet), v) \cdot \weil(\rec(\bullet),\psi)^{\kappa(\bar y)} \cdot \phi(a^\dual).$$
It follows that 
$$\weil(\rec(\bullet), \psi)^{\kappa(\bar y)} \cdot \phi(a^\dual) \in \Cont(\Weil_F / \Inertia_F) \Leftrightarrow \weil(\rec(\bullet), {}^v \psi)^{\kappa(\bar y)} \cdot \phi(a^\dual) \in \Cont(\Weil_F / \Inertia_F).$$

Again since $p$ is odd, we find that
$$\Hilb_2(\rec(\Fr), v) = \Leg(\bar v),$$
where $\bar v$ is the image of $v \in R^\times$ in the residue field $\FF_q^\times$.
Now for $\psi' = {}^v \psi$, and any $a \in \BB$, we have 
$$a^\dual(g_{\psi', \kappa}) = a^\dual(g_{\psi, \kappa}) \cdot \Leg(\bar v)^{\kappa(\bar y)}.$$
Since $\kappa(\bar y) \ident Q(y)$, mod $2$, we find that
$$a^\dual(g_{\psi', \kappa}) = a^\dual(g_{\psi, \kappa}) \cdot \Leg(\bar v)^{Q(y)}.$$
The function $y \mapsto \Leg(\bar v)^{Q(y)}$ defines the central element $z = z_Q(v) \in \alg{Z}(\alg{\tilde G}\LVee)(\ZZ)$.  

Moreover, for any $a \in \BB$, with $y = \lt{a}$,
\begin{eqnarray*}
a^\dual(z g_{\psi, \kappa}) & = & \sum_{a_1, a_2} m_{a_1, a_2}^a  a_1^\dual(z) a_2^\dual(g_{\psi, \kappa}) \\
& = & \Leg(\bar v)^{Q(y)} \cdot a^\dual(g_{\psi, \kappa}) \\
& = & a^\dual(g_{\psi', \kappa}).
\end{eqnarray*}
Hence $g_{\psi', \kappa} = z \cdot g_{\psi, \kappa}$ as desired.
\qed

When $\psi$ is a nontrivial continuous additive character of $F$, the {\em conductor} of $\psi$ is the largest integer $m$ for which $\psi(x) = 1$ for all $x \in \varpi^{-m} R$.  
\begin{corollary}
\label{UnrChoices}
Suppose that $p$ is odd and $n=2$ as before.  Then the answer to the question ``Is $\rho_{\psi, \kappa}$ unramified?'' depends only upon $\kappa$ and the parity of the conductor of $\psi$.  Moreover, for any fixed $\kappa$ and parity ($0$ or $1$), the element $g_{\psi, \kappa}$ depends only on $\weil(\varpi, \psi)$ and $\rho$.
\end{corollary}
\proof
We have seen that the unramified condition is preserved under twisting the additive character by $R^\times$.  Since the Weil index is unchanged after twisting the additive character by squares, we find that the unramified condition is preserved under twisting the additive character by $R^\times \times \varpi^{2 \ZZ}$.  Hence the unramified condition depends only on the parity of the conductor of $\psi$.

The element $g_{\psi, \kappa}$ is uniquely determined by $\rho$, $\psi$, and $\kappa$, via the formula
$$a^\dual(g_{\psi, \kappa}) = \weil(\varpi, \psi)^{\kappa(\bar y)} \cdot \phi(a^\dual)(\Fr),$$
where $\Fr$ is any Frobenius element satisfying $\varpi = \rec(\Fr)$.  Hence knowing the Weil index $\weil(\varpi, \psi)$ together with $\rho$ suffices to determine $g_{\psi, \kappa}$ uniquely.
\qed

\begin{remark} 
By twisting $\psi$, one can force $\weil(\varpi, \psi)$ to be either square root of $\Hilb_2(\varpi, \varpi)$.  Hence ``knowing the Weil index $\weil(\varpi, \psi)$'' boils down to choosing a square root of $-1$ if $\Hilb_2(\varpi, \varpi) = -1$.  This choice of square root of $\Hilb_2(\varpi, \varpi)$ can be found in \cite{Sav}.  On the other hand, if $q \ident 1$, mod $4$, or equivalently, if $\Hilb(\varpi, \varpi) = 1$, then one can choose $\psi$ in such a way that $\weil(\varpi, \psi) = 1$.  This simplification reduces the number of choices necessary in the work of McNamara \cite{McN}.  
\end{remark}
 
\subsubsection{Compatibility with Savin's correspondence}
Now let $\alg{G}$ be a simply-connected, simply-laced (not type $A_1$), simple split algebraic group over the $p$-adic field $F$.  In \cite{Sav}, Savin studies the ``unramified'' representations of a central extension $\tilde G$:
$$1 \rightarrow \mu_n \rightarrow \tilde G \rightarrow G \rightarrow 1.$$
This covering corresponds to the Weyl-invariant quadratic form $Q\colon Y \rightarrow \ZZ$ satisfying $Q(\alpha_i^\vee) = 1$ for all simple coroots $\alpha_i^\vee$.

Savin's construction of the cover $\tilde G$ in \cite[\S 2]{Sav} inspired our definition of a bisector; his construction relies on an ordering of the simple coroots (which generate the $\ZZ$-module $Y$).  

In \cite[Theorem 7.8]{Sav}, Savin finds an isomorphism from the Iwahori Hecke algebra of $\tilde G$ to the Iwahori Hecke algebra of $G_n$, where $\alg{G}_n$ is the algebraic group obtained as the quotient of $\alg{G}$ by the $n$-torsion subgroup of its center.  This is precisely the group with root datum $\tilde \Psi$, and thus suggests an ``unramified parameterization'':
$$\Irr_\epsilon^{unr}(\tilde G) \rightarrow \Par^{unr}({}^L \alg{\tilde G}\LVee).$$

In order to define this isomorphism of Iwahori Hecke algebras, no choices are necessary if $n$ is odd.  But if $n = 2$, Savin makes a crucial choice in \cite[\S 4]{Sav}: he requires a function $\gamma\colon \bar Y \rightarrow \mu_4$ such that
$$\gamma(\bar y_1) \gamma(\bar y_2) = \Hilb_2(\varpi, \varpi)^{C(y_1, y_2)} \gamma(\bar y_1 + \bar y_2).$$
While Savin makes such a choice using an independently interesting observation \cite[Proposition 4.2]{Sav}, we observe that the choice of a tetractor together with a choice of a fourth root of unity (if $\Hilb_2(\varpi, \varpi) = -1$) is also sufficient.  Indeed, one does not need the entire data of an additive character, if one only cares about unramified representations.  In this way, Savin's correspondence is compatible with what our L-group suggests, not only in the ``dual group'' $\alg{\tilde G}\LVee$, but also in the choices necessary in Corollary \ref{UnrChoices}.

\subsubsection{Compatibility with McNamara's parameterization}

In \cite{McN}, McNamara works with a fairly general class of split reductive groups $\alg{G}$ over local nonarchimedean fields $F$.  McNamara makes a ``mild'' assumption on $n$, which ends up drastically simplifying the L-group and parameterization from our perspective.  Namely, McNamara assumes that $p$ does not divide $n$, and that $F^\times$ contains $2n$ distinct $2n^{th}$ roots of unity.  The latter assumption implies that $(u,u)_n = (u,u)_{2n}^2 = 1$ for all $u \in F^\times$.  This entirely trivializes our ``second twist'' $\chi$ in the construction of the L-group.  The double-twist becomes a single-twist and the L-group ${}^L \alg{\tilde G}\LVee$ depends only on the quadratic form $Q$ and not on the bisector $C$.

Under McNamara's assumptions, our tetractor $\kappa$ only matters mod $2$ instead of mod $4$ (it always occurs in an exponent of $\pm 1$), and since $\kappa(y) \ident Q(y)$, mod $2$, the choice of tetractor also becomes irrelevant.  

Fixing the parity of an additive character $\psi$ to be zero, we find by Corollary \ref{UnrChoices} that the question ``Is $\rho_{\psi, \kappa}$ unramified?'' has answer independent of any other choices.  Fixing a uniformizer $\varpi$ and additive character $\psi$ satisfying $\weil(\varpi, \psi) = 1$, no further choices are necessary to give a bijection between the set of unramified Weil parameters and the set of elements of $\tilde G\LVee$.

This is compatible with the metaplectic Satake isomorphism \cite[Theorem 10.1]{McN}, and the isomorphism of Iwahori Hecke algebras of \cite[Corollary 13.1]{McN}.  This also explains why McNamara's assumptions on $F$ and $n$ allow him to make so few choices in order to determine such Hecke algebra isomorphisms.
\appendix
\section{Locally convex Hopf algebras}

When $G$ is a Hausdorff topological group, let $\Cont(X)$ be the ring of continuous functions from $X$ to $\CC$, with pointwise addition and multiplication.  If $G$ is a finite group (with the discrete topology), then the group structure on $G$ yields a Hopf algebra structure on $\Cont(G)$.  For generalization to locally compact groups, the most popular approaches involve $C^\ast$-algebras or von Neumann algebras.  The definition of a Hopf $C^\ast$-algebra requires some awkward bouncing between different spaces of functions and measures; we refer to \cite{V-V} for a survey.

The approach via von Neumann, or $W^\ast$, or Kac algebras suffers fewer complications than the $C^\ast$-algebra approach.  But still, defining the counit in this setting presents difficulties, typically solved with a ``Haar weight''.  We refer to Majid \cite{Maj2} for a discussion of Hopf algebras in this setting.  Even with these difficulties, the theory of Hopf algebras in the setting of $C^\ast$-algebras or von Neumann algebras has enjoyed successes in the past decade, particularly in generalizing some harmonic analysis and the theory of unitary representations.

For the purposes of this paper, it is most natural to consider the locally convex topological vector space $\Cont(G)$ of continuous functions from $G$ to $\CC$, with the compact-open topology.  Using a bit of functional analysis (particularly, the completed injective tensor product), most of which has been around for over 40 years, we can view $\Cont(G)$ as a Hopf algebra in a symmetric monoidal category of complete locally convex Hausdorff topological vector spaces over $\CC$.

\subsection{Topological algebras}

In what follows, we assume some functional analysis over the complex numbers; we refer to K\"othe \cite{Kot} for relevant information about topological vector spaces.  Let $\C{Vec}$ denote the category of locally convex, Hausdorff, complete topological vector spaces over $\CC$ and continuous linear maps.  All vector spaces here will be over $\CC$, and all linear maps will be $\CC$-linear.  For $V \in \C{Vec}$, we write $V'$ for the vector space of continuous linear maps from $V$ to $\CC$.

The category $\C{Vec}$ is additive with respect to the ordinary direct sum of vector spaces.  It is a symmetric monoidal category, with respect to the completed injective tensor product (sometimes called the completed $\epsilon$-tensor product) defined below.  
\begin{definition}[{\cite[Ch.\ 7, \S 44.1]{Kot}}]
For $V,W \in \C{Vec}$, let $V \hat \otimes W$ denote the completion of $V \otimes W$ with respect to the weakest topology for which
\begin{enumerate}
\item
The canonical map from $V \times W$ to $V \otimes W$ is separately continuous;
\item
For all $\lambda \in V'$ and $\mu \in W'$, the element $\lambda \otimes \mu$ is in $(V \otimes W)'$, i.e., $\lambda \otimes \mu$ is continuous;
\item
If $\Lambda \subset V'$ is equicontinuous on $V$ and $M \subset W'$ is equicontinuous on $W$, then $\Lambda \otimes M \subset V' \otimes W'$ is equicontinuous on $V \otimes W$.
\end{enumerate}
\end{definition}

Let $\Cat{Top}_{lc}$ be the category of locally compact Hausdorff spaces and continuous maps.  For $X \in \Cat{Top}_{lc}$, let $\Cont(X)$ denote the space of continuous $\CC$-valued functions on $X$, with the compact-open topology.  Then $\Cont(X) \in \C{Vec}$.  More generally, when $E \in \C{Vec}$, we write $\Cont(X,E)$ for the space of continuous $E$-valued functions on $X$, with the compact-open topology.

When $f \colon X \rightarrow Y$ is a morphism in $\Cat{Top}_{lc}$, composition with $f$ yields a continuous $\CC$-linear homomorphism $\Cont(f)\colon \Cont(Y) \rightarrow \Cont(X)$.  In this way, we define a contravariant functor 
$$\Cont\colon \Cat{Top}_{lc} \rightarrow \C{Vec}.$$

The following proposition is a direct consequence of \cite[Ch.\ 8, \S 44.7]{Kot} and can also be found in \cite[Ch.\ XI, \S 1, Corollary 1.1]{Mal}.
\begin{proposition}
Let $X$ and $Y$ be locally compact Hausdorff topological spaces.  Consider the $\CC$-bilinear map
$$\odot\colon \Cont(X) \times \Cont(Y) \rightarrow \Cont(X \times Y), \quad [f \odot g](x,y) = f(x) \cdot g(y).$$
Then $\odot$ extends uniquely to an isomorphism in $\C{Vec}$ (i.e., a continuous isomorphism with continuous inverse):
$$\hat \odot\colon \Cont(X) \hat \otimes \Cont(Y) \rightarrow \Cont(X \times Y).$$
\end{proposition}

\begin{corollary}
Pointwise multiplication, a priori a $\CC$-bilinear map
$$\cdot\colon \Cont(X) \times \Cont(X) \rightarrow \Cont(X),$$ 
extends uniquely to a continuous $\CC$-linear map:
$$\mu\colon \Cont(X) \hat \otimes \Cont(X) \rightarrow \Cont(X).$$
\end{corollary}
\proof
The previous proposition yields a continuous isomorphism in $\C{Vec}$:
$$\hat \odot\colon \Cont(X) \hat \otimes \Cont(X) \rightarrow \Cont(X \times X).$$
The diagonal inclusion $d\colon X \hookrightarrow X \times X$ yields a continuous pullback:
$$\Cont(d)\colon \Cont(X \times X) \rightarrow \Cont(X).$$
Observing that
$$[\Cont(d) \circ \odot](f,g) (x) = f(x) g(x),$$
the desired map $\mu$ is $\Cont(d) \circ \hat \odot$.
\qed

We write $\One$ for the element of $Hom(\CC,\Cont(X))$ satisfying $[\One(z)](x) = z$ for all $z \in \CC$ and $x \in X$.  The triple $(\Cont(X), \mu, \One)$ is our first example of a complete algebra.
\begin{definition}
A {\em complete algebra} is an algebra (always associative with unit) in the symmetric monoidal category $\C{Vec}$, i.e., a triple $(\OO, \mu, \One)$, with $\OO \in \C{Vec}$, and where 
$$\mu\colon \OO \hat \otimes \OO \rightarrow \OO, \quad \One\colon \CC \rightarrow \OO$$
are morphisms in $\C{Vec}$ making the diagrams for unit and associativity commute.  A morphism of complete algebras means a continuous unital $\CC$-algebra homomorphism.  We write $\C{Alg}$ for the resulting category of complete algebras.
\end{definition}

The previous corollary implies that the space $\Cont(X)$ becomes a commutative complete algebra, with respect to pointwise multiplication of functions (extended to $\mu$) and unit morphism $\One$ sending constants to constant functions.
\begin{proposition}
The functor $\Cont\colon \Cat{Top}_{lc} \rightarrow \C{Alg}$ is fully faithful.
\label{FFF}
\end{proposition}
\proof
Suppose that $X$ and $Y$ are locally compact Hausdorff spaces, and $\phi\colon \Cont(Y) \rightarrow \Cont(X)$ is a continuous ring homomorphism.  Let $\upsilon Y$ be the {\em realcompactification} of $Y$, as defined by Hewitt \cite{Hew}; $Y$ is a dense subspace of $\upsilon Y$ and there is an associated (algebraic, but not necessarily continuous) $\CC$-algebra isomorphism $\upsilon\colon \Cont(Y) \rightarrow \Cont(\upsilon Y)$ (\cite[Theorem 56]{Hew}).  If $\Res_Y$ denotes the continuous restriction map from $\Cont(\upsilon Y)$ to $\Cont(Y)$, then
$$\Res_Y \circ \upsilon = \Id\colon \Cont(Y) \rightarrow \Cont(Y).$$
In other words every $g \in \Cont(Y)$ extends uniquely to a continuous $\upsilon g \in \Cont(\upsilon Y)$.

By \cite[Theorem 10.8]{G-J}, there exists a unique continuous function $f\colon X \rightarrow \upsilon Y$ such that for all $x \in X$, and all $g \in \Cont(Y)$, $[\upsilon g](f(x)) = [\phi(g)](x)$.
This can be rephrased as an equality of continuous $\CC$-algebra homomorphisms:
$$\ev_{f(x)} \circ \upsilon = \ev_x \circ \phi \in \Hom_{\C{Alg}}(Y, \CC).$$
But every element of $\Hom_{\C{Alg}}(Y, \CC)$ coincides with $\ev_y$ for a unique $y \in Y$, by \cite[Theorem 2]{Edw}; thus $\ev_{f(x)} \circ \upsilon = \ev_y = ev_y \circ \Res_Y \circ \upsilon$.
Since $\upsilon$ is a ring isomorphism, $\ev_{f(x)} = \ev_y \circ \Res_Y$.  Hence evaluation at $f(x)$ coincides with evaluation at $y$, so $y = f(x)$.  In other words the continuous function $f\colon X \rightarrow \upsilon Y$ has image in $Y$ (under our hypothesis that $\phi$ is continuous).    

For any continuous ring homomorphism $\phi\colon \Cont(Y) \rightarrow \Cont(X)$, we find that there exists a unique continuous map $f\colon X \rightarrow Y$ such that $\phi = \Cont(f)$.  Therefore the functor from $\Cat{Top}_{lc}$ to $\C{Vec}$ is fully faithful.
\qed

\begin{remark}
The full faithfulness of $\Cont$ also follows from the spectral theorems in \cite{Mal}.
\end{remark}

\begin{corollary}
\label{ImageMap}
Let $f\colon X \rightarrow Y$ be a continuous map of locally compact Hausdorff spaces.  Let $U$ be an open subset of $Y$.  Then $\Im(f) \subset U$ if and only if the ring homomorphism $\Cont(f)\colon \Cont(Y) \rightarrow \Cont(X)$ factors through $\Cont(U)$ (via the natural restriction $\Cont(Y) \rightarrow \Cont(U)$).
\end{corollary}
\proof
One direction is clear; if $\Im(f) \subset U$, then $\Cont(f)$ factors through $\Cont(U)$.  

Let $\Inc_U\colon U \rightarrow Y$ denote the inclusion map, so $\Cont(\Inc_U)$ is the restriction map from $\Cont(Y)$ to $\Cont(U)$.  Suppose that there exists a continuous homomorphism of complete algebras $\phi\colon \Cont(U) \rightarrow \Cont(X)$ such that $\Cont(f) = \phi \circ \Cont(\Inc_U)$.  Then, since $U$ is open in $Y$, $U$ is locally compact, and hence $\phi = \Cont(\tilde f)$ for some continuous map $\tilde f\colon X \rightarrow U$.  Since $\Cont(f) = \Cont(\tilde f) \circ \Cont(\Inc_U)$, the faithfulness of $\Cont$ implies that $f = \Inc_U \circ \tilde f$, so $\Im(f) \subset U$.
\qed

\subsection{Algebraic varieties and complete algebras}

Let $\alg{X}$ be a {\em variety} over $\CC$, by which we mean a reduced and separated {\em affine} scheme over finite type over $\CC$.  We write $\OO = \OO(\alg{X})$ for the coordinate ring of $\alg{X}$.  Thus $\OO$ is a finitely-generated $\CC$-algebra, and  hence $\OO$ has countable dimension as a $\CC$-vector space.  We endow $\OO$ with the inductive limit topology, as a countable direct sum of one-dimensional spaces.  It follows that $\OO$ is a complete locally convex topological vector space; $\OO$ is a (DF)-space.  

\begin{proposition}
Let $\alg{X}$ and $\alg{Y}$ be two varieties over $\CC$, and $\alg{X} \times \alg{Y}$ their fibre product over $\CC$.  Then 
$$\OO(\alg{X}) \hat \otimes \OO(\alg{Y}) = \OO(\alg{X}) \otimes \OO(\alg{Y}) \isom \OO(\alg{X} \times \alg{Y}),$$
where the final isomorphism is given by $\odot$ as before:
$$[f \odot g](x,y) = f(x) \cdot g(y),$$
for every $f \in \OO(\alg{X})$ and every $g \in \OO(\alg{Y})$.
\end{proposition}
\proof
Choose bases $\{ P_i \}_{i = 1}^\infty$ and $\{ Q_j \}_{j=1}^\infty$ of $\OO(\alg{X})$ and $\OO(\alg{Y})$, respectively.  By \cite[Ch.\ 8, \S 44.5]{Kot}, since $\OO(\alg{X})$ is a (DF)-space, we find isomorphisms
$$\OO(\alg{X}) \hat \otimes \OO(\alg{Y}) \isom \bigoplus_{j=1}^\infty \OO(\alg{X}) \hat \otimes \CC Q_j \isom \bigoplus_{i,j} \CC (P_i \otimes Q_j)
.$$
These are the isomorphisms desired in the proposition.
\qed

As in the case of locally compact Hausdorff topological spaces, we find that $\OO(\alg{X})$ is a complete algebra whenever $\alg{X}$ is a variety over $\CC$.  Note that, when $\alg{X} \rightarrow \alg{Y}$ is a morphism of varieties $\CC$, the induced map of complete algebras $\OO(\alg{Y}) \rightarrow \OO(\alg{X})$ is continuous.  In fact, {\em every} $\CC$-algebra morphism from $\OO(\alg{Y})$ to $\OO(\alg{X})$ is continuous.  This implies the following algebraic analogue of Proposition \ref{FFF}.
\begin{proposition}
$\OO$ defines a fully faithful contravariant functor from the category $\Cat{Var}$ of varieties over $\CC$ to the category $\C{Alg}$ of complete algebras.
\end{proposition}

In this way, the topology on $\OO(\alg{X})$ does not tell us anything new about affine varieties over $\CC$.  Rather, it demonstrates that the rings of algebraic geometry -- the rings $\OO(\alg{X})$ -- can be treated on equal footing with the rings of continuous functions on locally compact Hausdorff spaces.  Importantly for what follows, we may mix algebra and topology to consider continuous maps from locally compact topological spaces to (the complex points of) algebraic varieties.

\begin{proposition}
\label{App}
Let $X$ be a locally compact Hausdorff topological space.  Let $\alg{Y}$ be a variety over $\CC$.  Let $Y = \alg{Y}(\CC)$ be the set of complex points of $\alg{Y}$, endowed with its natural topology via any closed embedding into $\CC^n$.  Let $\rho\colon  \OO(\alg{Y}) \hookrightarrow \Cont(Y)$ be the inclusion of the space of polynomial functions into the space of continuous functions.  Then $\rho$ is a continuous $\CC$-linear map, and the restriction map:
$$\rho^\ast\colon \Hom_{\C{Alg}}(\Cont(Y), \Cont(X)) \rightarrow \Hom_{\C{Alg}}(\OO(\alg{Y}), \Cont(X))$$
is a bijection.
\end{proposition}
\proof
For the continuity of $\rho$, observe that any finite-dimensional subspace of $\Cont(Y)$ has the usual topology of a finite-dimensional complex vector space.  Hence the restriction of $\rho$ to any finite-dimensional subspace of $\OO(\alg{Y})$ is continuous.  Since $\OO(\alg{Y})$ has the inductive limit topology from its finite-dimensional subspaces, the map $\rho$ is continuous.

To demonstrate that $\rho^\ast$ is injective, we must verify that a continuous algebra homomorphism $\phi$ from $\Cont(Y)$ to $\Cont(X)$ is uniquely determined by its restriction to polynomial functions on $Y$.  Such a continuous algebra homomorphism is uniquely determined by the functorially associated continuous map $f\colon X \rightarrow Y$; namely, if $x \in X$, then $\ev_x \circ \phi = \ev_y$ precisely when $y = f(x)$.

Now, if $x \in X$, then $\ev_x \circ \phi \vert_{\OO(\alg{Y})}$ is a $\CC$-algebra homomorphism from $\OO(\alg{Y})$ to $\CC$.  Since $\alg{Y} = \Spec(\OO(\alg{Y}))$, there exists a unique $y \in Y = \alg{Y}(\CC)$ such that
$$\ev_x \circ \phi \vert_{\OO(\alg{Y})} = \ev_y.$$
Hence the restriction of $\phi$ to $\OO(\alg{Y})$ determines the continuous map $f\colon X \rightarrow Y$, which in turn determines the algebra homomorphism $\phi$.

To demonstrate that $\rho^\ast$ is surjective, begin with a continuous map of algebras $\tau\colon \OO(\alg{Y}) \rightarrow \Cont(X)$.  As $\OO(\alg{Y})$ is a finitely-generated algebra over $\CC$, let us choose a presentation $\CC[T_1, \ldots, T_d] \twoheadrightarrow \OO(\alg{Y})$, and let $t_1, \ldots, t_d$ denote the images of $T_1, \ldots, T_d$ in $\OO(\alg{Y})$.  These functions provide a topological embedding $\vec t\colon Y \rightarrow \CC^d$.  

On the other hand, let $p_i = \tau(t_i)$ for $1 \leq i \leq d$.  Then each $p_i$ is a continuous function from $X$ to $\CC$, and we write $\vec p\colon X \rightarrow \CC^d$ for the resulting continuous map.  As $\tau$ is a ring homomorphism, we find that $\vec p$ has image in $Y$; we find a continuous function $\vec p\colon X \rightarrow Y$ in this way.  If $\phi\colon \Cont(Y) \rightarrow \Cont(X)$ is the associated continuous $\CC$-algebra homomorphism, then $\rho^\ast(\phi) = \tau$ as required.
\qed

The theory of locally compact Hausdorff spaces and affine schemes overlaps when discussing compact totally disconnected spaces.  When $X$ is a compact totally disconnected Hausdorff space (a Stone space), the complete algebra $\Cont(X)$ is a boolean algebra, and is naturally an inductive limit (over open covers) of finite-dimensional subalgebras.  For such a space $X$, we write $\alg{X} = \Spec(\Cont(X))$ for the resulting affine scheme over $\CC$.  Thus $\Cont(X) = \OO(\alg{X})$; the continuous functions on $X$ can be identified with the regular functions on $\alg{X}$.  
  
\subsection{Locally compact groups and Hopf algebras}

\begin{definition}
A complete coalgebra is a coalgebra in the symmetric monoidal category $\C{Vec}$, i.e., a triple $(\PP, \Delta, \epsilon)$, where $\PP \in \C{Vec}$, and 
$$\Delta\colon \PP \rightarrow \PP \hat \otimes \PP, \quad \epsilon\colon \PP \rightarrow \CC,$$
are morphisms in $\C{Vec}$ making the diagrams for counit and coassociativity commute.
\end{definition}
In most examples of interest to us, $\Delta$ does {\em not} arise from a continuous $\CC$-linear map from $\PP$ to $\PP \otimes \PP$ (the algebraic tensor product).

If $(\OO_1, \mu_1, \One_1)$ and $(\OO_2, \mu_2, \One_2)$ are complete algebras, there is a unique complete algebra structure on $\OO_1 \hat \otimes \OO_2$ extending the ordinary (algebraic) tensor product of algebras $\OO_1 \otimes \OO_2$.  Similarly, when $(\PP_1, \Delta_1, \epsilon_1)$ and  $(\PP_2, \Delta_2, \epsilon_2)$ are complete coalgebras, there is a natural coalgebra structure on $\PP_1 \hat \otimes \PP_2$.  Indeed, all of this makes sense for algebras and coalgebras in any symmetric monoidal category.

\begin{definition}
A complete bialgebra is a bialgebra $\BB$ in the symmetric monoidal category $\C{Vec}$, i.e., a quintuple $(\BB, \Delta, \epsilon, \mu, \One)$ for which $(\BB, \Delta, \epsilon)$ is a complete coalgebra, $(\BB,\mu,\One)$ is a complete algebra, and these structures satisfy compatibility diagrams which can be summarized in either of the following equivalent ways:
\begin{enumerate}
\item
$\mu$ and $\One$ are morphisms of complete coalgebras.
\item
$\Delta$ and $\epsilon$ are morphisms of complete algebras.
\end{enumerate}
\end{definition}

Finally, a Hopf algebra is a bialgebra with an antipode.
\begin{definition}
A complete Hopf algebra is a Hopf algebra $\HH$ in the symmetric monoidal category $\C{Vec}$, i.e., a sextuple $(\HH,\Delta,\epsilon,\mu,\One,S)$ for which $(\HH, \Delta, \epsilon, \mu, \One)$ is a complete bialgebra and $S\colon \HH \rightarrow \HH$ is a morphism in $\C{Vec}$ making the antipode diagrams commute.
\end{definition}

Uniqueness of the antipode in a Hopf algebra implies that the category $\C{Hopf}$ of complete Hopf algebras and continuous Hopf algebra homomorphisms is a full subcategory of the category $\C{Bialg}$ of complete bialgebras and continuous bialgebra homomorphisms.

The most important first example of a complete Hopf algebra is given by the following.
\begin{corollary}
Let $G$ be a locally compact Hausdorff topological group, with structure maps,
$$m\colon G \times G \rightarrow G, \quad \iota\colon G \rightarrow G, \quad e\colon \{ pt \} \rightarrow G.$$
Then the complete algebra $\Cont(G)$ becomes a complete Hopf algebra with comultiplication $\Delta = \Cont(m)$, counit $\epsilon = \Cont(e)$, and antipode $\Cont(\iota)$.

This extends to a fully faithful functor from the category $\Cat{Gp}_{lc}$ of locally compact Hausdorff groups and continuous group homomorphisms to the category $\C{Hopf}$ of complete Hopf algebras.
\end{corollary}
\proof
This is a formal consequence of Proposition \ref{FFF}.
\qed

The analogous statement holds for affine group schemes.
\begin{corollary}
Let $\alg{G}$ be a affine group scheme of finite type over $\CC$.  Then the ring $\OO(\alg{G})$ is naturally a complete Hopf algebra.
\end{corollary}

For profinite groups, one may apply Stone duality.
\begin{corollary}
Let $\Gamma$ be a profinite topological group, and $\alg{\Gamma} = \Spec(\Cont(\Gamma))$.  Then $\Cont(\Gamma) = \OO(\alg{\Gamma})$ is naturally a complete Hopf algebra.
\end{corollary}

\bibliographystyle{halpha}
\bibliography{MetaHopf}
 
\end{document}